\documentclass[a4paper]{amsart}
\makeatletter
\usepackage{fullpage}
\usepackage{geometry}
\newgeometry{margin=1.7in}
\usepackage{amsmath, amssymb, amsfonts, amstext, verbatim, amsthm, mathrsfs, mathtools}
\usepackage[mathcal]{eucal}
\usepackage{graphicx}
\usepackage[all]{xy}
\usepackage{tikz} 

\usepackage{enumerate}
\usepackage{enumitem}

\usepackage[normalem]{ulem}
\usepackage{soul}

\usepackage{marginnote}

\usepackage{xspace}
\usepackage[modulo]{lineno}
\usepackage{aliascnt} 
\usepackage[colorlinks=true,linkcolor=blue,citecolor=blue,urlcolor=blue,citebordercolor={0 0 1},urlbordercolor={0 0 1},linkbordercolor={0 0 1}]{hyperref} 
\usepackage[nameinlink]{cleveref}
\crefname{diagram}{Diagram}{Diagram}

\numberwithin{equation}{section}

\def\makeCal#1{%
\expandafter\newcommand\csname c#1\endcsname{\mathcal{#1}}}
\def\makeBB#1{%
\expandafter\newcommand\csname b#1\endcsname{\mathbb{#1}}}
\def\makeFrak#1{%
\expandafter\newcommand\csname f#1\endcsname{\mathfrak{#1}}}

\count@=0
\loop
\advance\count@ 1
\edef\y{\@Alph\count@}%
\expandafter\makeCal\y
\expandafter\makeBB\y
\expandafter\makeFrak\y
\ifnum\count@<26
\repeat

\theoremstyle{plain}
\newtheorem{thm}{Theorem}[section]
\newtheorem{cor}[thm]{Corollary}
\newtheorem{lem}[thm]{Lemma}

\newtheorem{prop}[thm]{Proposition}

\newtheorem*{thm*}{Theorem}

\newtheorem{thmx}{Theorem}

\theoremstyle{definition}
\newtheorem{rem}[thm]{Remark}
\newtheorem{remark}[thm]{Remark}

\newtheorem{defn}[thm]{Definition}

\newtheorem{ex}[thm]{Example}
\newtheorem{example}[thm]{Example}

\newtheorem{warn}[thm]{Warning}

\newcommand{\op}[1]{\!\!\mathop{\rm ~#1}\nolimits}

\newcommand{\epf}{\qed \vspace{+10pt}}

\DeclareMathOperator{\Spec}{Spec}
\DeclareMathOperator{\ST}{ST}
\DeclareMathOperator{\GL}{GL}
\DeclareMathOperator{\ev}{ev}
\DeclareMathOperator{\id}{id}
\DeclareMathOperator{\Hom}{Hom}
\DeclareMathOperator{\Aut}{Aut}
\DeclareMathOperator{\Frac}{Frac}

\DeclareMathOperator{\red}{red}

\DeclareMathOperator{\Bun}{Bun}

\newcommand{\hookarr}{\hookrightarrow}
\newcommand{\uHom}{\underline{{\rm Hom}}}
\newcommand{\uCoh}{\underline{{\rm Coh}}}

\newcommand{\ilim}{\varprojlim}

\DeclareMathOperator{\D}{D}

\DeclareMathOperator{\Filt}{Filt}

\newcommand{\uMap}{\underline{{\rm Map}}}
\DeclareMathOperator{\Grad}{Grad}

\DeclareMathOperator{\Flag}{Flag}
\DeclareMathOperator{\Mod}{Mod}
\DeclareMathOperator{\Ind}{Ind}

\DeclareMathOperator{\colim}{colim}
\DeclareMathOperator{\coker}{coker}
\DeclareMathOperator{\Tor}{Tor}

\DeclareMathOperator{\Coh}{Coh}
\DeclareMathOperator{\QCoh}{QCoh}
\DeclareMathOperator{\rank}{rank}
\DeclareMathOperator{\gr}{gr}

\newcommand{\co}{\colon}

\newcommand{\oh}{\cO}
\newcommand{\tensor}{\otimes}

\DeclareMathOperator{\End}{End}
\DeclareMathOperator{\num}{num}
\DeclareMathOperator{\Isom}{Isom}

\renewcommand{\ss}{{\rm ss}}
\newcommand{\spref}[1]{\href{http://stacks.math.columbia.edu/tag/#1}{Tag #1}}
\newcommand{\iso}{\stackrel{\sim}{\to}}
\DeclareMathOperator{\Stab}{Stab}

\DeclareMathOperator{\Lie}{Lie}
\DeclareMathOperator{\Map}{Map}
\DeclareMathOperator{\Higgs}{Higgs}
\DeclareMathOperator{\Ram}{Ram}
\DeclareMathOperator{\ideg}{ideg}
\DeclareMathOperator{\ad}{ad}
\DeclareMathOperator{\Gr}{Gr}
\DeclareMathOperator{\HN}{HN}
\DeclareMathOperator{\car}{car}

\DeclareMathOperator{\alg}{Alg}

\DeclareMathOperator{\wt}{wt}
\DeclareMathOperator{\Rees}{Rees}
\def\cg{\mathfrak{g}}
\def\cp{\mathfrak{p}}
\def\oST{\overline{\ST}}

\def\Scomplete{\textsf{S}-complete}
\def\tX{\widetilde{X}}
\def\tZ{\widetilde{Z}}
\def\p{p}
\def\kbar{\overline{k}}
\newcommand{\map}[1]{\stackrel{#1}{\longrightarrow}}
\newdir{|>}{!/4.5pt/@{|}*:(1,-.2)@^{>}*:(1,+.2)@_{>}}

\usepackage{pbox}
\definecolor{forrestgreen}{RGB}{34,139,34}

\author[J.~Alper]{Jarod Alper}
\address{University of Washington, Seattle, Washington, USA.}
\email{jarod@uw.edu}

\author[D.~Halpern-Leistner]{Daniel Halpern-Leistner}
\address{Cornell University, Ithaca, New York, USA.}
\email{daniel.hl@cornell.edu}

\author[J.~Heinloth]{Jochen Heinloth}
\address{Universität Duisburg-Essen, Essen, Germany}
\email{jochen.heinloth@uni-due.de}

\begin{document}

\title{Existence of moduli spaces for algebraic stacks}

\begin{abstract}
We provide necessary and sufficient conditions for when an algebraic stack admits a good moduli space and prove a semistable reduction theorem for points of algebraic stacks equipped with a $\Theta$-stratification. These results provide a generalization of the Keel--Mori theorem to moduli problems whose objects have positive dimensional automorphism groups and give criteria on the moduli problem to have a separated or proper good moduli space. To illustrate our method, we apply these results to construct proper moduli spaces parameterizing semistable $\mathcal{G}$-bundles on curves and moduli spaces for objects in abelian categories.
\end{abstract}

\maketitle

\tableofcontents

\section{Introduction}

The construction of moduli spaces is a recurring problem in the study of moduli problems in algebraic geometry. Given a moduli problem, described by an algebraic stack $\cX$, the ideal solution would be for $\cX$ to be representable by a scheme or an algebraic space. This is never the case when objects parameterized by $\cX$ have non-trivial automorphism groups. In this case one hopes for the existence of a universal map to an algebraic space $q \co \cX \to X$ with useful properties.

For algebraic stacks with finite automorphism groups the Keel--Mori theorem \cite{keelmori} gave a satisfactory existence result from the intrinsic perspective. It states that if $\cX$ is an algebraic stack of finite type over a noetherian base whose inertia stack is finite over $\cX$, then there is a \emph{coarse moduli space} $q \co \cX \to X$, which in addition to being a universal map to an algebraic space is bijective on geometric points.

The restriction to the case of finite automorphism groups is not necessary for the construction of moduli spaces using GIT. Furthermore in many examples, such as the moduli of vector bundles or coherent sheaves on a projective variety, one must consider objects with positive dimensional automorphism groups in order to construct moduli spaces that are proper. 

In \cite{alper-good}, the first author introduced the notion of a \emph{good moduli space} for an algebraic stack $\cX$ as an intrinsic formulation of many of the useful properties of the notion of a \emph{good quotient} \cite{seshadri_quotient}, a specific type of GIT quotient including all GIT quotients in characteristic $0$.  
By definition, a good moduli space is a map $q \co \cX \to X$ to an algebraic space such that the pushforward $q_\ast$ of quasi-coherent sheaves is exact, and such that the canonical map $\cO_X \to q_\ast(\cO_\cX)$ is an isomorphism. This simple definition leads to many useful properties, including that $q$ is universal for maps to an algebraic space, and that the fibers of $q$ classify orbit-closure equivalence classes of points in $\cX$. 

Our main result gives necessary and sufficient conditions under which an algebraic stack admits a good moduli space, and can be seen as uniting the main theorem of geometric invariant theory with the intrinsic perspective of the Keel--Mori theorem. In characteristic $0$ the main result admits the following simple formulation.
\begin{thmx} \label{T:existence} \label{T:A}
Let $\cX$ be an algebraic stack of finite presentation, with affine {stabilizers} and separated diagonal over a noetherian algebraic space $S$ of characteristic $0$. Then $\cX$ admits a separated good moduli space $X$ if and only if $\cX$ is $\Theta$-reductive (\Cref{D:theta-reductive}) and  {\textsf{S}}-complete  (\Cref{D:S-complete}).

Moreover, $X$ is proper if and only $\cX$ satisfies the existence part of the valuative criterion for properness.
\end{thmx}

{In the main text we prove a more general characteristic independent statement of this result, \Cref{T:existence-separated}, which provides conditions for the existence of both good and adequate moduli spaces over a general base. We also give a similar result, \Cref{T:existence-general}, which allows for potentially non-separated moduli spaces under the stronger hypothesis that $\cX$ has affine diagonal.}

Let us give an informal explanation of the above conditions.
The first condition, $\Theta$-reductivity, is the geometric analog of the statement that filtrations by semistable objects extend under specialization. This is formulated in terms of maps from the stack $\Theta:=[\bA^1/\bG_m]$ into $\cX$. 
The second condition, \Scomplete ness,
where the \textsf{S} stands for ``Seshadri,'' is a geometric property that is reminiscent of classical methods of establishing separatedness of moduli spaces. 
{More precisely, we introduce a geometric notion of an \emph{elementary modification} (\Cref{D:elementary_mod}) of a family over a DVR, and  {\textsf{S}}-completeness states that any two families over a DVR that are isomorphic at the generic point differ by an elementary modification.}  It turns out that {\textsf{S}}-completeness has many desirable consequences:  namely, {\textsf{S}}-completeness implies that automorphism groups of closed points are reductive and moreover, that connected components of automorphism groups specialize to closed points. The latter property is our analog for the finite inertia hypothesis appearing in the Keel-Mori theorem.

Ultimately, both  {\textsf{S}}-completeness and $\Theta$-reductivity are local criteria in the sense that each is equivalent to a filling condition for $\bG_m$-equivariant families over a suitable punctured regular 2-dimensional scheme.

Our second main theorem is an analog of Langton's semistable reduction theorem \cite{langton} for moduli of bundles, that works for a large class of algebraic stacks equipped with a notion of stability that induces a  $\Theta$-stratification, a geometric analog of the notion of Harder--Narasimhan--Shatz stratifications. As in Langton's theorem, the statement is that if a family of objects parametrized by a DVR specializes to a point that is more unstable than the generic fiber of the family, then one can modify the family along the closed point to get a family that has the same stability properties as the generic fiber. Surprisingly the existence of modifications can be obtained from the local geometry of $\Theta$-stratifications. The formal statement is the following. 

\begin{thmx}[\Cref{T:Langton-Algorithm}] \label{T:B} 
	Let $\cX$ be an algebraic stack locally of finite type with affine diagonal over a noetherian algebraic space $S$, and let $\cS \hookrightarrow \cX$ be a $\Theta$-stratum (\Cref{D:theta-stratum}). Let $R$ be a DVR with fraction field $K$ and residue field $\kappa$. Let $\xi \co \Spec(R) \to \cX$ be an $R$-point such that the generic point $\xi_K$ is not mapped to $\cS$, but the special point $\xi_{\kappa}$ is mapped to $\cS$: 
	$$\xymatrix{
		\Spec(K) \ar@{^(->}[r]\ar[d]^{\xi_K} & \Spec(R)\ar[d]^{\xi_R} & \Spec(\kappa)\ar@{_(->}[l]\ar[d]^{\xi_{\kappa}} \\
		\cX-\cS \ar@{^(->}[r]^{j}& \cX & \ar@{_(->}[l]_{\iota}\cS.}$$
Then there exists an extension $R \to R'$ of DVRs with $K \to K'=\Frac(R')$ finite and an elementary modification (\Cref{D:elementary_mod}) $\xi' \co \Spec(R') \to \cX$ of $\xi|_{R'}$ that lands in $\cX -\cS$.
\end{thmx}

We may apply the above results to the semistable locus $\cX^{\ss} \subset \cX$ defined by a class $\ell \in H^2(\cX;\bR)$ via the Hilbert--Mumford criterion (see \Cref{D:semistable-locus}).  As many properties of $\cX$ are inherited by the semistable locus, we can provide conditions on $\cX$ ensuring that the semistable locus $\cX^{\ss}$ admits a separated good moduli space and a further condition ensuring that the good moduli space is proper.  To summarize, we have:

\begin{thmx}[\Cref{C:valuative_criterion_semistable}, \Cref{P:semistable_locus}, \Cref{P:moduli_space_semistable_locus}] \label{T:C}
Let $\cX$ be an algebraic stack locally of finite type with affine diagonal over a quasi-separated and locally noetherian algebraic space $S$, and let $\ell \in H^2(\cX;\bR)$ be a class defining a semistable locus $\cX^{\ss} \subset \cX$ which is part of a well-ordered $\Theta$-stratification of $\cX$ compatible with $\ell$.\footnote{See \Cref{P:semistable_locus} for the precise compatibility condition between $\ell$ and the $\Theta$-stratification.}
Then if $\cX$ is either $\Theta$-reductive,  {\textsf{S}}-complete, or satisfies the existence part of the valuative criterion for properness, then the same is true for $\cX^{\ss}$.

In particular, if in addition $S$ has characteristic $0$, $\cX \to S$ is {\textsf{S}}-complete and $\Theta$-reductive, and $\cX^{\ss} \to S$ is quasi-compact, then there exists a good moduli space $\cX^{\ss} \to X$ such that $X$ is separated over $S$ (and proper over $S$ if $\cX \to S$ satisfies the existence part of the valuative criterion for properness).
\end{thmx}

We expect that in the semistable reduction theorem (\Cref{T:B}), weak $\Theta$-strata (that only require canonical filtrations to exist after a purely inseparable extension) should be sufficient, and these are available in greater generality in positive characteristic. Similarly, in positive characteristic, we expect that \Cref{T:C} holds with ``good moduli space" replaced with ``adequate moduli space."
The main obstruction for these generalizations is a version of the local structure theorem where the embedding of a stratum is replaced with a radicial map.

\subsubsection*{Applications}

To illustrate our results we give some applications that may be of independent interest. First, we use the semistable reduction theorem to give a proof that the Hitchin fibration for semistable $G$-Higgs bundles is proper if the characteristic of the ground field is not too small (\Cref{C:higgs_properness}). This result is of course expected, but we could not find it in the literature. 

Second we apply our existence theorem to construct some new good moduli spaces. Namely we construct proper good moduli spaces for semistable $\cG$-bundles for a Bruhat--Tits groups scheme $\cG$ over a smooth geometrically connected projective curve over a field of characteristic $0$, generalizing work of Balaji and Seshadri \cite{BalajiSeshadri} (\Cref{T:GoodModuliBunG}). We also construct proper good moduli spaces for objects in abelian categories in \Cref{T:abelian_moduli_space_1} and \Cref{T:abelian_semistable_GMS}.  As a special case, we construct proper moduli spaces of semistable complexes with respect to a Bridgeland stability condition on a smooth projective variety $X$ over a field of characteristic $0$.  Whereas in these examples the lack of a convenient global quotient description of the corresponding moduli problems seems to pose a serious obstruction to a construction using GIT, the verification of the conditions of our main theorems turns out to be surprisingly simple.

\subsubsection*{Acknowledgements}

We would like to thank Arend Bayer, Brian Conrad, Johan de Jong, Maksym Fedorchuk, Jack Hall, Alexander Polishchuk, David Rydh, Michael Thaddeus, Yukinobu Toda and Xiaowei Wang for helpful conversations related to this project. We thank our referees for their constructive comments and thoughtful corrections. We also thank Victoria Hoskins for organizing the workshop {\it New Techniques in Geometric Invariant Theory} where some of these ideas were first discussed.  The first author was partially supported by NSF grant DMS-1801976. The second author was partially supported by NSF grant DMS-1762669. The third author was partially supported by Sonderforschungsbereich/Transregio 45 of the DFG.

\section{Preliminaries}

Throughout we will fix a base $S$ that will be a quasi-separated algebraic space. The most interesting case for many readers will be when $S=\Spec(k)$ is the spectrum of a field.

As our arguments build on the one hand on local structure theorems and on the other hand on notions that came up in the study of notions of stability on algebraic stacks, we briefly recall these results in this section.

\subsection{Reminder on local structure theorems for algebraic stacks}

\begin{defn} \label{D:local-quotient-presentation} If $\cX$ is an algebraic stack and $x \in |\cX|$ is a point, a \emph{quotient presentation around $x$} is a pointed, flat morphism 	
	$f \co (\cW, w) \to (\cX,x)$ of algebraic stacks such that: 1) $\cW \cong [\Spec(A)/ \GL_N]$ for some $N$; and 2) $f$ induces an isomorphism of stabilizer groups at $w$, i.e., $\Aut_\cW(p) = \Aut_\cX(f(p))$ for any point $p \in \cX(k)$ representing $w$. An \emph{\'etale (resp. representable, separated, affine) quotient presentation around $x$} is a quotient presentation $f \co (\cW, w) \to (\cX,x)$ such that $f$ has the corresponding property.
\end{defn}

The \'etale local structure theorem of Alper, Hall and Rydh shows that algebraic stacks having enough linearly reductive closed points admit \'etale quotient presentations:

\begin{thm} \label{T:ahr} \cite[Thm.~1.1]{ahr2} Let $S$ be a quasi-separated algebraic space.  Let $\cX$ be an algebraic stack locally of finite presentation and quasi-separated over $S$, with affine stabilizers. If $x \in |\cX|$ is a point with image $s \in |S|$ such that the residue field extension $\kappa(x)/\kappa(s)$ is finite and the stabilizer of $x$ is linearly reductive, then there exists an \'etale quotient presentation $f \co (\cW, w) \to (\cX,x)$ around $x$. If $\cX$ has separated (resp. affine) diagonal, then there exists a representable (resp. affine), \'etale quotient presentation $f \co (\cW, w) \to (\cX,x)$ around $x$.
\end{thm}

\begin{remark} If $S$ is the spectrum of an algebraically closed field and $\cX$ has affine diagonal, the above theorem follows from \cite[Thm.~1.2]{ahr}.  In this case, one can arrange that there is an \'etale quotient presentation $(\cW,w) \to (\cX,x)$ with $\cW \cong [\Spec(A) / G_x]$, the quotient of an affine scheme by the stabilizer $G_x = \Aut_{\cX}(x)$ of $x$.
\end{remark}

\begin{remark} While $\GL_N$ is linearly reductive in characteristic 0, it is not linearly reductive in positive or mixed characteristic.  For the same reason, the morphism $[\Spec(A) / \GL_N] \to \Spec(A^{\GL_N})$ will only be an adequate moduli space (and not a good moduli space) in general.
\end{remark}

Our methods also apply more generally to algebraic stacks that admit a local description as quotient stacks for actions of reductive groups. We formalize this as follows.

\begin{defn} \label{D:locally-reductive}
A quasi-separated algebraic stack $\cX$ with affine stabilizers is \emph{locally reductive} if every point of $\cX$ specializes to a closed point and every closed point admits  an \'etale quotient presentation.
\end{defn}

\begin{remark} \label{R:locally-reductive}
	{If $\cX$ is locally reductive, then every closed point of $\cX$ has geometrically reductive stabilizer.  Indeed, there is an \'etale quotient presentation $([\Spec A/\GL_n],w) \to (\cX,x)$.  Since $w$ is necessarily a closed point, the stabilizer $G_w \cong G_x$ is geometrically reductive.}	
	However, a finite type non-closed point may have a non-geometrically reductive stabilizer.
	{On the other hand}, if $\cX$ is finitely presented over a quasi-separated algebraic space and closed points of $\cX$ have linearly reductive stabilizers, then $\cX$ is locally reductive by the local structure theorem (\Cref{T:ahr}). 
\end{remark}

\begin{remark} \label{rmk:local-structure-refinements}
	 If $\cX$ is locally reductive and has separated diagonal (resp. affine diagonal), then {by \cite[Proposition 12.5]{ahr2}} every closed point admits a representable (resp. affine) \'etale quotient presentation.
\end{remark}

To prove the semistable reduction theorem, we will need a relative version of the above local structure theorem where we fix a subgroup isomorphic to the multiplicative group $\bG_m$ of the stabilizer $G_x = \Aut_{\cX}(x)$, but do not assume $G_x$ to be linearly reductive. A general result of this form is the following theorem.

\begin{thm} {\cite[Thms.~1.3 and 1.6]{ahhlr}} \label{T:relativeslice}
Let $\cX$ be a quasi-separated algebraic stack with affine stabilizers that is locally of finite presentation over a quasi-separated algebraic space $S$. If $\cY \subset \cX$ is a closed substack, then any smooth (resp. \'etale) morphism $f_0\colon [\Spec(B)/\bG_m^n] \to \cY$ extends to {a} smooth (resp. \'etale) morphism $$f\colon [\Spec(A)/\bG_m^n] \to \cX,$$ i.e., we have $[\Spec(A)/\bG_m^n] \times_{\cX} \cY \cong [\Spec(B)/\bG_m^n]$. Furthermore, if $f_0$ is representable and $\cX$ has separated inertia, then one can arrange that $f$ is representable as well.
\end{thm}

\begin{remark}
{For stacks $\cX$ over an algebraically closed field $k$ with smooth stabilizer groups (a condition that is automatic in characteristic $0$), we will only apply the above result in the situation where $\cY=B_k G_x \subset \cX$ is the canonical inclusion of the residual gerbe of a point $x \in \cX(k)$ with automorphism group $G_x$, $\bG_m^r\subset G_x$ a subgroup and $[\Spec(k)/\bG_m^r] \to \cY$ the induced morphism. In this case the above result is a special case of \cite[Thm.~1.1]{ahr}.}
\end{remark}

\subsection{Good and adequate moduli spaces}
Throughout this paper, we use the concepts of good moduli spaces \cite[Defn.\ 4.1]{alper-good} and adequate moduli spaces \cite[Defn.\ 5.1.1]{alper-adequate}: a quasi-compact and quasi-separated morphism $\pi \co \cX \to X$ of algebraic stacks, where $X$ is an algebraic space, is a {\it good moduli space} (resp.\ an {\it adequate moduli space}) if $\pi$ is cohomologically affine (resp.\ adequately affine) and $\oh_X \to \pi_* \oh_{\cX}$ is an isomorphism.  Recall that $\pi$ is \emph{cohomologically affine} (resp. \emph{adequately affine}) if $\pi_*$ is exact on the category of quasi-coherent $\oh_{\cX}$-modules (resp.\ if for every surjection $\cA \to \cB$ of quasi-coherent $\oh_{\cX}$-algebras, then any section $s$ of $\pi_*(\cB)$ over a smooth morphism $\Spec A \to \cY$ has a positive power that lifts to a section of $\pi_*(\cA)$).

When a locally reductive stack $\cX$ admits a good or adequate moduli space $X$, there are quotient presentations \'etale-local{ly} on $X$.

\begin{prop} \label{P:local-on-ams}
	Let $\cX$ be a locally reductive algebraic stack with separated diagonal. Assume that there is an adequate moduli space $\pi \co \cX \to X$.
	For any closed point $x \in |\cX|$, there exists an \'etale neighborhood $(\Spec B, x') \to (X,\pi(x))$ and a cartesian diagram
	\[
	\xymatrix{
		[\Spec A / \GL_N]\ar[r] \ar[d]^{\pi'}		
			& \cX \ar[d]^{\pi} \\
		\Spec B \ar[r]							& X.
	} 
	\]
	where $\pi'$ is an adequate moduli space (i.e., $B = A^{\GL_N}$). 
\end{prop}

\begin{proof}
The statement follows from the existence of a representable \'etale quotient presentation around $x$ and Luna's fundamental lemma \cite[Thm. 3.14]{ahr2}.
\end{proof}

\subsection{Reminder on mapping stacks and filtrations}

As in \cite{hlinstability} we will denote by $\Theta:=[\bA^1/\bG_m]$ the quotient stack defined by the standard contracting action of the multiplicative group on the affine line and by $B\bG_m=[\Spec(\bZ)/\bG_m]$, the classifying stack of the group $\bG_m$. Both stacks are defined over $\Spec(\bZ)$ and therefore pull back to any base $S$.
Note that since $\bG_m$ is a linearly reductive group {scheme over $\Spec \bZ$}, the structure morphisms $\Theta \to \Spec(\bZ)$ and $B\bG_m \to \Spec(\bZ)$ are good moduli spaces.

Maps from $\Theta$  into a stack are the key ingredient to define stability notions on algebraic stacks \cite{hlinstability,heinloth-stability}, and we need to recall some of their properties.

By definition for any stack $\cX$ and point $\Spec(k) \to S$ a map $B\bG_{m,k} \to \cX$ is a point $x\in \cX(k)$ together with a cocharacter $\bG_{m,k} \to \Aut_{\cX}(x)$. As the action of $\bG_m$ on a vector space is the same as a grading on the vector space, we often think of a morphism $B\bG_{m}\to \cX$ as a point of $\cX$ equipped with a grading.

Similarly, a vector bundle on $\Theta=[\bA^1/\bG_m]$ is the same as a $\bG_m$ equivariant bundle on $\bA^1$ and these are the same as vector spaces equipped with a filtration. So we think of morphisms $f\co \Theta_k \to \cX$ as an object of $x_1\in \cX(k)$ (the object $f(1)$) together with a filtration of $x_1$ and as $f(0)=x_0$ as the associated graded object. 

In examples it is often easy to see that once one has found that some moduli problem is described by an algebraic stack, the stacks of filtered or graded objects are again algebraic. This turns out to be a general phenomenon, which we recall next. 
For algebraic stacks $\cX$ and $\cY$ over $S$, we denote by 
$$\uMap_S(\cY, \cX)$$ 
the stack over $S$ parameterizing $S$-morphisms $\cY \to \cX$.   If $\cY$ is defined over $\Spec(\bZ)$, we will use the convention that $\uMap_S(\cY, \cX)$ denotes the mapping stack $\uMap_S(\cY \times S, \cX)$.

That these mapping stacks are again algebraic if $\cY=\Theta$ or $\cY=B\bG_m$  for quite general $\cX$ follows from general results established in \cite[Thm.~14.9]{ahr2}, \cite[Prop.~1.1.2]{hlinstability} and \cite[Thm.~1.6]{hlp}:  if $\cX$ is locally of finite presentation and quasi-separated over a quasi-separated algebraic space $S$, with affine stabilizers, and $\cY$ is of finite presentation and with affine diagonal over $S$ such that $\cY \to S$ is flat and a good moduli space, then $\uMap_S(\cY, \cX)$ is an algebraic stack locally of finite presentation and quasi-separated over $S$, with affine stabilizers.  Moreover, if $\cX \to S$ has affine (resp.\ quasi-affine, resp.\ separated) diagonal, then so does $\uMap_S(\cY,\cX)$.

The \emph{stack of graded points in $\cX$} is defined to be $\uMap_S(B\bG_m,\cX)$, and it is denoted $\Grad(\cX)$. This terminology highlights that a morphism $B\bG_m \to \cX$ is the same as an object of $\cX$ together with a morphism of $\bG_m$ to its automorphism group, which in practice often amounts to the choice of a grading of the object. The \emph{stack of filtrations in $\cX$} is defined to be $\uMap_S(\Theta,\cX)$, and it is denoted $\Filt(\cX)$ in \cite{hlinstability}.

\subsection{The example of quotient stacks}\label{S:drinfeld} To compute examples we recall that   stacks of filtrations and graded objects have a concrete description for quotient stacks.
If $\cX=[X/G]$ is a quotient stack locally of finite type over a field $k$, where $G$ is a smooth {affine} algebraic group acting on a quasi-separated algebraic space $X$, these mapping stacks have a classical interpretation \cite[Thm.~1.4.8]{hlinstability}. To state this recall that given a cocharacter $\lambda \co  \bG_m \to G$, one defines
$$
L_{\lambda} =  \{l \in G \, | \,  l = \lambda(t) l \lambda(t)^{-1} \, \, \forall t\} \quad  \text{and} \quad
P^+_{\lambda} =  \{p \in G \, | \lim_{t \to 0} \lambda(t) p \lambda(t)^{-1} \text{ exists} \}. 
$$
If $G$ is geometrically reductive, then $P^+_{\lambda} \subset G$ is a parabolic subgroup.
There is a surjective homomorphism 
$P^+_{\lambda} \to L_{\lambda},$ defined by $p \mapsto \lim_{t \to 0} \lambda(t) p \lambda(t)^{-1}$.  

Similarly, one defines the functors:
$$\begin{aligned}
	X^0_{\lambda}		& := \uMap_k^{\bG_m}(\Spec(k), X) & \quad & \text{(the fixed locus)}\\
	X^+_{\lambda}  		& := \uMap_k^{\bG_m}(\bA^1, X) & \quad & \text{(the attractor)}
\end{aligned}$$
By  \cite[Thm.~1.4.2]{drinfeld}, 
these functors are representable by algebraic spaces.  Moreover, there are the following natural morphisms:  
a closed immersion $X^0_{\lambda} \hookarr X$,  
an unramified morphism $X^+_{\lambda} \to X$ (given by evaluation at $1$)  and  
an affine \cite[Thm.~2.22]{ahr} morphism $X^+_{\lambda} \to X^0_{\lambda}$ (given by evaluation at $0$). If $X$ is separated, then $X^+_{\lambda} \to X$ is a monomorphism.

The $k$-points of $X^0_{\lambda}$ are simply the $\lambda$-fixed points, and if $X$ is separated, the $k$-points of $X^+_{\lambda}$ are the points $x \in X(k)$ such that $\lim_{t \to 0} \lambda(t) \cdot x$ exists. The algebraic space $X^0_{\lambda}$ inherits an action of $L_{\lambda}$ and $X^+_{\lambda}$ inherits an action of $P^+_{\lambda}$ such that the evaluation map $X^+_{\lambda} \to X^0_{\lambda}$ taking $x \mapsto \lim_{t\to 0}\lambda(t) \cdot x$ is equivariant with respect to the surjection $P^+_{\lambda} \to L_{\lambda}$.

{To reinterpret this in terms of mappings stacks for quotient stacks, recall that a map $B\bG_m \to BG$ corresponds to a conjugacy class of a cocharacter $\lambda$ and that a 2-automorphism of this map corresponds to an element in $L_{\lambda}$.  More generally, a map $B\bG_m \to [X/G]$ is a conjugacy class of a cocharacter $\lambda$ and a point $x \in X^0_{\lambda}$ fixed by $\lambda$, and a 2-automorphism is an element in $L_{\lambda}$ stabilizing $x$.  Similarly, a map $[\bA^1/\bG_m] \to [X/G]$ is the data of a conjugacy class of a cocharacter $\lambda$ and a point $x \in X^+_{\lambda}$ (c.f. \cite[Lem.~1.6]{heinloth-stability}), and a 2-automorphism is an element in $P_{\lambda}$ stabilizing $x$ . This can be further upgraded to a description of the mapping stacks for quotient stacks:}
\begin{prop} \cite[Thm.~1.4.8]{hlinstability}  \label{P:theta-description}
	Let $X$ be a quasi-separated algebraic space locally of finite type over a field $k$ equipped with an action of a smooth {affine} algebraic group $G$ over $k$ with a split maximal torus. Let $\Lambda$ be a complete set of conjugacy classes of cocharacters $\bG_m \to G$. Then there are isomorphisms
	$$\begin{aligned} \,
	\Grad([X/G]):=\uMap_k(B\bG_m, [X/G])		&\cong \bigsqcup_{\lambda \in \Lambda} [X^0_{\lambda}/L_{\lambda}]; \\
	\Filt([X/G]):=\uMap_k(\Theta, [X/G] ) 	   	&\cong  \bigsqcup_{\lambda \in \Lambda} [X^+_{\lambda}/P^+_{\lambda}].
	\end{aligned}$$  
	Moreover, the morphism $\ev_1 \co \uMap_k(\Theta, [X/G]) \to [X/G]$ is induced by the $(P^+_{\lambda} \to G)$-equivariant morphism $X^+_{\lambda} \to X$.  The morphism $\gr \co \uMap_k(\Theta,[X/G]) \to \uMap(B\bG_m, [X/G])$  is induced by the $(P^+_{\lambda} \to L_{\lambda})$-equivariant morphism $X^+_{\lambda} \to X^0_{\lambda}$.
\end{prop}

The notation $\Filt(\cX)$ is chosen to indicate that we like to think of morphisms $\Theta \to \cX$ as filtered objects that admit an associated graded object. In the above description of $\Filt([X/G])$ this filtration appears in the form of the group $P_\lambda$.


\section{Valuative criteria for stacks}

In this section we introduce and study the valuative criteria appearing in the main theorem.  Two of these -- $\Theta$-reductive morphisms (\Cref{D:theta-reductive}) and {\textsf{S}}-complete morphisms (\Cref{D:S-complete}) -- concern the degeneration properties of filtered objects, i.e., maps from $\Theta$ into a stack $\cX$. The second of these criteria already includes a geometric notion of elementary modifications for $R$-valued points of a stack that generalizes the notion of elementary modifications of vector bundles. 
	
The third criterion -- unpunctured inertia -- concerns the specialization properties of connected components of automorphism groups of objects. This is the key property needed to glue good moduli spaces.


\subsection{Morphisms of stacks of filtrations}
It will be important to understand the behavior of the stacks $\uMap(\Theta,\cX)$ under morphisms $\cX \to \cY$, i.e., study the behavior of filtrations on objects under morphisms.

\begin{lem} \label{L:theta-morphism} Let $S$ be a quasi-separated algebraic space.
	Let $f \co \cX \to \cY$ be a morphism of algebraic stacks, locally of finite presentation and quasi-separated over $S$, with affine stabilizers.  If $f$ satisfies one of the following properties
	\begin{enumerate}[label=(\alph*)]
		\item representable; 
		\label{L:theta-morphism-representable}
		\item monomorphism; 
		\label{L:theta-morphism-monomorphism}
		\item separated; 
		\label{L:theta-morphism-separated}
		\item unramified; 
		\label{L:theta-morphism-unramifed}
		\item \'etale;
		\label{L:theta-morphism-etale}
		\item  \'etale, surjective and representable; or
		\label{L:theta-morphism-representable-etale-surjective}
		\item {open immersion};
		\label{L:theta-morphism-open}
	\end{enumerate}
	then  $\uMap_S(\Theta,\cX) \to \uMap_S(\Theta,\cY)$ has the same property.
\end{lem}	 

\begin{proof}
	Properties \ref{L:theta-morphism-representable} and \ref{L:theta-morphism-monomorphism} are clear.  
	For \ref{L:theta-morphism-separated}, we verify the valuative criterion by showing that if $R$ is a valuation ring with fraction field $K$ and $h_1,h_2 \co \Theta_R \to \cX$ are morphisms agreeing over $\cY$, then an isomorphism $h_1|_{\Theta_K} \to h_2|_{\Theta_K}$ extends {uniquely} to an isomorphism $h_1 \iso h_2$.  The restriction of the map $(h_1, h_2) \co \Theta_R \to \cX \times_{\cY} \cX$ to the complement $j \co \Theta_R \smallsetminus 0 \hookrightarrow \Theta_R$ of the unique closed point lifts to $\cX$ under the diagonal $\cX \to \cX \times_{\cY} \cX$.  Since the diagonal is finite {thus in particular affine} and $j_* \oh_{\Theta_R \smallsetminus 0} = \oh_{\Theta_R}$, {the fact that rational maps between affine schemes automatically extend over codimension $2$ subsets implies } that $(h_1, h_2) \co \Theta_R \to \cX \times_{\cY} \cX$ lifts to $\cX$ ({the argument for affine morphisms reappears in more detail in the proof of \Cref{P:theta-reductive}\eqref{P:theta-reductive-affine}}).
	
	Properties \ref{L:theta-morphism-unramifed} and \ref{L:theta-morphism-etale} follow from the formal lifting criterion and descent.
	For \ref{L:theta-morphism-representable-etale-surjective}, it remains to show that $\uMap_S(\Theta,\cX) \to \uMap_S(\Theta,\cY)$ is surjective.  Let $h \co \Theta_k \to \cY_s$ be a morphism over a geometric point $s\co \Spec(k) \to S$. We will use Tannaka duality to construct a lift to $\cX$.  As any \'etale representable cover of $B\bG_{m,k}$ admits a section, we may choose a lift $B \bG_{m,k} \to \cX_s$ of $B \bG_{m,k} \hookarr \Theta_k \xrightarrow{h} \cY_s$.  
	Let $\Theta^{[n]}_k = [\Spec(k[x]/x^{n+1}) / \bG_m]$ be the $n$th nilpotent thickening of $B\bG_m \hookarr \Theta$.   Since $f$ is \'etale, there exist compatible lifts $\Theta^{[n]}_k \to \cX_s$ of $\Theta^{[n]}_k \hookarr \Theta_k \xrightarrow{h} \cY_s$. Since $\Theta_k$ is coherently complete along $B\bG_{m,k}$, by \cite[Cor.~3.6]{ahr}, there is an equivalence of categories $\Map_k(\Theta_k, \cX_s) = \ilim_n \Map_k(\Theta^{[n]}_k, \cX_s)$.  This constructs the desired lift $\Theta_k \to \cX_s$ of $h$.  See also \cite[Lem.~3.2.10]{hlinstability}. {Property \ref{L:theta-morphism-open} is clear, as for an open immersion $\cX\to \cY$ a morphism from $\Theta$ to $\cY$ factors through $\cX$ if and only if the image of the closed point lies in the open substack $\cX$.}
\end{proof}

Property \ref{L:theta-morphism-representable-etale-surjective} is not preserved if the representability hypothesis is dropped.  For instance, if $\cX = B\bG_m \to B\bG_m = \cY$ is induced by $\bG_m \to \bG_m, t \to t^d$ for $d > 1$, then $\uMap_S(\Theta,\cX) \to \uMap_S(\Theta,\cY)$ is not surjective. 

The concrete description of $\uMap_k(\Theta,[X/G])$ for quotient stacks often allows to reduce the study of the local geometry of these spaces to the case of torus quotients. This is summarized in the following lemma, whose proof relies on the local structure theorem (\Cref{T:relativeslice}).

\begin{lem}{\cite[Lem.~4.4.6]{hlinstability}}  \label{L:filtration_atlas}
{Let $\cX$ be an algebraic stack with affine stabilizers that is finitely presented over a noetherian algebraic space $S$.} Then there is an affine scheme $X$ with a $\bG_m^n$ action for some $n\geq 0$, and a smooth, surjective and representable morphism $[X/\bG^n_m] \to \cX$ such that the morphism $\uMap_S(\Theta,[X/\bG^n_m]) \to \uMap_S(\Theta,\cX)$ is smooth, surjective and representable.
\end{lem}

\subsection{Property \texorpdfstring{$\Theta$-$\cP$}{Theta-P}} 
If $f \co \cX \to \cY$ is a morphism of algebraic stacks over an algebraic space $S$, we denote by $\ev(f)_1$ the induced morphism of stacks
$$
\ev(f)_1 \co \uMap_S(\Theta,\cX)  \to \cX \times_{\cY, \ev_1} \uMap_S(\Theta,\cY), \qquad \lambda \mapsto (\ev_1(\lambda), f \circ \lambda),
$$
i.e., this morphism takes an object together with a filtration in $\cX$ and remembers the object together with the induced filtration of the image in $\cY$.
It sits in a commutative diagram:
\begin{equation} 
\label[diagram]{D:relative}
\begin{split}
\xymatrix{
\uMap_S(\Theta,\cX)   \ar@/^1pc/[rrd]^-{f\circ -}  \ar@/_1pc/[rdd]_{\ev_1}  \ar[rd]^{\ev(f)_1}  \\
			& \cX \times_{\cY} \uMap_S(\Theta,\cY) \ar[r]^{p_2} \ar[d]^{p_1}		& \uMap_S(\Theta,\cY)  \ar[d]^{\ev_1}\\
			& \cX \ar[r]^f								& \cY
}
\end{split}
\end{equation}

\begin{defn}  \label{D:theta}
Let $\cP$ be a property of morphisms of algebraic stacks.  We say that a morphism $f \co \cX \to \cY$ of algebraic stacks, locally of finite presentation and quasi-separated over a quasi-separated algebraic space $S$, with affine stabilizers, 
has \emph{property $\Theta$-$\cP$} if $\ev(f)_1 \co \uMap_S(\Theta,\cX)  \to \cX \times_{\cY, \ev_1} \uMap_S(\Theta,\cY)$ has property $\cP$.   We say that $\cX$ 
 has \emph{property $\Theta$-$\cP$} if  $\cX \to \Spec(\bZ)$ does.
\end{defn}

For example, a morphism $f\co \cX \to \cY$ is $\Theta$-surjective if one can lift filtrations of any point $f(x)$ to filtrations of $x$.

The assignment $f \mapsto \ev(f)_1$ behaves well with respect to compositions and base change.  Namely, given a composition $g \circ f \co \cX \xrightarrow{f} \cY \xrightarrow{g} \cZ$ of morphisms of algebraic stacks over $S$, then $\ev(g \circ f)_1$ is naturally isomorphic to the composition
$$\begin{aligned}
\uMap_S(\Theta,\cX) &\xrightarrow{\ev(f)_1} \cX \times_{\cY} \uMap_S(\Theta,\cY) \\
&\xrightarrow{\id \times \ev(g)_1} \cX \times_{\cY} (\cY \times_{\cZ} \uMap_S(\Theta, \cZ)) \cong \cX \times_{\cZ} \uMap_S(\Theta, \cZ),
\end{aligned}$$
and if
$$\xymatrix{
\cX' \ar[r]^{f'} \ar[d]		& \cY' \ar[d] \\
\cX \ar[r]^f				& \cY
}$$
is a Cartesian diagram of algebraic stacks over $S$, then
\begin{equation}
\label[diagram]{D:theta-base-change}
\begin{split}
\xymatrix{
\uMap_S(\Theta,\cX') \ar[r]^{\ev(f')_1 \quad} \ar[d]	& \cX' \times_{\cY'} \uMap_S(\Theta,\cY') \ar[d]  \ar[r]& \uMap_S(\Theta,\cY') \ar[d]\\
\uMap_S(\Theta,\cX) \ar[r]^{\ev(f)_1 \quad}					& \cX \times_{\cY} \uMap_S(\Theta,\cY) \ar[r] & \uMap_S(\Theta,\cY)
}
\end{split}
\end{equation}
is Cartesian. We conclude:

\begin{prop} \label{P:composition-basechange}
Let $\cP$ be a property of morphisms of algebraic stacks. If $\cP$ is stable under composition and base change, then so is the property $\Theta$-$\cP$.  If $\cP$ is stable under fppf (resp.\ smooth, resp.\ \'etale) descent, then $\Theta$-$\cP$ is stable under descent by morphisms $\cY' \to \cY$ such that $\uMap_S(\Theta,\cY') \to \uMap_S(\Theta,\cY)$ is fppf (resp.\ smooth and surjective, resp.\ \'etale and surjective). \epf
\end{prop}

\begin{lem} \label{L:descent}
Let $\cP$ be a property of representable morphisms of algebraic stacks.  
If $\cP$ is stable under \'etale descent, then $\Theta$-$\cP$ is stable under descent by representable, \'etale and surjective morphisms. 
\end{lem}

\begin{proof} This follows immediately from \Cref{P:composition-basechange} and \Cref{L:theta-morphism}\ref{L:theta-morphism-representable-etale-surjective}.
\end{proof}

\begin{lem} \label{L:quasi-finite-factoring-theta-reductive}
{Let $S$ be a quasi-separated algebraic space. Let $\cX$ be an algebraic stack locally of finite presentation and quasi-separated over $S$} with quasi-finite inertia.
If $T$ is an algebraic space {over $S$}, then any morphism $\Theta_T \to \cX$ factors uniquely through $\Theta_T \to T$.
\end{lem}

\begin{proof} This follows from \cite[Lem.~1.3.13]{hlinstability}. The key point of this is to observe that $\bG_m$, that is the automorphism group of the closed point of $\Theta$, cannot map non-trivially to the finite automorphism groups of objects in $\cX$. 
	{Since the question is \'etale-local on $T$, we can assume that $T$ is an affine scheme.  Using that $\cX$ is locally of finite presentation over $S$, limit methods further reduce us to when $T$ is of finite type over $\bZ$. }
	It then follows from Tannaka duality (\cite[Theorem 1.1]{hallj_dary_coherent_tannakian_duality} or \cite{BhattHL}) that the morphism has to factor through the good moduli space $\Theta_T\to T$.
\end{proof}

\begin{lem} \label{L:basic}
Let $S$ be a quasi-separated algebraic space.
Let $f \co \cX \to \cY$ be a morphism of algebraic stacks, locally of finite presentation {and quasi-separated} over $S$, with separated diagonals and with affine stabilizers. 
\begin{enumerate} 
	\item The morphism $\ev(f)_1$ is representable. 
		 \label{L:basic-representable}
	\item If $f$ is separated, then so is $\ev(f)_1$.
		\label{L:basic-separated}
	\item  If $f$ is representable and separated, then $\ev(f)_1$ is a monomorphism.
		\label{L:basic-monomorphism}
	\item If {the relative inertia $I_{\cX/\cY}$ is quasi-finite over $\cX$}, then $\ev(f)_1$ is an isomorphism.
		\label{L:basic-isomorphism}
	\item If $f$ is \'etale, then so is $\ev(f)_1$.
		\label{L:basic-etale}	
	\item If $f$ is representable, \'etale, and separated, then $\ev(f)_1$ is an open immersion.
		\label{L:basic-open-immersion}
\end{enumerate}
\end{lem}
\begin{proof}
For \eqref{L:basic-representable}, by diagram \eqref{D:relative}, it suffices to show that $\ev_1 \co \uMap_S(\Theta,\cX) \to \cX$ is representable, which is \cite[Prop.~1.1.13]{hlinstability}. {This is not hard to see, as we only need to show that the induced morphism on automorphism groups of objects is injective, but an automorphism of an object over $\theta$ is an equivariant section of the relative automorphism group of the corresponding family over $\bA^1$ and as we assumed that the diagonal of $\cX$ is separated, any such section is uniquely determined by its restriction to $\bA^1\smallsetminus \{0\}$.} 

Part \eqref{L:basic-separated} follows from  \Cref{L:theta-morphism}\ref{L:theta-morphism-separated}.

For \eqref{L:basic-monomorphism}, to show that $\ev(f)_1$ is a monomorphism, we have to show that for an affine scheme $\Spec(R)$ over $S$ any two morphisms $h_1,h_2\colon \Theta_R \to \cX$ that become isomorphic after applying $\ev_1(f)$ have to be isomorphic, i.e., we need to show that any commutative diagram of solid arrows
$$\xymatrix{
\Spec(R) \ar[r] \ar@{^(->}[d]			& \cX \ar[d]^{\Delta} \\
 \Theta_R \ar[r]_-{(h_1,h_2)}	\ar@{-->}[ur]				& \cX \times_{\cY} \cX
}$$
can be filled in with a dotted arrow.  As $f$ is representable and separated, the base change $\cX \times_{\cX \times_{\cY} \cX} \Theta_R \to \Theta_R$ is a closed immersion containing the dense set $\Spec(R)$; it is therefore an isomorphism.

Part \eqref{L:basic-etale} follows directly from  \Cref{L:theta-morphism}\ref{L:theta-morphism-etale} using diagram \eqref{D:relative}.  Part \eqref{L:basic-open-immersion} follows directly from Parts \eqref{L:basic-monomorphism} and \eqref{L:basic-etale} as \'etale monomorphisms are open immersions.

For \eqref{L:basic-isomorphism}, {it suffices by diagrams \eqref{D:relative} and \eqref{D:theta-base-change} and descent} to show that $\ev_1 \co \uMap_S(\Theta,\cX) \to \cX$ is an isomorphism if $\cX$ has quasi-finite inertia which follows immediately from \Cref{L:quasi-finite-factoring-theta-reductive}.
\end{proof}

\begin{rem} \label{R:quasi-compact}
The morphism $\ev(f)_1$ is not in general quasi-compact.  For an example, if $f \co B \bG_{m,k} \to \Spec(k)$, the morphism $\ev(f)_1$ is the evaluation morphism $\ev_1 \co \uMap_S(\Theta,B\bG_{m,k}) = \bigsqcup_{n \in \bZ} B \bG_{m,k} \to B \bG_{m,k}$.  
\end{rem}

\begin{rem} \label{R:basic-not-monomorphism}
If $f$ is representable but not separated, then $\ev(f)_1$ is not necessarily a monomorphism.  
\end{rem}

\subsection{\texorpdfstring{$\Theta$}{Theta}-reductive morphisms} \label{S:theta-reductive}
In this section, we study the class of $\Theta$-reductive morphisms as introduced in \cite{hlinstability}. As before, we set $\Theta := [\bA^1/ \bG_m]$  defined over $\Spec(\bZ)$. 
If $R$ is a DVR with fraction field $K$, we  set $0 \in \Theta_R := \Theta \times \Spec(R)$ to be the unique closed point.  Observe that a morphism $\Theta_R \smallsetminus 0 \to \cX$ is the data of morphisms $\Spec(R) \to \cX$ and $\Theta_K \to \cX$ together with an isomorphism of their restrictions to $\Spec(K)$.

\begin{defn} \label{D:theta-reductive}
A morphism $f \co \cX \to \cY$ of locally noetherian algebraic stacks is \emph{$\Theta$-reductive} if for every DVR $R$, 
any commutative diagram
\begin{equation}
\label[diagram]{E:theta-reductive}
\begin{split}
\xymatrix{
\Theta_R \smallsetminus 0 \ar[r] \ar[d]		& \cX  \ar[d]^f \\
\Theta_R \ar[r]	\ar@{-->}[ur]			& \cY
}
\end{split}
\end{equation}
of solid arrows can be uniquely filled in.
\end{defn}

\begin{remark} \label{R:theta-reductive-equiv}
Let $S$ be a noetherian algebraic space and $f \co \cX \to \cY$ be a morphism of algebraic stacks, locally of finite type over $S$, with quasi-compact and separated diagonal over $S$, and with affine stabilizers.
Then $f$ is $\Theta$-reductive if and only if $\ev(f)_1 \co \uMap_S(\Theta,\cX)  \to \cX \times_{\cY, \ev_1} \uMap_S(\Theta,\cY)$ satisfies the valuative criterion for properness with respect to DVRs, that is,  
for every DVR $R$ with fraction field $K$, any diagram
$$\xymatrix{
\Spec(K) \ar[r] \ar[d]					& \uMap_S(\Theta,\cX)  \ar[d]^{\ev(f)_1} \\
 \Spec(R) 	 \ar@{-->}[ur] \ar[r]				& \cX \times_{\cY, \ev_1} \uMap_S(\Theta,\cY)
}$$
of solid arrows can be uniquely filled in.  Note that the morphism $\ev(f)_1$ is always representable (\Cref{L:basic}(\ref{L:basic-representable})) and locally of finite type. However, the morphism $\ev(f)_1$ is not in general quasi-compact (see \Cref{R:quasi-compact}) and therefore $\ev(f)_1$ is not in general proper.
\end{remark}

\begin{remark} \label{R:theta-reductive-absolute}
In the context of the previous remark, when $\cY$ has quasi-finite inertia the morphism $\ev_1 : \uMap_S(\Theta,\cY) \to \cY$ is an equivalence (\Cref{L:quasi-finite-factoring-theta-reductive}), and $\ev(f)_1$ is isomorphic to $\ev_1 : \uMap_S(\Theta,\cX) \to \cX$. Therefore, $f$ is $\Theta$-reductive if and only if $\cX$ is $\Theta$-reductive in the absolute sense (i.e. $\cX \to \Spec(\bZ)$ is $\Theta$-reductive). In order to be consistent with the terminology of $\Theta$-reductivity introduced in of \cite[Def.~5.1.1]{hlinstability}, we have deviated from the ``property $\Theta$-$\cP$'' naming convention.
\end{remark}

\subsubsection{Examples illustrating $\Theta$-reductivity	}
In the following examples, we work over a field $k$. The following proposition gives a criterion using the notation from \S \ref{S:drinfeld} for when a quotient stack $[X/G]$ is $\Theta$-reductive. 
\begin{prop} \label{P:quotient-theta-complete}
	Let $\cX = [X/G]$ be a quotient stack, where $X$ is a quasi-separated algebraic space locally of finite type over 
a field $k$ and $G$ is a (smooth but not necessarily connected) split reductive algebraic group over $k$.  Then $\cX$ is $\Theta$-reductive if and only if for every cocharacter $\lambda \co \bG_m \to G$, the morphism $X^+_{\lambda} \to X$ is proper.
\end{prop}

\begin{remark} If $X$ is separated, then $X^+_{\lambda} \to X$ is proper if and only if it is a closed immersion.
\end{remark}

\begin{proof}   This follows easily from the explicit description of the mapping stack $\uMap_S(\Theta,\cX)$ in \Cref{P:theta-description}.  Indeed, there is a factorization
	$$\ev_1 \co [X^+_{\lambda} / P^+_{\lambda}] \to [X / P^+_{\lambda}] \to [X/G]$$
	and since $G$ is reductive,  each $P^+_{\lambda} \subset G$ is a parabolic subgroup.  Since the quotient $G / P^+_{\lambda}$ is projective, the morphism $[X / P^+_{\lambda}] \to [X/G]$ is proper.  Thus properness of $\ev_1$ is equivalent to properness of $X^+_{\lambda} \to X$.
\end{proof}

In order to develop some intuition for $\Theta$-reductivity, we use this result to provide some basic examples and counterexamples of $\Theta$-reductivity.   For an integer $n$, we denote by $\lambda_n \co \bG_m \to \bG_m$ the cocharacter defined by $t \mapsto t^n$; in this way, the integers $\bZ$  index the cocharacters of $\bG_m$.
\begin{example}[Affine quotients] \label{E:A2} Consider the action of $\bG_m$ on $X=\bA^2$ via $t \cdot (x,y) = (tx, t^{-1}y)$.  Then   
 $$ X^+_{\lambda_n} = \left\{ 
\begin{aligned}  
	V(y)		&	& \text{if $n > 0$} \\
	\bA^2	& 	& \text{if $n = 0$}  \\
	V(x)		& 	& \text{if $n < 0$} 
\end{aligned}\right.
$$
On the component indexed by $\lambda_n$, the evaluation morphism $[X^+_{\lambda_n} / \bG_m] \to [X/\bG_m]$ is induced by the inclusion $X^+_{\lambda_n} \to X$.  We see directly that $[X/\bG_m]$ is $\Theta$-reductive.

More generally, if $X = \Spec(A)$ is an affine scheme of finite type over $k$ with an action of a reductive algebraic group $G$, then $[X/G]$ is $\Theta$-reductive.  Indeed, if $\lambda \co \bG_m \to G$ is a cocharacter, then $A$ inherits a $\bZ$-grading $A = \bigoplus_{n \in \bZ} A_n$.  If $I^-_{\lambda}$ denotes the ideal generated by homogeneous elements of strictly negative degree, then it is easy to see that $X^+_{\lambda} = V(I^-_{\lambda})$; see \cite[\S 1.3.4]{drinfeld}.  Thus, $X^+_{\lambda} \to X$ is a closed immersion and the conclusion follows from the characterization in \Cref{P:quotient-theta-complete}. \end{example}

\begin{example} \label{E:thetaR} In contrast, quotients of schemes that are not affine are not always $\Theta$-reductive.
	Consider the action of $\bG_m$ on $X=\bA^2 \smallsetminus 0$ via $t \cdot (x,y) = (tx, y)$.  Then
 $$ X^+_{\lambda_n} = \left\{
\begin{aligned}  
	\{y \neq 0\}	&	& \text{if $n > 0$} \\
	X		& 	& \text{if $n = 0$}  \\
	V(x)			& 	& \text{if $n < 0$} 
\end{aligned}\right.
$$
and we see that $[X/\bG_m]$ is \emph{not} $\Theta$-reductive as $X^+_{\lambda_n} \to X$ is not proper for $n > 0$.    Similarly, for a DVR $R$, the algebraic stack $\Theta_R \smallsetminus 0$ is not $\Theta$-reductive.  These are the prototypical examples of non-$\Theta$-reductive stacks.

Another example is given by projective quotients. Consider the multiplication action of $\bG_m$ on $X = \bP^1$ via $t \cdot [x:y] = [tx: y]$ and on the nodal cubic $C \subset \bP^2$ such that the normalization $\bP^1 \to X$ is $\bG_m$-equivariant.
	Then 
$$ X^+_{\lambda_n} = \left\{ 
\begin{aligned}  
	\bP^1 \smallsetminus \{0\} \sqcup \{0\}			&	& \text{if $n > 0$} \\
	\bP^1								& 	& \text{if $n = 0$}  \\
	\bP^1 \smallsetminus \{\infty\} \sqcup \{\infty\}		& 	& \text{if $n < 0$} 
\end{aligned}\right.
\quad \text{and}
\quad
C^+_{\lambda_n} = \left\{ 
\begin{aligned}  
	\bP^1 \smallsetminus \{0\}			&	& \text{if $n > 0$} \\
	C									& 	& \text{if $n = 0$}  \\
	\bP^1 \smallsetminus  \{\infty\}		& 	& \text{if $n < 0$} 
\end{aligned}\right.
$$
where the maps $C^+_{\lambda_n} \to C$ for $n \neq 0$ are induced by the normalization.
We see that $ [\bP^1 / \bG_m]$ and $[C/\bG_m]$ are \emph{not} $\Theta$-reductive.

\end{example}

\subsubsection{Properties of $\Theta$-reductive morphisms}
We now give a few properties of $\Theta$-reductive morphisms.  First observe
that $\Theta$-reductive morphisms are stable under composition and base change.
We first show that one can check the lifting criterion of \eqref{E:theta-reductive} after taking extensions of the DVR, and in some situations it suffices to check only essentially finite type DVRs.

\begin{prop} \label{P:theta-reductive-extensions}
Let $ \cX \to \cY$ be a morphism of locally noetherian algebraic stacks, and consider a diagram of the form \eqref{E:theta-reductive}.  There exists a unique dotted arrow filling in the diagram if either
\begin{enumerate}
\item there exists a unique filling after passing to an unramified extension $R \subset R'$ of DVRs which is an isomorphism on residue fields, such as the completion of $R$, or
\item $\cX \to \cY$ has affine diagonal, and there exists a filling after an arbitrary extension of DVRs $R \subset R'$.
\end{enumerate}
In particular, to verify that $\cX \to \cY$ is $\Theta$-reductive, it suffices to check the lifting criterion \eqref{E:theta-reductive} for complete DVRs.  
\end{prop}

\begin{proof} The first statement follows from an explicit descent argument similar to \cite[Rmk.~2.5]{heinloth-stability}.  Alternatively, 
if $R \subset R'$ is an unramified extension of DVRs with isomorphic residue fields, then
$$\xymatrix{
\Theta_{R'} \smallsetminus 0 \ar[r] \ar[d]		& \Theta_{R} \smallsetminus 0 \ar[d] \\
\Theta_{R'} \ar[r]	& \Theta_R
}$$
is a flat Mayer--Vietoris square (\cite[Defn. 1.2]{hall-rydh-mv}) and thus by \cite[Thm. A]{hall-rydh-mv} is a pushout in the 2-category of algebraic stacks.  This establishes the first statement.

For the second statement, we begin with the observation that if $\cX \to \cY$ has affine diagonal and $j \co \cU \to \cT$ is an open immersion of algebraic stacks over $\cY$ with $j_* \oh_{\cU} = \oh_{\cT}$, then any two extensions $f_1, f_2 \co \cT \to \cX$ of a $\cY$-morphism $\cU \to \cX$ are canonically 2-isomorphic.  Indeed, since $\Isom_{\cT}(f_1, f_2) \to \cT$ is affine, the section over $\cU$ induced by the 2-isomorphism $f_1|_{\cU} \iso f_2|_{\cU}$ extends uniquely to a section of $\cT$.

Consider a diagram \eqref{E:theta-reductive}, an extension of DVRs $R \subset R'$ and a lifting $\Theta_{R'} \to \cX$.  The open immersion $j \co \Theta_R \smallsetminus 0 \to \Theta_R$ satisfies $j_* \oh_{ \Theta_R \smallsetminus 0} = \oh_{\Theta_R}$ and by flat base change, the same property holds for the morphisms obtained by base changing $j$ along  $\Theta_{R'} \to \Theta_{R}$, $\Theta_{R'} \times_{\Theta_{R}} \Theta_{R'} \to \Theta_{R}$, and $\Theta_{R'} \times_{\Theta_{R}} \Theta_{R'}  \times_{\Theta_{R}} \Theta_{R'} \to \Theta_{R}$. By the above observation, there exists a canonical 2-isomorphism between the two extensions $\Theta_{R'} \times_{\Theta_{R}} \Theta_{R'}  \rightrightarrows \Theta_{R'} \to \cX$ which necessarily satisfies the cocycle condition.  By fpqc descent, the lifting $\Theta_{R'} \to \cX$ descends to a lifting  $\Theta_{R} \to \cX$.
\end{proof}

\begin{prop} \label{P:theta-reductive-essentially-finite-type-DVRs}
{Let $S$ be a quasi-excellent and quasi-separated algebraic space.}
Let $f \co \cX \to \cY$ be a quasi-compact morphism of algebraic stacks locally of finite type and quasi-separated {over $S$ with separated diagonals} and affine stabilizers. Assume that $\cX$ admits a representable \'etale cover by quotient stacks.
		To verify that $\cX \to \cY$ is $\Theta$-reductive, it suffices to check the lifting criterion for DVRs $R$ essentially of finite type over $S$.
	\end{prop}

To prove the result we need that the irreducible components of the mapping stacks appearing in the definition of $\Theta$-reductivity are quasi-compact. As the general results only provide components locally of finite type, we restrict to stacks covered by quotients stacks in the above statement to get a more explicit hold on the components.

\begin{proof}
It suffices to prove this for base change along every quasi-compact open substack of $\cY$, and thus we may assume that $\cY$ and $\cX$ are quasi-compact. The morphism $$\ev(f)_1 \colon \uMap_S(\Theta,\cX) \to \uMap_S(\Theta,\cY) \times_{\cY} \cX$$ is representable, but not necessarily quasi-compact. If $[X/G] \to \cX$ is an \'etale, surjective, representable morphism, then so is $\uMap_S(\Theta,[X/G]) \to \uMap_S(\Theta,\cX)$ by \Cref{L:theta-morphism}(f). We can use this to show that the irreducible components of $\uMap_S(\Theta,\cX)$ are quasi-compact as follows:

By \cite[Cor.~1.1.7]{hlinstability} the canonical map is an isomorphism
\[
\uMap_S(B\bG_m,[X/G]) \cong \uMap_S(B\bG_m,\cX) \times_\cX [X/G],
\]
so only finitely many irreducible components of $\uMap_S(B\bG_m,[X/G])$ lie over each irreducible component of $\uMap_S(B\bG_m,\cX)$. On the other hand, \cite[Lem.~1.3.8]{hlinstability} that the restriction to $0$ map $\uMap_S(\Theta,\cX) \to \uMap_S(B\bG_m,\cX)$ is finite type, and the morphism $\uMap_S(\Theta,[X/G]) \to \uMap_S(\Theta,\cX)$ is compatible with these projections. It follows that only finitely many components of $\uMap_S(\Theta,[X/G])$ map to each component of $\uMap_S(\Theta,\cX)$. The explicit description of $\uMap_S(\Theta,[X/G])= \coprod [X_\lambda^+/P_{\lambda}^+]$ in \Cref{P:theta-description} (\cite[Thm.~1.4.7]{hlinstability} over a general base) then implies that each irreducible component is quasi-compact, and hence each irreducible component of $\uMap_S(\Theta,\cX)$ admits a cover by finitely many quasi-compact stacks.

Given any lifting diagram for a DVR $R$ with fraction field $K$
	$$\xymatrix{
		\Spec(K) \ar[d] \ar[r] &  \uMap_S(\Theta,\cX) \ar[d]^{\ev(f)_1} \\ \Spec(R) \ar@{-->}[ur] \ar[r] & \uMap_S(\Theta,\cY) \times_{\cY} \cX
	}$$
it suffices to verify the lifting criterion separately on each connected component of $\uMap_S(\Theta,\cX)$. The map $\Spec(R) \to \uMap_S(\Theta,\cY) \times_{\cY} \cX$ factors through an affine scheme $\Spec(A)$ of finite type over $S$ (because the stack is of finite presentation and thus limit-preserving), and it suffices to prove the lifting criterion for the base change of $\ev(f)_1$ to $\Spec(A)$. As $\ev(f)_1$ is representable this reduces the problem to a lifting problem for a quasi-compact morphism of algebraic spaces, where it has a solution if and only if it has a solution for DVRs essentially of finite type by \Cref{lem:ess_fintype_criterion}.

\end{proof}

$\Theta$-reductivity satisfies the following two descent properties.  The second property is not used in this paper and is only included for thoroughness.

\begin{prop} \label{P:theta-reductiveheta-reductive-descent} Let $\cX \to \cY$ be a morphism of locally noetherian algebraic stacks.
\begin{enumerate}
\item  If $\cY' \to \cY$ is an \'etale, representable and surjective morphism, then $\cX \to \cY$ is $\Theta$-reductive if and only if $\cX \times_{\cY} \cY' \to \cY'$ is $\Theta$-reductive.
\item If $\cX' \to \cX$ is a finite, \'etale and surjective morphism, then $\cX \to \cY$ is $\Theta$-reductive if and only if $\cX' \to \cY$ is $\Theta$-reductive.
\end{enumerate}
\end{prop}

\begin{remark} If $\cX \to \cY$ has affine diagonal, then (2) also holds, with a similar proof, if one replaces the words `finite, \'etale' with `quasi-compact, universally closed.'
\end{remark}

\begin{proof}
For (1), to check $\Theta$-reductivity of $\cX \to \cY$, by \Cref{P:theta-reductive-extensions} we may assume we have a diagram \eqref{E:theta-reductive} where $R$ is a complete DVR. As any \'etale, representable cover of $\Theta_R$ has a section after a finite \'etale extension  $R \subset R'$, we may lift the composition $\Theta_{R'} \to \Theta_R \to \cY$ to $\Theta_{R'} \to \cY'$. The $\Theta$-reductivity of $\cX' := \cX \times_{\cY} \cY' \to \cY'$ shows that the lift $\Theta_{R'} \smallsetminus 0 \to \cX'$ extends uniquely to a morphism $\Theta_{R'} \to \cX'$. This implies that the lift $\Theta_{R'} \smallsetminus 0 \to \cX$ extends uniquely to a morphism $\Theta_{R'} \to \cX$ as well, because both extension problems can be rephrased in terms of sections of $\Theta_{R'} \times_{\cY'} \cX' \simeq \Theta_{R'} \times_{\cY} \cX$.

Finally, the first few levels of the Cech nerve for the \'etale cover $\Theta_{R'} \to \Theta_R$ have the form
\[
\xymatrix{ \ldots \ar@<-1.5ex>[r] \ar@<-.5ex>[r] \ar@<.5ex>[r] \ar@<1.5ex>[r]  & \bigsqcup_j \Theta_{R'''_j} \ar@<-1ex>[r] \ar[r] \ar@<1ex>[r]  & \bigsqcup_i \Theta_{R''_i} \ar@<.5ex>[r] \ar@<-.5ex>[r]& \Theta_{R'} },
\]
for some complete DVRs $R''_i$ and $R'''_j$. The argument of the previous paragraph shows that for any of the DVRs $A=R',R''_i,R'''_j$ the lift $\Theta_A \smallsetminus 0 \to \cX$ extends uniquely to a lift $\Theta_A \to \cX$ of the map $\Theta_A \to \cY$. \'Etale descent now implies that the original lift $\Theta_{R} \smallsetminus 0 \to \cX$ extends uniquely to a morphism $\Theta_R \to \cX$.

For (2), the `only if' direction follows since finite morphisms are $\Theta$-reductive.
Conversely, given a diagram \eqref{E:theta-reductive}, we may find a finite \'etale extension $R \subset R'$ with fraction field $K \subset K'$ such that the composition $\Spec(K) \to \Theta_R \smallsetminus 0 \to \cX$ lifts to a map $\Spec(K') \to \cX'$.  As $\cX' \to \cX$ is finite, we may extend this morphism uniquely to a map $\Theta_{R'} \smallsetminus 0 \to \cX'$ lifting $\Theta_{R'} \smallsetminus 0 \to \Theta_{R} \smallsetminus 0 \to \cX$.  By $\Theta$-reductivity of $\cX' \to \cY$, this map extends uniquely to a morphism $\Theta_{R'} \to \cX'$. The composition $\Theta_{R'} \to \cX' \to \cX$ is an extension of $\Theta_{R'} \smallsetminus 0 \to \Theta_{R} \smallsetminus 0 \to \cX$, and is unique since $\cX' \to \cX$ is $\Theta$-reductive.   By an \'etale descent argument similar to the one given in Part (1), this descends uniquely to the desired lift $\Theta_{R} \to  \cX$.
\end{proof}

We now provide some important classes of $\Theta$-reductive morphisms.

\begin{prop} \label{P:theta-reductive}
 \qquad
\begin{enumerate}
	\item \label{P:theta-reductive-affine}
		An affine morphism of locally noetherian algebraic stacks is $\Theta$-reductive.
	\item \label{P:theta-reductive-quotient}
		Let $S$ be a locally noetherian scheme.  Let $G \to S$ be a geometrically reductive and \'etale-locally embeddable group scheme (e.g. reductive) acting on a locally noetherian scheme $X$ affine over $S$.  Then the morphism $[X/G] \to S$ is $\Theta$-reductive.
	\item \label{P:theta-reductive-ams} If $\cX$ is a locally reductive, locally noetherian algebraic stack with separated diagonal, then a locally noetherian\footnote{Unlike for good moduli spaces, $X$ is not guaranteed to be noetherian, and $\Theta$-reductivity is only defined for morphisms of locally noetherian stacks.} adequate moduli space $\cX \to X$ is $\Theta$-reductive.
\end{enumerate} 
\end{prop}

\begin{remark}
When $G$ is a smooth split reductive algebraic group over a field $S=\Spec(k)$, the claim in Part \eqref{P:theta-reductive-quotient} follows from the explicit calculation in \Cref{E:A2}.
\end{remark}

\begin{proof}
For \eqref{P:theta-reductive-affine}, since $0 \in \Theta_{R}$ has codimension $2$ and $\Theta_R$ is regular for a DVR $R$, we have that 
$(\Theta_R \smallsetminus 0 \to \Theta_R)_* \oh_{\Theta_R \smallsetminus 0} = \oh_{\Theta_R}$.  Given an affine morphism $f \co \cX \to \cY$,  we have canonical isomorphisms
$$\begin{aligned}
\Map_{\cY} (\Theta_{R} \smallsetminus 0, \cX)   & \cong \Map_{\oh_{\cY}-\text{alg}}(f_* \oh_{\cX}, (\Theta_R \smallsetminus 0 \to \cY)_* \oh_{\Theta_R \smallsetminus 0}) \\
	& \cong \Map_{\oh_{\cY}-\text{alg}}(f_* \oh_{\cX}, (\Theta_R \to \cY)_* \oh_{\Theta_R}) \\
	& \cong \Map_\cY(\Theta_R, \cX).
\end{aligned}$$
See also \cite[Prop.~1.3.2]{hlinstability}, which shows that $\ev(f)_1$ is a closed immersion when $f$ is affine.

For \eqref{P:theta-reductive-quotient}, since $\Theta$-reductive morphisms descend under representable, \'etale and surjective morphisms (\Cref{P:theta-reductiveheta-reductive-descent}), we may assume that $S$ is an affine noetherian scheme and that $G$ is a closed subgroup of $\GL_{N,S}$ for some $N$.  We first show that $B_{\bZ} \GL_{N} = [\Spec(\bZ) / \GL_N]$ is $\Theta$-reductive, which implies that $B_S \GL_{N} = [S/\GL_{N,S}]$ is also $\Theta$-reductive. A morphism $\Theta_R \smallsetminus 0 \to \cX$ corresponds to a vector bundle $\cE$ on $\Theta_R \smallsetminus 0$.  If $\widetilde{\cE}$ is any coherent sheaf on $\Theta_R$ extending $\cE$, then the double dual $\widetilde{\cE}^{\vee \vee}$ is a vector bundle extending $\cE$.  This provides the desired extension $\Theta_R \to \cX$.  Since $\GL_{N,S} / G$ is affine \cite[Thm.~9.4.1]{alper-adequate}, $B_S G \to B_S \GL_N$ is affine.  By Part \eqref{P:theta-reductive-affine}, $B_S G$ is $\Theta$-reductive.  Since $X$ is affine over $S$,  $[X/G] \to B_SG$ is affine which implies using again Part \eqref{P:theta-reductive-affine} that $[X/G]$ is $\Theta$-reductive.

For \eqref{P:theta-reductive-ams}, we may assume that $X$ is quasi-compact. By \Cref{P:local-on-ams}, 
there exists an \'etale cover $\Spec(B) \to X$ such that $\cX \times_X \Spec(B) \cong [\Spec(A) / \GL_N]$ for some $N$ and $B =A^{\GL_N}$.  Since $\Theta$-reductive morphisms descend under representable, \'etale and surjective morphisms, this reduces the claim to the statement that $[\Spec(A) / \GL_N] \to \Spec(A^{\GL_N})$ is $\Theta$-reductive which follows from Part \eqref{P:theta-reductive-quotient}. 
\end{proof}

\begin{prop}  \label{P:quasi-finite-theta-reductive}
	{Let $S$ be a quasi-separated algebraic space.}
	Let $f \co \cX \to \cY$ be a morphism of algebraic stacks {locally of finite presentation and quasi-separated over $S$.}  {If $\cX$ and $\cY$ have quasi-finite inertia, then $f$} is $\Theta$-reductive.
\end{prop}

\begin{proof}
This follows from \Cref{L:quasi-finite-factoring-theta-reductive}.
\end{proof}

\subsubsection{Specialization of $k$-points}
Next we provide general criteria for when specialization of $k$-points can be realized by a morphism from $\Theta_k$.

\begin{lem} \label{L:isotrivial}
Let $\cX$ be a locally reductive algebraic stack locally of finite type over {an algebraically closed field $k$}.
Then any specialization $x \rightsquigarrow x_0$ of $k$-points, where $x_0$ is a closed point, is realized by a morphism $\Theta_k \to \cX$.
\end{lem}

\begin{proof}
The condition can be checked on an \'etale neighborhood of $x_0$. As $\cX$ is locally reductive, we may assume that $\cX=[\Spec(A)/G]$ is a quotient by a reductive group. {In this case the closed point $x_0$ defines a closed $G$-orbit $Z=G\tilde{x_0}\subset \Spec(A)$ and $x$ defines a $G$-orbit $G\tilde{x}\subset \Spec(A)$ whose closure contains $Z$, thus the Hilbert--Mumford criterion \cite[Thm.~4.2]{kempf1978instability} defines a cocharacter $\lambda \colon \bG_m \to G$ such that $\lim_{t \to 0}\lambda(t)\tilde{x_0}\in Z$, i.e., $\lambda$ defines a morphism $\Theta_{k}=[\bA^1/\bG_m] \to [\Spec(A)/G]$ that realizes the specialization $x\rightsquigarrow x_0$.}
\end{proof}

The topology of $k$-points of $\Theta$-reductive stacks is analogous to the topology of quotient stacks arising from GIT.

\begin{lem} \label{L:unique_closed_point}
Let $\cX$ be a locally reductive algebraic stack locally of finite type over a field $k$. If $\cX$ is $\Theta$-reductive, then the closure of any $k$-point $p$ contains a unique closed point $x$.
\end{lem}
\begin{proof}
Assume that $x$ and $x'$ are two closed points in the closure of $p$. After replacing $k$ with an extension if necessary, we may assume that $k$ is {algebraically closed}, and that $x$ and $x'$ are $k$-rational.  It follows from \Cref{L:isotrivial} that these specializations come from two filtrations $f,f' : \Theta_k \to \cX$ with $f(1) \simeq f'(1) \simeq p$, $f(0) \simeq x$ and $f'(0) \simeq x'$. The maps $f$ and $f'$ glue to define a map $[\bA^2_k - \{(0,0)\} / (\bG_m^2)_k]$, 
and choosing one of the two $\bG_m$ factors we can apply $\Theta$-reductivity {to the local ring of $\bA^1$ at $0$} to extend this morphism to a map $[\bA^2/\bG_m]\to \cX$. Then $\gamma(0,0)$ is a specialization of both $x \simeq \gamma(1,0)$ and $x' \simeq \gamma(0,1)$, which because $x$ and $x'$ are closed implies that $x \simeq \gamma(0,0) \simeq x'$.
\end{proof}

\subsection{\texorpdfstring{$\Theta$}{Theta}-surjective morphisms} \label{S:theta-surjective}
In this section, we study the class of $\Theta$-surjective morphisms.  We will observe that $\Theta$-surjective morphisms between locally reductive algebraic stacks necessarily map closed points to closed points (\Cref{L:closed-to-closed}).  This notion will play a fundamental role in our proof of \Cref{T:existence}; namely, we will use  $\Theta$-reductivity to ensure that we can find \'etale quotient presentations which are $\Theta$-surjective (\Cref{P:quotient-presentations}\eqref{P:quotient-presentations-theta-surjective}).

By \Cref{D:theta}, a morphism $f \co \cX \to \cY$ (of algebraic stacks, locally of finite presentation over a quasi-separated algebraic space $S$, with quasi-compact and separated diagonal over $S$, and with affine stabilizers) is \emph{$\Theta$-surjective} if
$$\ev(f)_1 \co \uMap(\Theta,\cX)  \to \cX \times_{\cY, \ev_1} \uMap(\Theta,\cY)$$
is surjective.  From \Cref{P:composition-basechange} and \Cref{L:descent}, $\Theta$-surjective morphisms are stable under composition and base change, and they descend under representable, \'etale and surjective morphisms.

\begin{rem}  \label{R:theta-surjective-valuative-criterion}
The condition of $\Theta$-surjectivity translates into the following lifting criterion:  For a field $k$, denote by $j \co \Spec(k) \hookarr \Theta_k$ the open immersion.  Then $f \co \cX \to \cY$ is $\Theta$-surjective if and only if for any algebraically closed field $k$, any commutative diagram
\begin{equation} \label[diagram]{E:theta-surjective}
\begin{split}
\xymatrix{
\Spec(k) \ar[d]^{j} \ar[r]		& \cX   \ar[d]^f\\
 \Theta_k \ar[r]	\ar@{-->}[ur]		& \cY
} \end{split}
\end{equation}
of solid arrows can be filled in with a dotted arrow.
\end{rem}

\begin{rem} \label{R:theta-surjective-integral-morphism}
If $f$ is representable and separated, it follows from \Cref{L:basic}\eqref{L:basic-monomorphism} that there is at most one lift in diagram \eqref{E:theta-surjective}, that is, $f$ is $\Theta$-universally injective (or equivalently $\Theta$-radicial). This fails for non-separated morphisms. 

We also note that if $f$ is representable, separated, quasi-compact and universally closed, then the valuative criterion implies that there exists a unique lift in the above diagram.  In particular, representable proper morphisms are $\Theta$-universally bijective.
\end{rem}

If $\cX$ is an algebraic stack over a quasi-separated algebraic space $S$ and $s \in |S|$, let $\cX_s$ be the fiber product $\cX \times_S \Spec(\kappa(s))$, where $\kappa(s)$ is the residue field of $s$.

\begin{lem} \label{L:closed-to-closed} 
Let $S$ be a quasi-separated algebraic space and $f \co \cX \to \cY$ be a morphism of algebraic stacks, locally of finite presentation and quasi-separated over $S$, with affine stabilizers.  Suppose that $\cY$ is locally reductive and $f$ is $\Theta$-surjective.  If $x \in |\cX|$ is a point with image $s \in |S|$ such that $x \in |\cX_s|$ is closed, then $f(x) \in |\cY_s|$ is closed.
\end{lem}

\begin{proof}
We immediately reduce to the case when $S$ is the spectrum of an algebraically closed field $k$ and $x \in |\cX|$ is a closed point.  If $f(x)$ is not closed, then there exists a specialization $f(x) \rightsquigarrow y_0$ of $k$-points to a closed point.  By \Cref{L:isotrivial}, there exists a morphism $\Theta_k \to \cY$ realizing $f(x) \rightsquigarrow y_0$.  As the diagram
$$\xymatrix{
\Spec(k) \ar[d]^{j} \ar[r]	^{x}	& \cX   \ar[d]^f\\
 \Theta_k \ar[r]	\ar@{-->}[ur]^h		& \cY
}$$
can be filled in with a morphism $h$ and $x \in |\cX|$ is closed, $h(0) = h(1)$.  It follows that $f(x) = y_0$ is closed.
\end{proof}

\begin{rem}
The converse of  \Cref{L:closed-to-closed} is not true; see \Cref{E:relative5}.  
\end{rem}
For the construction of good moduli spaces we will need a variant of the above properties. 
Let $f \co \cX \to \cY$ be a morphism of algebraic stacks, locally of finite presentation and quasi-separated over a quasi-separated algebraic space $S$, with affine stabilizers.
Define $\Sigma_f \subset |\cX|$ to be the set of points $x \in |\cX|$ where $f$ is \emph{not} $\Theta$-surjective at $x$, i.e., points $x \in |\cX|$ where there exists a representative $\Spec(k) \to \cX$ of $x$ with $k$ algebraically closed and a commutative diagram as in diagram \eqref{E:theta-surjective} which cannot be filled in. 
By definition, $\Sigma_f$ is the image under $p_1$ of the complement of the image of $\ev(f)_1$, i.e.,
\begin{equation} \label{E:sigma-equality}
\Sigma_f = p_1\bigg( \big( \cX \times_{\cY}  \uMap_S(\Theta,\cY) \big) \, \smallsetminus \,  \ev(f)_1 (\uMap_S(\Theta,\cX) ) \bigg) \subset  |\cX|.
\end{equation}

\begin{lem} \label{L:sigma-closed}
Let $\cX$ and $\cY$ be algebraic stacks, locally of finite {type and quasi-separated over a noetherian} algebraic space $S$, with separated diagonal and affine stabilizers. Let $f \co \cX \to \cY$ be a representable, quasi-finite, and separated morphism.  Suppose that either
	\begin{enumerate}
		\item $\cY$ admits a good moduli space; or
		\item $\cY$ is locally reductive and admits an adequate moduli space.
	\end{enumerate}
	Then the locus $\Sigma_f \subset |\cX|$ is closed.
\end{lem}

\begin{proof}
	Zariski's Main Theorem \cite[Thm.~16.5]{lmb} provides a factorization $f \co \cX \xrightarrow{j} \widetilde{\cY} \xrightarrow{\nu} \cY$ where $j$ is an open immersion and $\nu$ is a {finite} morphism.  Since $\nu$ is $\Theta$-surjective by \Cref{R:theta-surjective-integral-morphism} and $\Sigma_j = \Sigma_f$, it suffices to assume that $f$ is an open immersion.  Let $\cZ \subset \cY$ be the reduced complement of $\cX$  and let $\pi \co \cY \to Y$ denote the adequate moduli space.  We claim that $\Sigma_f = \pi^{-1} (\pi(|\cZ|)) \cap |\cX|$.  This would finish the proof as $\pi^{-1}(\pi(|\cZ|)) \subset |\cY|$ is closed.
	
	Indeed, the inclusion ``$\subset$" is clear:  the morphism $\cY \smallsetminus \pi^{-1} (\pi(|\cZ|)) \hookarr \cY$ is the base change of the $\Theta$-surjective morphism $Y \smallsetminus \pi(|\cZ|) \hookarr Y$ of algebraic spaces.  For the inclusion ``$\supset$," let $x \in \pi^{-1} (\pi(|\cZ|)) \cap |\cX|$ and let $\overline{x} \co \Spec(k) \to \cX$ be a representative of $x$, where $k$ is algebraically closed, with image $s \co \Spec(k) \to S$.   Let $x_s \in |\cX_s|$ be the image of $\Spec(k) \to \cX_s$ and $z \in |\cZ_s|$ be the unique closed point in the closure of $x_s$.  If $\cY$ admits a good moduli space, it is in particular locally reductive. In either case (1) or (2), we may apply  \Cref{L:isotrivial} to obtain a morphism $\Theta_k \to \cY_s$ realizing the specialization $x_s \rightsquigarrow z$.  Since the commutative diagram
	$$\xymatrix{
		\Spec(k) \ar[r]^{\overline{x}} \ar[d]						& \cX \ar[d]^f 	\\
		\Theta_k \ar[r] \ar@{-->}[ur]					& \cY
	}$$
	does not admit a lift, $x \in \Sigma_f$. 
\end{proof}

\begin{prop} \label{P:sigma-constructible} 
Let $\cX$ and $\cY$ be algebraic stacks of finite {type and quasi-separated} over a {noetherian} algebraic space $S$, with separated diagonal and affine stabilizers. Let $f \co \cX \to \cY$ be a representable, quasi-finite and separated morphism. If $\cY$ is locally reductive, then $\Sigma_f \subset |\cX|$ is locally constructible.
\end{prop}

\begin{proof}
	We may assume that $S$ is quasi-compact.
	The hypotheses imply that there exists a representable, \'etale and surjective morphism $g \co \cY' \to \cY$ where $\cY' \cong [\Spec(A) / \GL_N]$.
	Let $\cX' = \cX \times_{\cY} \cY'$ with projections  $g' \co \cX' \to \cX$ and $f' \co \cX' \to \cY'$.
	By \Cref{L:theta-morphism}\ref{L:theta-morphism-representable-etale-surjective}, the morphism $\uMap_S(\Theta,\cY') \to \uMap_S(\Theta,\cY)$ is surjective.
	Since $\Sigma_f$ is the locus of points $x \in |\cX|$ where there is a diagram \eqref{E:theta-surjective} that does not lift, $\Sigma_f = g'(\Sigma_{f'})$. By Chevalley's Theorem and \Cref{L:sigma-closed}, the locus $\Sigma_f$ is constructible.

Alternatively, one could invoke \cite[Lem.~4.4.8]{hlinstability}, which implies that the image in $|\cX|$ of any open and closed substack of $\cX \times_\cY \uMap_S(\Theta,\cY)$ is  locally constructible.
\end{proof}

Let us give some simple examples and non-examples of $\Theta$-surjectivity.  In these examples, we work over a field.

\begin{lem} \label{E:relative4} 
	Let $\pi \co \cX \to X$ be an adequate moduli space of finite type {such that $X$ is a noetherian algebraic space and $\cX$ is a} locally reductive stack.  If $\cU \subset \cX$ is an open substack, then $\cU$ is saturated (i.e. $\pi^{-1}(\pi(\cU)) = \cU$) if and only if $\cU \hookarr \cX$ is $\Theta$-surjective.  In this case, $\cU \hookarr \cX$ is even a $\Theta$-isomorphism. 
\end{lem}

\begin{proof}
The ``$\implies$" direction is clear.  {For the converse suppose that $u\in \cU(\kbar)$ is a geometric point and  $x_0\in \cX(\kbar)$ is a geometric point over a closed point in the  fiber of $\pi^{-1}(\pi(u))$. Since $\cX$ is locally reductive and $\pi \co \cX \to X$ is of finite type, \Cref{L:isotrivial} implies that the specialization $u \rightsquigarrow x_0$ is realized by a morphism $\Theta_{\kbar} \to \cX$. As $\cU\hookrightarrow \cX$ is $\Theta$-surjective this implies that $x_0\in \cU(\kbar)$.} 
\end{proof}

\begin{example} \label{E:relative1}
The open immersion $\Spec(k) \hookarr  [\bA_k^1 / \bG_m]$ is $\Theta$-reductive but \emph{not} $\Theta$-surjective.  Indeed, this is the prototypical example of a morphism that does not send closed points to closed points.
\end{example}

\begin{example} \label{E:relative2}
	 Consider the action of $\bG_m$ on $X=\bA^2 \smallsetminus 0$ via $t \cdot (x,y) = (tx, y)$ (as in \Cref{E:thetaR})) and the open immersion $f \co \bA^1 \hookarr  [X/\bG_m]$
of the locus where $x$ is non-zero.  Then
	 $$\ev(f)_1 \co \bA^1 = \uMap(\Theta,\bA^1) \to \uMap(\Theta,[X/\bG_m]) =  \bA^1 \sqcup \big( \bigsqcup_{n < 0}  \bA^1 \smallsetminus 0 \big) $$
	 which is the inclusion onto the first factor.  Again, $f$ is affine and hence $\Theta$-reductive but not $\Theta$-surjective.
\end{example}

\begin{example} \label{E:relative5}
	{The nodal cubic  $C$ can be described as taking the union of two coordinate axes $W=\Spec(k[x,y]/xy)$ and identifying the open subsets $\Spec(k[x,x^{-1})], \Spec(k[y,y^{-1}])\subset W$ via $x=y^{-1}$. In particular the $\bG_m$-action with weights $1$ and $-1$ on $x$ and $y$ on $W$ defines a $\bG_m$-action on $C$ for which the node is the only fixed point.
		
 The étale presentation $f\colon [W/\bG_m]\to [C/\bG_m]$ maps closed points to closed points but it is \emph{not} $\Theta$-surjective, because there is no lift in the diagram
$$\xymatrix{
\Spec(k)=[\Spec(k[x,x^{-1}])/\bG_m] \ar[r] \ar[d]^{x\mapsto y^{-1}}	& [\Spec( k[x,y]/xy) / \bG_m]  \ar[d]^{f}
\ar[d]\\
\Theta\cong[\Spec(k[y])/\bG_m] \ar[r] \ar@{-->}[ur]	& [C/\bG_m]
},$$
i.e., the open point $\Spec(k)\subset [C/\bG_m]$ has two preimages in $[W/\bG_m]$, but for each of these only one of the two possible extensions to a map $\Theta \to [C/\bG_m]$ lifts to $[W/\bG_m]$.}
\end{example}

\subsection{Elementary modifications and  {\textsf{S}}-complete morphisms}

\subsubsection{Modifications and elementary modifications}

As in \cite[\S 2.B]{heinloth-stability} the following stack, which depends on a choice of DVR $R$, plays an important role in our analysis of criteria for separatedness of good moduli spaces.
\begin{equation} \label{E:ST}
\oST_R := [\Spec \big(R[s,t] / (st-\pi)\big)  / \bG_m],
\end{equation}
where $s$ and $t$ have $\bG_m$-weights $1$ and $-1$ respectively, and $\pi$ is a choice of uniformizer for $R$. A different choice of $\pi$ results in an isomorphic stack. 

Observe that $\oST_R \smallsetminus 0 \cong \Spec(R) \cup_{\Spec(K)} \Spec(R)$, where $K$ is the fraction field of $R$, because the locus where $s \neq 0$ in $\oST_R$ is isomorphic to $[\Spec \big(R[s,t]_s/(t-\pi/s)\big) / \bG_m] \cong [\Spec(R[s]_s)/\bG_m] \cong \Spec(R)$ and the locus where $t \neq 0$ has a similar description. 
A morphism $h \co \oST_R \smallsetminus 0 \to \cX$ to an algebraic stack is the data of two morphisms $\xi,\xi' \co \Spec(R) \to \cX$, where $\xi := h|_{\{s \neq 0\}}$ and $\xi' := h|_{\{t\neq 0\}}$, together with an isomorphism $\xi_K \simeq \xi'_K$. 

\begin{defn} \label{D:elementary_mod} \label{D:modifications}
Let $\cX$ be an algebraic stack and let  $\xi\co \Spec(R) \to \cX$ be a morphism where $R$ is a DVR with fraction field $K$.
\begin{enumerate}
	\item A \emph{modification} of $\xi$ is the data of a morphism $\xi' \co \Spec(R) \to \cX$ along with an isomorphism between the restrictions $\xi|_K \simeq \xi'|_K$.
	\item An \emph{elementary modification} of $\xi$ is the data of a morphism $h \co \oST_{R} \to \cX$ along with an isomorphism $\xi \simeq h|_{\{s\neq 0\}}$.
\end{enumerate}
\end{defn}

An elementary modification is clearly also a modification.

\begin{remark} \label{R:langton}
The terminology here is inspired by the terminology of \cite{langton}, but does not exactly coincide. Langton's notion of ``elementary modifications" of families of vector bundles over a DVR are examples of the notion of elementary modification above which flip two-step filtrations. 
To see this, let $X$ be a noetherian scheme and $\uCoh(X)$ the stack of coherent sheaves on $X$. Let $R$ be a DVR with fraction field $K$ and residue field $\kappa$. A quasi-coherent sheaf on $X \times \oST_R$ corresponds to a $\bZ$-graded coherent sheaf $\bigoplus_{n \in \bZ} F_n$ on $X_R$ together with a diagram of maps
  \begin{equation}
\xymatrix{\cdots \ar@/^/[r]^s & F_{n-1} \ar@/^/[r]^s \ar@/^/[l]^t & F_{n} \ar@/^/[r]^s \ar@/^/[l]^t & F_{n+1}  \ar@/^/[r]^s \ar@/^/[l]^t & \ar@/^/[l]^t \cdots },
\end{equation}
such that $st=ts=\pi$.  Moreover, $F$ is coherent if each $F_n$ is coherent, $s \co F_{n-1} \to F_{n}$ is an isomorphism for $n \gg 0$ and $t \co F_n \to F_{n-1}$ is an isomorphism for $n \ll 0$.  The sheaf $F$ is a flat over $\oST_R$ if and only if the maps $s$ and $t$ are injective, and the induced map $t \co F_{n+1}/sF_n \to F_n / sF_{n-1}$ is injective. 
(See \Cref{C:reese} for a proof that these properties characterize coherence and flatness of $F$.)

Suppose that we have a coherent sheaf $E$ on $X_R$ which is flat over $R$ whose restriction $E_\kappa$ to $X$ fits into a short exact sequence
\begin{equation} \label{E:exact1}
	0 \to B \to E_\kappa \to G \to 0.
\end{equation}
(In Langton's algorithm, one takes $B \subset E_\kappa$ to be the maximal destabilizing subsheaf.)
Let $E' = \ker(E \to E_\kappa \to G).$  Then $E'$ is flat over $R$ and $E_K = E'_K$.  Moreover, we have that $\pi E \subset E' \subset E$ with $E/E' = G$ and $E'/\pi E = B$; this implies that $E'_\kappa$ fits into a short exact sequence
\begin{equation} \label{E:exact2}
	0 \to G \to E'_\kappa \to B \to 0.
\end{equation}

This data defines a coherent sheaf on $X \times \oST_R$ flat over $\oST_R$ as follows: 
set $F_n$ to be $E$ if $n \ge 0$ and $E'$ if $n < 0$.  Let $s$ and $t$ act via
$$\xymatrix{
\cdots \ar@/^/[r]^{\pi} & 	 E' \ar@/^/[r]^{\pi} \ar@/^/[l]^1 & E' \ar@/^/[r]^{1} \ar@/^/[l]^1 & E \ar@/^/[r]^{1} \ar@/^/[l]^{\pi} & E \ar@/^/[r]^1 \ar@/^/[l]^{\pi} & \cdots \ar@/^/[l]^{\pi}
}$$

In general, the restriction of $F$ to $\{s \neq 0\} = \Spec(R)$ is the colimit over the $\bZ$-sequence of maps $s \co F_n \to F_{n+1}$.  In our case, this restriction is $E$.  Likewise, the restriction of $F$ to $\{t \neq 0\}$ is the colimit over the $t$ maps so its $E'$ in our case.
In general, the restriction of $F$ to $\{s = 0\} = \Theta_\kappa$ is the (generalized) $\bZ$-filtration
$$\cdots \leftarrow  F_n / sF_{n-1} \xleftarrow{t} F_{n+1} / sF_{n+2}  \leftarrow \cdots$$
which in our case corresponds to $E'_\kappa \supseteq G$ of \eqref{E:exact1}. Similarly, the restriction of $E$ to $\{t=0\}$ corresponds to the filtration $B \subset E_\kappa$ of \eqref{E:exact2}.  

An analogous construction shows that any finite sequence of steps in Langton's algorithm can be realized by a \emph{single} elementary modification.
\end{remark}

\subsubsection{\Scomplete\ morphisms}

\begin{defn} \label{D:S-complete}
We say that a morphism $f \co \cX \to \cY$ of locally noetherian algebraic stacks is \emph{\Scomplete} if for any DVR $R$ and any commutative diagram
\begin{equation} \label[diagram]{E:S-complete}
\begin{split}
\xymatrix{
\oST_R  \smallsetminus 0 \ar[r] \ar[d]					& \cX \ar[d]^f \\
\oST_R \ar[r] \ar@{-->}[ur]		& \cY
}
\end{split} \end{equation}
of solid arrows, there exists a unique dotted arrow filling in the diagram.
\end{defn}

\begin{remark} The motivation for the terminology ``\Scomplete" comes from Seshadri's work on the {\sf S}-equivalence of semistable vector bundles.  Namely, if $\cX$ is the moduli stack of semistable vector bundles over a smooth projective curve $C$ over $k$, then $\cX$ is {\textsf{S}}-complete (see e.g., \Cref{L:BunGScomplete}).  If $R$ is a DVR with fraction field $K$ and residue field $k$, and $\cE, \cF$ are two families of semistable vector bundles on $C_R$ which are isomorphic over $C_K$, then {\textsf{S}}-completeness implies that the special fibers $\cE_0$ and $\cF_0$ on $C$ are {\textsf{S}}-equivalent.
\end{remark}

\begin{remark} {\textsf{S}}-complete morphisms are stable under composition and {also under} base change {provided that the locally noetherian assumption is preserved, e.g. along morphisms essentially of finite type.}  A morphism of quasi-separated and locally noetherian algebraic spaces is {\textsf{S}}-complete if and only if it is separated (\Cref{P:quasi-finite-S-complete}).  While affine morphisms are always {\textsf{S}}-complete (\Cref{P:S-complete}(\ref{P:S-complete-affine})), it is not true that separated, representable morphisms are {\textsf{S}}-complete.  For instance, the open immersion $\oST_R \smallsetminus 0 \to \oST_R$ is not {\textsf{S}}-complete.  This example also shows that {\textsf{S}}-complete morphisms do not satisfy smooth descent; however, {\textsf{S}}-completeness does descend along representable, \'etale and surjective morphisms (\Cref{P:S-complete-descent}).
\end{remark}

We now state properties of {\textsf{S}}-completeness analogous to \Cref{P:theta-reductive-extensions}, \Cref{P:theta-reductive-essentially-finite-type-DVRs}, \Cref{P:theta-reductiveheta-reductive-descent}, and \Cref{P:theta-reductive}.  In each case, the proof is either identical or requires a mild modification.
\begin{prop} \label{P:S-complete-extensions}
Let $ \cX \to \cY$ be a morphism of locally noetherian algebraic stacks, and consider a diagram of the form \eqref{E:S-complete}, then there exists a unique dotted arrow filling the diagram if either
\begin{enumerate}
\item there exists a unique filling after passing to an unramified extension $R \subset R'$ of DVRs which is an isomorphism on residue fields, such as the completion of $R$, or
\item $\cX \to \cY$ has affine diagonal, and there exists a filling after an arbitrary extension of DVRs $R \subset R'$.
\end{enumerate}
In particular, to verify that a morphism of locally noetherian algebraic stacks is {\textsf{S}}-complete, it suffices to check the lifting criterion \eqref{E:S-complete} for complete DVRs. \epf
\end{prop}

\begin{prop} \label{P:S-complete-essentially-finite-type-DVRs}
{Let $S$ be a {quasi-}excellent and quasi-separated algebraic space $S$.}
Let $f \co \cX \to \cY$ be a quasi-compact morphism of algebraic stacks locally of finite type and with affine diagonal over {$S$.} {Assume that $\cX$ admits a representable \'etale cover by quotient stacks.} To verify that $\cX \to \cY$ is \textsf{S}-complete, it suffices to check the lifting criterion \eqref{E:S-complete} for DVRs $R$ essentially of finite type over $S$.
\end{prop}

\begin{proof}
We first reformulate {\textsf{S}}-completeness as a lifting problem for DVRs. Let $[\bA^2/\bG_m]$ be the quotient of $\bA^2$ by the anti-diagonal action of $\bG_m$ and $p\colon [\bA^2/\bG_m]\to \bA^1$ the good moduli space, so that for any DVR $R$ with local parameter $\pi$ we get $\oST_R= [\bA^2/\bG_m] \times_{\bA^1} \Spec(R)$.

Then the mapping stack $$\mathbf{S}(\cX):=\uMap_{\bA^1}([\bA^2/\bG_m],\cX \times \bA^1)$$is an algebraic stack locally of finite type with affine diagonal \cite[Thm. 14.9]{ahr2} over $\bA^1$ and a morphism $\oST_R \to \cX$ is a section
$$\xymatrix{ & \mathbf{S}(\cX) \ar[d]\\ 
	\Spec(R) \ar[r]^{\pi \to t}\ar[ur] & \bA^1.}$$
The analog of $\ev(f)_1$ is the morphism 
$$\ev(f)_{(x,1),(1,y)} \colon \mathbf{S}(\cX) \to \mathbf{S}(\cY) \times_{\cY \times \cY} (\cX \times \cX)$$
which is given by composition on the first factor and evaluation on the two sections $\bA^1 \to [\bA^2/\bG_m]$ on the second factor.

The defining lifting problem of \textsf{S}-completeness \eqref{E:S-complete} then translates into the lifting problem:
\begin{equation}\label{E:S-completeness_lifting_problem}
\xymatrix{
	\Spec K \ar[d] \ar[r] &  \mathbf{S}(\cX) \ar[d]^{\ev(f)_{(x,1),(1,y)}} \\ 
	\Spec R \ar[r]\ar@{-->}[ur]  & \mathbf{S}(\cY) \times_{\cY \times \cY} (\cX \times \cX)
}
\end{equation}
As before the morphism $\ev_{(x,1),(1,y)}(f)$ is representable, but not necessarily quasi-compact, because the fiber over $0\in \bA^1$ is $\uMap_S(\Theta,\cX) \times_{\uMap_S(B\bG_m,\cX)} \uMap_S(\Theta,\cX)$.

First we prove the claim for a quotient stack $\cX = [X/G]$. By replacing $X$ with $(X \times \GL_n) / G$, we may assume that $G = \GL_n$. For any DVR $R$ we consider the composition of commutative squares
$$\xymatrix{
	\Spec K \ar[d] \ar[r] &  \mathbf{S}(\cX) \ar[d]^{\ev(f)_{(x,1),(1,y)}} \ar[r] &  \mathbf{S}(B\GL_n) \ar[d]^{\ev(B\GL_n \to S)_{(x,1),(1,y)}} \\ 
	\Spec R \ar[r]\ar@{-->}[ur]  & \mathbf{S}(\cY) \times_{\cY^2} \cX^2 \ar[r] & \bA^1 \times (B\GL_n)^2
}.$$
The outer square always admits a unique lift, because $B\GL_n$ is {\textsf{S}}-complete, so the lifting problem \eqref{E:S-completeness_lifting_problem} is equivalent to the lifting problem
\begin{equation} \label{E:S-completeness_lifting_problem_2}
\xymatrix{
	\Spec K \ar[d] \ar[r] &  \mathbf{S}(\cX) \ar[d] \\ 
	\Spec R \ar[r]\ar@{-->}[ur]  & (\mathbf{S}(\cY) \underset{\cY^2}{\times} \cX^2) \underset{(B\GL_n)^2}{\times}  \mathbf{S}(B\GL_n) }
\end{equation}

The fiber of $\textbf{S}(B\GL_n)$ over $0 \in \bA^1$ is
\[
\textbf{S}(B\GL_n)_0 \cong \uMap_S(\Theta,B\GL_n) \times_{\uMap_S(B\bG_m,B\GL_n)} \uMap_S(\Theta, B\GL_n),
\]
which is an infinite disjoint union of quasi-compact stacks corresponding to conjugacy classes of cocharacters of $\GL_n$ \cite[Thm.~1.4.7]{hlinstability}. For every such cocharacter $\lambda : \bG_m \to \GL_n$, the union of the corresponding component of $\mathbf{S}(B\GL_n)_0$ with the preimage of $\bA^1 \setminus 0$ in $\mathbf{S}(B\GL_n)$ defines a quasi-compact open substack $\cU_\lambda \subset \mathbf{S}(B\GL_n)$. The open the substack of $\mathbf{S}([X/\GL_n])$ lying over $\cU_\lambda$ has an explicit quasi-compact description in terms of $\lambda$-equivariant morphisms into $X$ (e.g., as in \cite[Theorem 2.2.2]{drinfeld}). Therefore the vertical morphism in \eqref{E:S-completeness_lifting_problem_2} is quasi-compact, and the lifting problem can be reduced to the lifting problem for a quasi-compact morphism of algebraic spaces, as in the proof of \Cref{P:theta-reductive-essentially-finite-type-DVRs}. The claim then follows from \Cref{lem:ess_fintype_criterion}.

For the general case, choose a separated representable \'etale cover of $\cX$ by quotient stacks $[X/G]$, and assume that $f$ satisfies the lifting criterion \eqref{E:S-complete} for all essentially finite type DVRs. By \Cref{P:S-complete-extensions}(2), it suffices to verify the lifting criterion \eqref{E:S-complete} for $\cX \to \cY$ and an arbitrary DVR $R$ after passing to an extension of $R$. The composition $[X/G] \to \cY$ satisfies the lifting criterion for essentially finite type DVRs because $[X/G] \to \cX$ is separated, and hence for all DVRs by the argument above. Therefore to complete the proof, it suffices to show that given a lifting diagram \eqref{E:S-complete} for an arbitrary DVR $R$, the morphism $\oST_R \setminus 0 \to \cX$ lifts to some $[X/G]$ after passing to an extension of $R$.

It suffices to find a point $\eta \in |\cX|$ that is a common specialization of the image of the closed point under the two maps $\Spec(R) \to \cX$, so that the image of $\oST_R \setminus 0 \to \cX$ is contained in the image of some $[X/G] \to \cX$. We argue as in the proof of \Cref{lem:ess_fintype_criterion} that given a lifting diagram \eqref{E:S-completeness_lifting_problem}, one can construct an equivalent lifting diagram by replacing $K$ with a subfield that is essentially finite type over $S$ and replacing $R$ with its intersection with this subfield. Then we write $R$ as a filtered union of essentially finite type rings $A_{\alpha}$ that all contain a uniformizer $\pi \in R$ and have fraction field $K$. The morphism $\Spec(R) \to \mathbf{S}(\cY) \times_{\cY^2} \cX^2$ must factor through some morphism $\Spec(A_\alpha) \to \mathbf{S}(\cY) \times_{\cY^2} \cX^2$. Then if $\mathfrak{p} \subset A_\alpha$ is the image of the closed point of $\Spec(R)$, the normalization of the localization of the affine blowup $A_\alpha[\frac{\mathfrak{p}}{\pi}]_{(\pi)}$ is an essentially finite type DVR $R'$ with fraction field $K$ that fits into a lifting diagram \eqref{E:S-completeness_lifting_problem}. By hypothesis a lift exists for \eqref{E:S-completeness_lifting_problem} and hence for \eqref{E:S-complete}, which implies that the images of the closed point under the two maps $\Spec(R') \to \cX$ have a common specialization $\eta \in |\cX|$. By construction, the images of the closed point of $\Spec(R')$ agree with those of the closed point of $\Spec(R)$ under the original maps, so $\eta$ is a common specialization of these points as well.
\end{proof}

\begin{prop} \label{P:S-complete-descent} Let $\cX \to \cY$ be a morphism of locally noetherian algebraic stacks.
\begin{enumerate}
\item  If $\cY' \to \cY$ is an \'etale, representable and surjective morphism, then $\cX \to \cY$ is {\textsf{S}}-complete if and only if $\cX \times_{\cY} \cY' \to \cY'$ is {\textsf{S}}-complete.
\item If $\cX' \to \cX$ is a finite, \'etale and surjective morphism, then $\cX \to \cY$ is {\textsf{S}}-complete if and only if $\cX' \to \cY$ is {\textsf{S}}-complete. \epf
\end{enumerate} 
\end{prop}

\begin{prop} \label{P:S-complete} \quad
\begin{enumerate}
\item  \label{P:S-complete-affine}
	An affine morphism of locally noetherian algebraic stacks is {\textsf{S}}-complete.
\item \label{P:S-complete-quotient}
		Let $S$ be a locally noetherian scheme.  Let $G \to S$ be a geometrically reductive and \'etale-locally embeddable group scheme (e.g. reductive) acting on a locally noetherian scheme $X$ affine over $S$.  Then the morphism $[X/G] \to S$ is {\textsf{S}}-complete.
\item \label{P:S-complete-gms}
	If $\cX$ is a locally reductive, locally noetherian algebraic stack with separated diagonal, then a locally noetherian adequate moduli space $\cX \to X$ is {\textsf{S}}-complete.	\epf
\end{enumerate} 
\end{prop}

We now detail additional important properties of {\textsf{S}}-completeness.

\begin{lem} \label{L:quasi-finite-factoring-S-complete}
	Let $S$ be a quasi-separated algebraic space. Let $\cX$ be an algebraic stack locally of finite presentation and quasi-separated over $S$ with quasi-finite inertia.  If $T$ is an algebraic space over $S$, then any morphism $\oST_T \to \cX$ factors uniquely through $\oST_T \to T$.
\end{lem}

\begin{proof} This can be established with the same method as \Cref{L:quasi-finite-factoring-theta-reductive}.
\end{proof}

\begin{prop}  \label{P:quasi-finite-S-complete}
	Let $S$ be a noetherian algebraic space, and 
	let $f \co \cX \to \cY$ be a morphism of algebraic stacks locally of finite presentation and quasi-separated over $S$ with quasi-finite inertia.  Then $f$ is {\textsf{S}}-complete if and only if $f$ is separated.
\end{prop}

\begin{proof}
Let $R$ be a DVR with fraction field $K$.  By \Cref{L:quasi-finite-factoring-S-complete}, any morphism from $\oST_R$ to $\cX$ or $\cY$ factors through $\oST_R \to \Spec(R)$. Because $\oST_R \smallsetminus 0 = \Spec(R) \bigcup_{\Spec(K)} \Spec(R)$, we see that the valuative criterion of \Cref{E:S-complete} is equivalent to the valuative criterion for separatedness.
\end{proof}

\begin{prop}\label{P:SCompleteReductive}  If $G$ is an {affine} algebraic group over a field $k$, then $G$ is geometrically reductive if and only if $B_k G$ is {\textsf{S}}-complete, and if $k$ is perfect it suffices to check the definition of \textsf{S}-completeness for any fixed DVR over $k$.  In particular, a closed point of an {\textsf{S}}-complete locally noetherian algebraic stack with affine stabilizers has a geometrically reductive stabilizer.
\end{prop}

\begin{proof} From \Cref{P:S-complete}\eqref{P:S-complete-quotient}, we know that if $G$ is geometrically reductive, then $B_k G$ is {\textsf{S}}-complete.  For the converse, we may assume that $k$ is algebraically closed. Note that $B\bG_a$ is not \textsf{S}-complete, because any morphism $\oST_R \to B\bG_a$ factors through the good moduli space $\Spec(R)$ of $\oST_R$, as any morphism $B\bG_m\to B\bG_a$ is constant. However, $\bG_a$-torsors on $\oST_R \smallsetminus 0 = \Spec(R)\cup_{\Spec(K)}\Spec(R)$ correspond to elements in $K/R$, and in particular there exist non-trivial torsors.

Now suppose that $G$ is not geometrically reductive.  Then by considering the unipotent radical $R_u(G)$ of the reduced group scheme $G^{\red}$, the induced morphism $B_k R_u(G) \to B_k G$ is affine.  Similarly, by taking a normal subgroup $\bG_a \subset R_u(G)$, there is an affine morphism $B_k \bG_a \to B_k R_u(G)$.  The composition $B_k \bG_a \to B_k R_u(G) \to B_k G$ is affine.  Since $B_k G$ is {\textsf{S}}-complete, by \Cref{P:S-complete}\eqref{P:S-complete-affine} so is $B_k \bG_a$, a contradiction.

If $k$ is perfect, a normal subgroup $\bG_a \subset R_u(G)$ already exists over $k$ (\cite[Corollary 14.3.10]{springer}) and the argument gives a non-extending torsor over any DVR.
\end{proof}

Expanding on \Cref{P:S-complete}\eqref{P:S-complete-gms}, we have the following criterion for when an adequate moduli space is separated.

\begin{prop} \label{P:gms-S-complete}
	Let $\cX$ be a locally reductive, {locally noetherian, and quasi-separated} algebraic stack with separated diagonal, and let $\cX \to X$ be a {locally noetherian} adequate moduli space.  Then
\begin{enumerate}
\item \label{P:gms-S-complete1}
the morphism $\cX \to X$ is {\textsf{S}}-complete;
\item \label{P:gms-S-complete2}
the morphism $X \to S$ is separated if and only if $\cX \to S$ is {\textsf{S}}-complete; and
\item \label{P:gms-S-complete3}
the  morphism $X \to S$  is proper if and only if $\cX \to S$ is of finite type, universally closed and {\textsf{S}}-complete.
\end{enumerate}
\end{prop}

\begin{proof} Part \eqref{P:gms-S-complete1} is \Cref{P:S-complete}\eqref{P:S-complete-gms}.
The implication `$\Rightarrow$' in Part \eqref{P:gms-S-complete2} follows from Part \eqref{P:gms-S-complete1} and the fact that separated algebraic spaces are {\textsf{S}}-complete.  Conversely, suppose $\cX$ is {\textsf{S}}-complete.  Suppose $f, g \co \Spec(R) \to X$ are two maps such that $f|_K = g|_K$.  After possibly an 
extension of $R$, we may choose a lift  $\Spec(K) \to \cX$ of $f|_K=g|_K$.  Since $\cX \to X$ is universally closed, after possibly further extensions of $R$, we may choose lifts $\widetilde{f}, \widetilde{g} \co \Spec(R) \to \cX$ of $f, g$ such that $\widetilde{f}|_K \cong \widetilde{g}|_K$.  By applying the {\textsf{S}}-completeness of $\cX$, we can extend $\widetilde{f}, \widetilde{g}$ to a morphism $\oST_R \to \cX$.  As $\oST_R \to \Spec(R)$ is a good moduli space and hence universal for maps to algebraic spaces \cite[Thm.~6.6]{alper-good}, the morphism $\oST_R \to \cX$ descends to a unique morphism $\Spec(R) \to X$ which necessarily must be equal to both $f$ and $g$.  We conclude that $X \to S$ is separated by the valuative criterion for separatedness.   Part  \eqref{P:gms-S-complete3} follows from Part \eqref{P:gms-S-complete1} using the fact that $X \to S$ is of finite type (resp. universally closed) if and only if $\cX \to S$ is \cite[Thms.~5.3.1 and 6.3.3]{alper-adequate}.
\end{proof}

\begin{cor} \label{C:modifications_GMS}
Let $\cX$ be a locally reductive, locally noetherian algebraic stack with separated diagonal, and let $\cX \to X$ be a locally noetherian adequate moduli space. Let $R$ be any DVR and consider two morphisms $\xi_0,\xi_1 \co \Spec(R) \to \cX$ with $(\xi_0)|_K \cong (\xi_1)|_K$  Then the following are equivalent:
\begin{enumerate}
\item $\xi_0$ and $\xi_1$ differ by an elementary modification,
\item $\xi_0$ and $\xi_1$ differ by a finite sequence of elementary modifications,
\item the compositions $\xi_i : \Spec(R) \to \cX \to X$ agree for $i=0,1$.
\end{enumerate}
\end{cor}

\begin{proof}
Clearly $(1) \Rightarrow (2)$. The projection $\oST_R \to \Spec(R)$ is a good moduli space and hence universal for maps to algebraic spaces \cite[Thm.~6.6]{alper-good}. It follows that any two maps which differ by an elementary modification induce the same $R$-point of $X$, {thus all the maps in a finite sequence of elementary modifications also induce the same $R$-point of $X$, i.e.,} $(2) \Rightarrow (3)$.  The implication $(3) \Rightarrow (1)$ follows from \Cref{P:S-complete}\eqref{P:S-complete-gms}.
\end{proof}

\begin{remark}
The above conditions are not equivalent to saying that $\xi_0$ and $\xi_1$ are modifications such that the closures of $\xi_0(0)$ and $\xi_1(0)$ intersect. This is because there exist algebraic spaces $X$ admitting  two distinct maps $\xi_0, \xi_1 \co \Spec(R) \to X$ with $\xi_0|_K = \xi_1|_K$ and $\xi_0(0) = \xi_1(0)$. An example of this is a bug-eyed cover,  the algebraic space $X$ obtained by taking the free $\bZ/2$-quotient of the non-separated affine line, where the action of $\bZ/2$ is via $x \mapsto -x$ and swaps the origins. Here the diagonal of $X$ is affine but not a locally closed immersion.
\end{remark}

\begin{remark}[Hartogs' principle]
Both $\Theta$-reductivity and {\textsf{S}}-completeness are conditions asserting the existence and uniqueness of extending morphisms along a codimension two locus.  One might be tempted to unify these two notions by defining that a morphism $f \co \cX \to \cY$ of locally noetherian algebraic stacks satisfies \emph{Hartogs' principle} if for any regular local ring $S$ of dimension $2$ with closed point $0 \in \Spec(S)$, there exists a unique dotted arrow filling in any commutative diagram
\begin{equation} \label[diagram]{E:Hartogs}
\begin{split}
\xymatrix{
\Spec(S) \smallsetminus 0 \ar[r] \ar[d]					& \cX \ar[d]^f \\
\Spec(S) \ar[r] \ar@{-->}[ur]		& \cY
}
\end{split} \end{equation}
of solid arrows.  Any such morphism is necessarily both $\Theta$-reductive and {\textsf{S}}-complete.  Moreover, the analogues of \Cref{P:theta-reductive} and \Cref{P:S-complete} hold for such morphisms.  However, many algebraic stacks (e.g. the stack $\uCoh(X)$ of coherent sheaves on a proper scheme $X$ over a field $k$) are both $\Theta$-reductive and {\textsf{S}}-complete but do \emph{not} satisfy Hartogs' principle.
\end{remark}

\subsection{Unpunctured inertia}

We now give the last of the properties that will turn out to be necessary for the existence of good moduli spaces.  

\begin{defn} \label{D:unpunctured}
	We say that a noetherian algebraic stack has \emph{unpunctured inertia} if for any closed point $x \in |\cX|$ and any {formally} smooth morphism $p \co (U,u) \to (\cX, x)$, where $U$ is the spectrum of a local ring {of a smooth neighborhood of $x$} with closed point $u$, each connected component of the inertia group scheme $\Aut_{\cX}(p) \to U$ has non-empty intersection with the fiber over $u$.
\end{defn}

\begin{remark}
The condition of unpuncturedness is related to the property of purity of the morphism $\Aut_{\cX}(p) \to U$ as defined in \cite[\S 3.3]{raynaud-gruson} and further studied in \cite[\spref{0CV5}]{stacks-project}.  If $U$ is the spectrum of a strictly henselian local ring, then purity requires that if $s \in U$ is {\it any} point and $\gamma$ is an associated point in the fiber $\Aut_{\cX}(p)_s$, then the closure of $\gamma$ in $\Aut_{\cX}(p)$ has non-empty intersection with the fiber over the closed point $u$ of $U$.
\end{remark}

In \Cref{S:unpunctured}, we will provide valuative criteria which can be used to verify that a stack has unpunctured inertia.
In this section though, we provide only a few situations in which this condition is easy to check. 

\begin{prop} \label{P:quasi-finite-unpunctured}
	If $\cX$ is a noetherian algebraic stack with quasi-finite inertia, then $\cX$ has unpunctured inertia if and only if $\cX$ has finite inertia.  
\end{prop}

\begin{proof}
	If $\cX$ has finite inertia, then $\Aut_{\cX}(p) \to U$ is finite so clearly the image of each connected component contains the unique closed point $u \in U$.  For the converse, we may assume that $U$ is the spectrum of a Henselian local ring in which case $\Aut_{\cX}(p) = G \sqcup H$ where $G \to U$ finite and the fiber of $H \to U$ over $u$ is empty.  If $\Aut_{\cX}(p)$ is not finite, then $H$ is non-empty and any connected component of $H$ will have empty intersection with the fiber over $u$.
\end{proof}

\begin{prop} \label{P:connected-unpunctured}
	Let $\cX$ be a noetherian algebraic stack.  If $\cX$ has connected stabilizer groups, then $\cX$ has unpunctured inertia.
\end{prop}

\begin{proof}
	This is clear, by definition all fibers of $\Aut_{\cX}(p)\to U$ are connected, so any connected component of $\Aut_{\cX}$ intersects the component containing the identity section. 
\end{proof}
The following example shows that unpuncturedness need not be preserved when passing to open substacks.
\begin{example} \label{E:ex1}
	Consider the action of $G=\bG_m \rtimes \bZ/2$ on $X=\bA^2$ via $t \cdot (a, b) = (t a, t^{-1} b)$ and $ -1 \cdot (a, b) = (b, a)$.   Note that every point $(a,b) \in X$ with $ab \neq 0$ is fixed by the order $2$ element $(a/b,-1) \in G$.  
	The algebraic stack $[(X \smallsetminus 0) / G]$ does not have unpunctured inertia by \Cref{P:quasi-finite-unpunctured}.  However, it will follow from \Cref{P:proper-over-affine-mod-G} that $[X/G]$ has unpunctured inertia.  
\end{example}

\section{Existence of good moduli spaces} \label{S:existence}
  
The following is our first generalization of \Cref{T:existence} providing necessary and sufficient conditions for an algebraic stack to admit a good (resp. adequate) moduli space.

\begin{thm} \label{T:existence-general}
	Let $\cX$ be an algebraic stack locally of finite type with affine diagonal over a quasi-separated and locally noetherian algebraic space $S$.  Then $\cX$ admits a good moduli space $X$ if and only if 
	\begin{enumerate}
		\item every closed point of $\cX$ has linearly reductive stabilizer;
		\item $\cX$ is $\Theta$-reductive  (\Cref{D:theta-reductive}); and
		\item $\cX$ has unpunctured inertia (\Cref{D:unpunctured}).
	\end{enumerate} 
	If $\cX$ is locally reductive (\Cref{D:locally-reductive}), then $\cX$ admits an adequate moduli space $X$ if and only if (2) and (3) hold. In both cases, if $S$ is excellent it suffices to check the filling condition in \Cref{D:theta-reductive} only for DVRs essentially of finite type over $S$.
\end{thm}

We prove the sufficiency of these conditions in this section and postpone the proof of the necessity until the next section.
The idea of the proof of the existence of the adequate (resp. good) moduli space is simple. We use the slice theorem (\Cref{T:ahr}) to reduce to quotient stacks and glue the resulting moduli spaces. We need to apply the slice theorem carefully so that the \'etale quotient presentations preserve stabilizer groups and the topology of finite type points in order to ensure that the \'{e}tale covering of the stack induces an \'etale covering on the level of good moduli spaces.
  
\subsection{Reminder on maps inducing \'etale maps on good moduli spaces}
 
 If $f \co \cX \to \cY$ is a morphism of algebraic stacks and $x \in |\cX|$, we say that $f$ is \emph{stabilizer preserving at $x$} if there exists a representative $\widetilde{x} \co \Spec(l) \to \cX$ of $x$ (equivalently, for all representatives of $x$), the natural map $\Aut_{\cX}(\widetilde{x}) \to \Aut_{\cY}(f \circ \widetilde{x})$ is an isomorphism.
 
 \begin{prop}  \label{P:etale-preserving}

Let $\cX$ and $\cY$ be locally reductive, locally noetherian algebraic stacks.
Consider a commutative diagram 
\begin{equation} \label[diagram]{diagram-luna}
\begin{split}
\xymatrix{ 
\cX \ar[r]^f \ar[d]^{\pi_{\cX}}		& \cY \ar[d]^{\pi_{\cY}} \\
X \ar[r]^g					& Y
}
\end{split}
\end{equation}
where $f$ is representable, \'etale and separated. Assume that either $\pi_{\cX}$ and $\pi_{\cY}$ are good moduli spaces or that they are adequate moduli spaces of finite type with $X$ and $Y$ locally noetherian.
If $f$ is $\Theta$-surjective and $f$ is stabilizer preserving at every closed point of $\cX$, then $g$ is \'etale and  \Cref{diagram-luna} is Cartesian.
\end{prop}
 
 \begin{proof}
	Observe that since $f$ is $\Theta$-surjective, the image $f(x) \in |\cY|$ of every closed point $x \in |\cX|$ is closed (\Cref{L:closed-to-closed}).  The statement follows from a stack-theoretic generalization of Luna's fundamental lemma \cite[Theorem 3.14]{ahr2}.	
 \end{proof}

 \subsection{Proof of the existence result} \label{S:existence-proof}
 We first provide conditions on an algebraic stack ensuring that there are \'etale quotient presentations which are $\Theta$-surjective and stabilizer preserving.  This is the key ingredient in the proof of \Cref{T:existence}.
 
\begin{prop}  \label{P:quotient-presentations}
Let $S$ be a quasi-separated and locally noetherian algebraic space.  
	Let $\cY$ be an algebraic stack of finite type over $S$ with quasi-compact and separated diagonal over $S$, and with affine stabilizers. Let $f \co (\cX, x) \to (\cY,y)$ be a quasi-compact, separated, representable and \'etale morphism from a {locally reductive} algebraic stack $\cX$ such that there exists an adequate moduli space $\pi \co \cX \to X$ and $f$ induces an isomorphism of stabilizer groups at $x$ (e.g., $f$ is a separated, representable and \'etale quotient presentation). 
\begin{enumerate}
	\item \label{P:quotient-presentations-theta-surjective}
	If $\cY$ is $\Theta$-reductive, then there exists an open subspace $U \subset X$ containing $\pi(x)$ such that the restriction $f|_{\pi^{-1}(U)}$ is $\Theta$-surjective.
	\item \label{P:quotient-presentations-unpunctured}
	If $\cY$ has unpunctured inertia, then there exists an open 
	subspace $U \subset X$ containing $\pi(x)$ such that the restriction $f|_{\pi^{-1}(U)}$ induces an isomorphism $I_{\pi^{-1}(U)} \to \pi^{-1}(U) \times_{\cY} I_{\cY}$.
\end{enumerate}
\end{prop}

\begin{proof}  
For \eqref{P:quotient-presentations-theta-surjective} let us first show that the image of
$\ev(f)_1\co \uMap_S(\Theta,\cX) \to \cX \times_{\cY} \uMap_S(\Theta,\cY)$ is open and closed.  Since $f$ is \'etale, so is $\ev(f)_1$ by \Cref{L:basic}\eqref{L:basic-etale}, and so its image is open. As $\cX$ admits an adequate moduli space, it is $\Theta$-reductive by \Cref{P:theta-reductive}, and since $\ev_1 \co \uMap_S(\Theta, \cY) \to \cY$ satisfies the uniqueness part of the valuative criterion (as $\cY$ is $\Theta$-reductive), it follows that $\ev(f)_1$ 
satisfies the valuative criterion for properness and that its image is closed under specialization, hence closed.

Let $\cZ \subset \cX \times_{\cY} \uMap_S(\Theta,\cY)$ be the open and closed complement of the image of $\ev(f)_1$.  By \eqref{E:sigma-equality}, the image $p_1(\cZ) \subset |\cX|$ consists of the points where $f$ is not $\Theta$-surjective.
Since $\cX$ and $\cY$ are of finite type over $S$, we may apply \Cref{P:sigma-constructible} to conclude that the image $p_1(\cZ) \subset |\cX|$ is locally constructible. On the other hand, since $\cY$ is $\Theta$-reductive, the image $p_1(\cZ)$ is closed under specializations, hence closed. 

We now show that $f$ is $\Theta$-surjective at $x$, i.e., $x \notin p_1(\cZ)$. Consider a diagram of solid arrows
$$\xymatrix{
\Spec(k) \ar[r]^x \ar[d]							&  \cX \ar[d]^f \\
 \Theta_k \ar[r]^{\lambda} \ar@{-->}[ur]					& \cY
}$$
where $k$ is an algebraically closed field. As $y=f(x) \in \cY(k)$ is a closed point,  $\lambda$ factors through the natural map $B G_y \to \cY$. Since $f$ is stabilizer preserving at $x$, the induced map $B G_x \to B G_y$ is an isomorphism
so $\lambda$ lifts to a morphism $\Theta_k \to BG_x \to \cX$ filling in the dotted arrow.

Setting $\cU := \cX \setminus \pi^{-1}( \pi(p_1(\cZ)))$, we claim that $f|_{\cU} \co \cU \to \cY$ is $\Theta$-surjective.  Indeed, this follows from the lifting criterion for $\Theta$-surjectivity (\Cref{R:theta-surjective-valuative-criterion}) because $f$ is $\Theta$-surjective at all points of $\cU$ and the inclusion $j \co \cU \hookarr \cX$ is a $\Theta$-isomorphism (i.e., $\ev(j)_1$ is an isomorphism) as $\cU$ is the preimage of $X \setminus \pi(p_1(\cZ))$.

For \eqref{P:quotient-presentations-unpunctured},  it suffices to find an open neighborhood $\cU \subset \cX$ of $x$ such that $f|_{\cU} \co \cU \to \cY$ induces an isomorphism $I_{\cU} \to \cU \times_{\cY} I_{\cY}$.  We have a  Cartesian diagram
$$\xymatrix{
I_{\cX} \ar[r] \ar[d]		& \cX \times_{\cY} I_{\cY} \ar[d] \\
\cX	\ar[r]				 & \cX \times_{\cY} \cX.
}$$
Since $f$ is \'etale and separated, 
the morphism $I_{\cX}  \to \cX \times_{\cY} I_{\cY}$ is finite \'etale. If we set  $\cZ \subset \cX \times_{\cY} I_{\cY}$ to be the open and closed substack over which $I_{\cX}  \to \cX \times_{\cY} I_{\cY}$ is not an isomorphism, then $x$ is not contained in the image of $\cZ$ under $p_1 \co \cX \times_{\cY} I_{\cY} \to \cX$, as $f$ is stabilizer preserving at $x$.

{We claim that $x \notin \overline{p_1(\cZ)}$. To see this, let $(U, u) \to (\cY,y) $ be a formally smooth morphism, where $U$ is the spectrum of a local ring with closed point $u$ obtained from a smooth presentation of $\cY$. {Applying the condition that $\cY$ has unpunctured inertia to the pull-back $(\cX\times_{\cY} U,(x,u)) \to (\cY,y)$ of $(U,u)\to (\cY,y)$ we obtain that} any component of the preimage of $\cZ$ in $\cX \times_{\cY} I_{\cY} \times_{\cY} U$ must intersect the fiber over $u$ and would thus contain a point mapping to $x$ under $p_1$. As $x\not \in p_1(\cZ)$ this implies that $\cZ$ has empty preimage in $\cX \times_{\cY} I_{\cY} \times_{\cY} U$, i.e., $x \notin \overline{p_1(\cZ)}$.  Therefore,  if we set $\cU = \cX \smallsetminus \overline{p_1(\cZ)}$, the induced morphism $I_{\cU} \to \cU \times_{\cY} I_{\cY}$ is an isomorphism.}
\end{proof}

\begin{lem} \label{L:glue} Let $\cX$ be a locally noetherian algebraic stack with affine diagonal. Suppose that $\{\cU_i\}_{i \in I}$ is a Zariski-cover of $\cX$ such that each $\cU_i$ admits an adequate (resp. good) moduli space $\pi_i \co \cU_i \to U_i$ of finite type and each inclusion $\cU_i \hookarr \cX$ is $\Theta$-surjective.  Then $\cX$ admits an adequate (resp. good) moduli space $\pi \co \cX \to X$ of finite type.
\end{lem}

\begin{proof}
	Since each inclusion $\cU_i \cap \cU_j \hookarr \cU_i$ is $\Theta$-surjective, there exist open subspaces $U_{i,j} \subset U_i$ with $\pi_i^{-1}(U_{i,j}) = \cU_i \cap \cU_j$ (\Cref{E:relative4}).  By universality of adequate moduli spaces \cite[Thm.~3.12]{ahr2}, there are isomorphisms $U_{i,j} \iso U_{j,i}$ providing gluing data for an algebraic space $X$.  The morphisms $\pi_i$ glue to produce an adequate (resp. good) moduli space $\pi \co \cX \to X$.
\end{proof}

\medskip

\noindent
\emph{Proof of sufficiency in \Cref{T:existence-general}.}  
	We first claim that it suffices to prove the existence of an adequate moduli space.  Indeed, the condition that every closed point has linearly reductive stabilizer implies that $\cX$ is locally reductive (\Cref{T:ahr}) and that an adequate moduli space $\cX \to X$ is necessarily a good moduli space \cite[Thm.~9.3]{ahr2}.  

	We follow the proof of \cite[Thm.~2.1]{afs-existence}.  First observe that by \Cref{L:glue} and \Cref{P:quotient-presentations}(\ref{P:quotient-presentations-theta-surjective}), it suffices to show that every closed point $x \in |\cX|$ has an  open neighborhood $\cU$ admitting an adequate moduli space.
	Since $\cX$ is locally reductive and has affine diagonal, there exists an affine and \'etale quotient presentation $f \co (\cX_1, x_1) \to (\cX,x)$ with $\cX_1=[\Spec(A) / \GL_N]$ (see \Cref{rmk:local-structure-refinements}).
	By \Cref{P:quotient-presentations}, after replacing $\cX_1$ with an open neighborhood of $x_1$, we may assume that $f$ is $\Theta$-surjective and $f$ induces an isomorphism $I_{\cX_1} \to \cX_1 \times_{\cX} I_{\cX}$.
    After replacing $\cX$ with $f(\cX_1)$, we may assume that $f$ is surjective.  Since $I_{\cX_1} \to \cX_1 \times_{\cX} I_{\cX}$ is an isomorphism, every closed point of $[\Spec(A) / \GL_n]$ has geometrically reductive stabilizer. 
	We let $\pi_1 \co \cX_1 \to X_1 := \Spec(A^{\GL_N})$ be the induced adequate moduli space. 

	Set $\cX_2 = \cX_1 \times_{\cX} \cX_1$. The projections $p_1, p_2\co   \cX_2 \to \cX_1$ are also   \'etale, affine, surjective, and $\Theta$-surjective morphisms that induce isomorphisms $I_{\cX_2} \to \cX_2 \times_{\cX_1} I_{\cX_1}$.  Since $f$ is affine, $\cX_2 \to X_2:= \Spec( \Gamma(\cX_2, \oh_{\cX_2}) )$ is an adequate moduli space.
	By \Cref{P:etale-preserving}, both commutative squares in the diagram
	\begin{equation} \label{E:existence-etale-equivalence-relation}
	\begin{split}
	\xymatrix{
	\cX_2 \ar@<.5ex>[r]^{p_1} \ar@<-.5ex>[r]_{p_2} \ar[d]^{\pi_2}  
	& \cX_1 \ar[r]^f \ar[d]^{\pi_1}    & \cX \\
	X_2 \ar@<.5ex>[r]^{q_1} \ar@<-.5ex>[r]_{q_2}  & X_1
	}
	\end{split}
	\end{equation}
	are Cartesian.  Moreover, by the universality of adequate moduli spaces \cite[Thm.~3.12]{ahr2}, the \'etale groupoid structure on $\cX_2 \rightrightarrows \cX_1$ induces a \'etale groupoid structure on $X_2 \rightrightarrows X_1$.  The fact that $f$ induces isomorphisms of stabilizer groups implies that $\Delta\co  X_2 \to X_1 \times X_1$ is a monomorphism (see the argument of \cite[Prop.~3.1]{afs-existence}). Thus, 
	$X_2 \rightrightarrows X_1$ is an \'etale equivalence relation and there exists an algebraic space quotient $X$.  It follows from descent that there is an induced morphism $\pi \co \cX \to X$ which is an adequate moduli space. 

It suffices to check the definition of $\Theta$-reductivity with respect to essentially finite type DVRs by \Cref{P:theta-reductive-essentially-finite-type-DVRs}, because $\cX$ is locally reductive.
\epf

We note another consequence of \Cref{P:quotient-presentations}, which will be used in \S \ref{S:unpunctured} below.

\begin{prop} \label{P:GMS_one_point}
Let $\cX$ be an algebraic stack which is of finite type with affine diagonal over a field $k$.  Suppose that $\cX$ is $\Theta$-reductive and there exists a single closed point $x \in |\cX|$.  If $\cX$ is locally reductive (resp. the stabilizer of $x$ is linearly reductive), then $\cX$ admits an adequate (resp. good) moduli space.

If $k$ is algebraically closed and the stabilizer $G_x$ is linearly reductive, then $\cX \cong [\Spec(A)/G_x]$, and if in addition $\cX$ is reduced, then $\cX \to \Spec(k)$ is the good moduli space.
\end{prop}

\begin{proof}
Choose an \'etale quotient presentation $f \co (\cX_1,x_1) \to (\cX,x)$ with $\cX_1=[\Spec(B)/\GL_n]$ such that $x_1 \in |\cX_1|$ is the unique point mapping to $x$. Since $\cX$ is $\Theta$-reductive, by \Cref{P:quotient-presentations}\eqref{P:quotient-presentations-theta-surjective}, we can assume that $f$ is $\Theta$-surjective.  This implies that $f$ sends closed points to closed points and both projections $\cX_2 = \cX_1 \times_{\cX} \cX_1 \rightrightarrows \cX_1$ send closed points to closed points. Since both $\cX$ and $\cX_1$ have a unique closed point and $f$ induces an isomorphism of residual gerbes $\cG_{x_1} \to \cG_x$, it follows that $\cX_2$ has a unique closed point and that both projections $\cX_2 \rightrightarrows \cX_1$ induce isomorphism of stabilizers at this point.  Moreover, there are adequate moduli spaces $\cX_1 \to X_1$ and $\cX_2 \to X_2$.  As in the proof of \Cref{T:existence}, \Cref{P:etale-preserving} implies that the induced groupoid $X_2 \rightrightarrows X_1$ is an \'etale equivalence relation, and the quotient $X_1/X_2$ is an adequate moduli space for $\cX$. If the stabilizer of $x$ is linearly reductive, then $\cX \to X$ is necessarily a good moduli space \cite[Thm.~9.3]{ahr2}.

The final statement follows from  \cite[Thms.~2.9 and ~3.14]{ahr2} coupled with the observation that if $\cX$ is reduced, so is its good moduli space.
\end{proof}


\section{Criteria for unpunctured inertia} \label{S:unpunctured}

In this section we establish criteria which imply that a stack has unpunctured inertia. They are ``valuative criteria'' in the sense that they apply to families over discrete valuation rings. We will need the following notion:

\begin{defn} \label{D:nearby_modification}
Let $\cX$ be an algebraic stack over an algebraic space $S$, let $R$ be a DVR over $S$ with fraction field $K$ and residue field $\kappa$, and let $\xi \co \Spec(R) \to \cX$ be a morphism.  A {\em nearby modification} of $\xi$ is a morphism $\xi' \co \Spec(R') \to \cX$ where $R'/R$ is an extension of DVRs with fraction field $K'$ along with an isomorphism $\xi'|_{K'} \simeq \xi|_{K'}$ such that the closures of $\xi'(0)$ and $\xi(0)$ in $|\cX \times_{S} \Spec(\kappa)|$ have nonempty intersection. If $\xi'(0)$ lies in the closure of $\xi(0)$ we say that $\xi'$ is a {\em specializing modification} of $\xi$.
\end{defn}

Unlike in the definitions of modifications and elementary modifications in \Cref{D:modifications}, we have allowed extensions of DVRs above. The notion of nearby modifications and elementary modifications where there is no extension of DVRs is stronger than that of a modification and weaker than that of an elementary modification.  Note however that in an  {\textsf{S}}-complete stack, the notion of nearby modifications and elementary modifications where there is no extension of DVRs coincides with those in \Cref{D:nearby_modification}.

We can now state the main results of this section:

\begin{prop}\label{prop:crit-unpunctured}
Let $\cX$ be an algebraic stack, locally of finite type with affine diagonal over a quasi-separated and locally noetherian algebraic space $S$.
Consider the following conditions:
\begin{enumerate}
\item For any essentially finite type DVR $R$, any morphism $\xi \co \Spec(R) \to \cX$, and any $g \in \Aut_\cX(\xi_K)$ of finite order, there is a nearby modification $\xi'\co \Spec(R')\to \cX$ of $\xi$ such that $g|_{K'}$ extends to an automorphism of $\xi'$; \label{I:valuative_A}
\item For any essentially finite type DVR $R$, any morphism $\xi \co \Spec(R) \to \cX$, and any geometrically connected component $H \subset \Aut_\cX(\xi_K)$, there is a nearby modification  $\xi'\co \Spec(R')\to \cX$ of $\xi$ and some $g \in \Aut_\cX(\xi')$ such that $g|_{K'}$ lies in $H$; \label{I:valuative_B}
\item For any essentially finite type DVR $R$, any morphism $\xi \co \Spec(R) \to \cX$, and any geometrically connected component $H \subset \Aut_\cX(\xi_K)$, there is a specializing modification  $\xi'\co \Spec(R')\to \cX$ of $\xi$ and some $g \in \Aut_\cX(\xi')$ such that $g|_{K'}$ lies in $H$; \label{I:valuative_C}
\item $\cX$ has unpunctured inertia. \label{I:unpunctured}
\end{enumerate}
Then we have (1) $\Rightarrow$ (2) and (3) $\Rightarrow$ (4). 
If $\cX$ is locally reductive and $\Theta$-reductive then (2) $\Rightarrow$ (3).
\end{prop}

Given a locally reductive algebraic stack $\cX$, we've already proven the sufficiency of the conditions in 
\Cref{T:existence-general}, i.e., if $\cX$ is $\Theta$-reductive with unpunctured inertia, then $\cX$ admits an adequate moduli space.  We now establish a converse.

\begin{thm}\label{T:main_unpunctured}
Let $\cX$ be an algebraic stack locally of finite type over a quasi-separated and locally noetherian algebraic space $S$. If $\cX$ is locally reductive with separated diagonal and either
\begin{enumerate}
	\item \label{T:main_unpunctured1} 
		$\cX$ has an adequate moduli space, or
	\item \label{T:main_unpunctured2} 
		$\cX$ is \textsf{S}-complete.
\end{enumerate}
then $\cX$ satisfies all the conditions of \Cref{prop:crit-unpunctured}, in particular $\cX$ has unpunctured inertia.
\end{thm}

{We will prove \Cref{T:main_unpunctured} shortly but we first  detail some of its consequences.  We finish the proofs of \Cref{T:existence-general} and \Cref{T:A}.  We also establish the following generalization of \Cref{T:A}.  While \Cref{T:A} had a characteristic 0 hypothesis, the following theorem is characteristic independent and as in \Cref{T:A} only requires the stabilizers to be affine (rather than that the diagonal is affine) as long as the closed points have linearly reductive stabilizers.
}

\begin{thm}\label{T:existence-separated}
	Let $\cX$ be an algebraic stack of finite presentation over a quasi-separated and locally noetherian algebraic space $S$, with affine stabilizers and separated diagonal. Then $\cX$ admits a good moduli space $X$ separated over $S$ if and only if
		\begin{enumerate}
		\item  every closed point of $\cX$ has linearly reductive stabilizer;
		\item $\cX  \to S$ is $\Theta$-reductive  (\Cref{D:theta-reductive}); and
		\item $\cX \to S$ is \Scomplete (\Cref{D:S-complete}).
	\end{enumerate} 
		If $\cX$ is locally reductive (\Cref{D:locally-reductive}) and has affine diagonal, then $\cX$ admits an adequate moduli space $X$ separated over $S$ if and only if (2) and (3) hold. In both cases, if $S$ is locally excellent and $\cX \to S$ has affine diagonal, it suffices to check the filling conditions of \Cref{D:theta-reductive} and \Cref{D:S-complete} only for DVRs that are essentially finite type over $S$.

Furthermore, in both cases $X \to S$ is proper if and only if $\cX \to S$ satisfies the existence part of the valuative criterion for properness.
\end{thm}

\begin{remark}
This theorem does not hold without the quasi-compactness assumption on $\cX$, as was incorrectly stated in the first version of this article.  For a counterexample, given an infinite field $k$ and a $k$-point $p \in \bA^1_{k}$, consider the quotient stack $\cX_p = [\Spec(k[x,s,t]/(x-p-st))/\bG_m]$ where $(x,s,t)$ have weights $(0,1,-1)$.  There is a good moduli space $\cX_p \to \bA^1_{k}$ defined by $(x,s,t) \mapsto x$ which is an isomorphism away from $p$ and admits two sections $\bA^1_k \to \cX_p$ each which is an open immersion. For an infinite sequence of distinct $k$-points $p_1, p_2, \ldots$, we set $\cX_1 = \cX_{p_1}$ and inductively define $\cX_n = \cX_{n-1} \bigcup_{\bA^1_k} \cX_{p_n}$.  The union $\cX_{\infty} = \bigcup_n \cX_n$  is locally reductive, $\Theta$-reductive and \textsf{S}-complete but is not quasi-compact and does not admit a good moduli space.
\end{remark}

\begin{proof}[Proof of \Cref{T:existence-general}]
		The sufficiency of the conditions was established in \S \ref{S:existence-proof}.
		Conversely suppose that $\cX$ admits a good moduli space $\pi \co \cX \to X$. Then the closed points of $\cX$ have linearly reductive stabilizers. If $\cX \to X$ is an adequate moduli space (e.g., a good moduli space), \Cref{P:theta-reductive}\eqref{P:theta-reductive-ams} implies that $\cX$ is $\Theta$-reductive and \Cref{T:main_unpunctured} implies that $\cX$ has unpunctured inertia.
\end{proof}

\begin{proof}[Proof of \Cref{T:existence-separated}]
By \Cref{T:main_unpunctured} \Scomplete ness implies that $\cX$ has unpunctured inertia. If $\cX$ has affine diagonal over $S$, then the existence of a good (resp. adequate) moduli space follows from \Cref{T:existence-general}. 

To handle the weaker condition of affine stabilizers when the closed points have linearly reductive stabilizers, we will follow the proof of existence of \Cref{T:existence-general}, bootstrapping on the affine diagonal case.

Since $\cX$ is locally reductive and has separated diagonal, around any closed point $x \in |\cX|$ there exists a representable and \'etale quotient presentation $f \co (\cX_1, x_1) \to (\cX,x)$ with $\cX_1=[\Spec(A) / \GL_N]$. Since both $\cX$ and $\cX_1$ are \Scomplete, so is $f$.  As $f$ is also representable, it follows from \Cref{P:quasi-finite-S-complete} that $f$ is separated. Since $f \co \cX_1 \to \cX$ is separated, representable and \'etale,  \Cref{P:quotient-presentations} implies that after replacing $\cX_1$ with an open neighborhood of $x_1$, we may assume that $f$ is $\Theta$-surjective and $f$ induces an isomorphism $I_{\cX_1} \to \cX_1 \times_{\cX} I_{\cX}$. Because $\cX$ is quasi-compact, we can can assume that $f$ is also surjective by repeating this construction at finitely many closed points and replacing $\cX_1$ with the disjoint union of the resulting local quotient presentations at each of these points.

Let $\pi_1 \co \cX_1 \to X_1 := \Spec(A^{\GL_N})$ be the adequate moduli space and define $\cX_2 = \cX_1 \times_{\cX} \cX_1$ as in Diagram \eqref{E:existence-etale-equivalence-relation}.  Since $\cX$ and $\cX_1$ are both $\Theta$-reductive and \textsf{S}-complete, the morphism  $\cX_1 \to \cX$ is $\Theta$-reductive and \textsf{S}-complete and thus so is $\cX_2$. Since $\cX_1 \to \cX$ is separated and $\cX_1$ has affine diagonal over $S$, $\cX_2$ also has affine diagonal over $S$. Since $\cX_1 \to \cX$ is $\Theta$-surjective and induces isomorphisms on stabilizer groups, it follows that the closed points of $\cX_2$ have linearly reductive stabilizers and thus $\cX_2$ is locally reductive by \Cref{T:ahr}. We may thus apply the affine diagonal case to construct an adequate moduli space $\cX_2 \to X_2$.  Finally, one argues as in the proof of sufficiency in \Cref{T:existence-general} that $X_2 \rightrightarrows X_1$ is an \'etale equivalence relation and that $X:=X_1/X_2$ is an adequate moduli space of $\cX$.

From \Cref{P:gms-S-complete}\eqref{P:gms-S-complete2} we know that \textsf{S}-completeness of $\cX \to S$ implies that $X \to S$ is separated. 

Conversely if there exists a adequate moduli space $X$ separated over $S$, then $\cX \to S$ is \textsf{S}-complete by \Cref{P:gms-S-complete}\eqref{P:gms-S-complete2} and $\Theta$-reductive by \Cref{P:theta-reductive}\eqref{P:theta-reductive-ams}.

When $S$ is locally excellent, it suffices to check the definition of $\Theta$-reductivity and {\textsf{S}}-completeness only for essentially finite type DVRs by \Cref{P:theta-reductive-essentially-finite-type-DVRs} and \Cref{P:S-complete-essentially-finite-type-DVRs}, because $\cX$ is locally reductive. The criterion for properness was explained in \Cref{P:gms-S-complete}\eqref{P:gms-S-complete3}.
\end{proof}

\begin{remark}
With the hypotheses of \Cref{T:existence-separated}, \cite[Cor.~13.11]{ahr2} implies that if $\cX$ admits a separated good moduli space, then $\cX$ necessarily has affine diagonal.
\end{remark}

\begin{proof}[Proof of \Cref{T:A}] We saw in \Cref{P:SCompleteReductive} that \Scomplete ness implies reductivity of stabilizer groups at closed points and these are linearly reductive in characteristic $0$. The result thus follows from \Cref{T:existence-separated}.
\end{proof}


We provide a few remarks and examples illustrating the valuative criteria \eqref{I:valuative_A}-\eqref{I:unpunctured} in \Cref{prop:crit-unpunctured}.

\begin{remark}
In the case that $\cX=[X/G]$ where $G$ is a reductive group acting on an affine scheme $X$ of finite type over a field $k$, then the quotient stack $\cX$ satisfies valuative criterion \eqref{I:valuative_A} if and only if for every map $\xi \co \Spec(R) \to X$ and $g \in G_{\xi_K} \subset G(K)$ of finite order, there exists after an extension $R \subset R'$  (with $K' = \Frac(R')$) an element $h \in G(K')$ such that $h \cdot \xi_{K'}$ and $h^{-1}g|_{K'}h$ both extend to $R'$-points. We are not aware of a completely elementary proof of this fact, despite the purely representation-theoretic nature of this property. This is the most challenging part of the proof of \Cref{T:main_unpunctured}.
\end{remark}

\begin{ex} \label{E:ex2}
To illustrate the subtlety of valuative criterion \eqref{I:valuative_A}, let us exhibit in the context of \Cref{E:ex1} a map from a DVR to $[\bA^2 / G]$ (where $G=\bG_m \rtimes (\bZ/2)$) where performing an extension and elementary modification allows a generic automorphism to extend.   
Let $R = k[\![z]\!]$ and $K = k(\!(z)\!)$.  
Consider $\xi \co \Spec(R) \to \bA^2$ via $z \mapsto (z^2,z)$. 
Then $g=(z,-1) \in G(K)$ stabilizes $\xi_K$ but does not extend to $G(R)$.  
Consider the degree 2 ramified extension $R \to R'$ with $R'=k[\![\sqrt{z}]\!]$ and $K' = k(\!(\sqrt{z})\!)$, 
and define $\xi' \co \Spec(R') \to X$ by $\sqrt{z} \mapsto ((\sqrt{z})^3,(\sqrt{z})^3)$. 
Over the generic point, $\xi'$ is isomorphic as a point in $[\bA^2/G]$ to the restriction $\xi_{K'}$, because {$\xi'_{K'} = (\sqrt{z},-1) \cdot \xi_{K'}$.} Under this isomorphism our generic automorphism $g$ becomes 
$g' =(\sqrt{z},-1)^{-1} \cdot g|_{K'}\cdot (\sqrt{z},-1)$.
which extends to an $R'$-point. 
\end{ex}

\begin{rem}
	The valuative criterion \eqref{I:valuative_B} does not imply the valuative criterion \eqref{I:valuative_A} without additional hypotheses. Consider the group $\bG_m \ltimes \bG_a$ given coordinates $(z,y)$ and the product rule $(z_1,y_1) \cdot (z_2,y_2) = (z_1 z_2, z_2 y_1 + y_2)$, and let $G \subset (\bG_m \ltimes \bG_a) \times \bA^1_t$ be the hypersurface cut out by the equation $ty = 1-z$. Then $G$ is in fact a smooth subgroup scheme over $\bA^1$ whose fiber over $0$ is $\bG_a$ and whose fiber everywhere else is $\bG_m$.
	
	Let $\cX = B_{\bA^1}G$ and consider the map $\xi \co \Spec(k[\![t]\!]) \to \cX$ which is just the completion of the canonical map $\bA^1_t \to \cX$ at the origin. Then all modifications of $\xi$ agree after composing with the projection $\cX \to \bA^1_t$, so after an extension of DVRs the automorphism group of $\xi$ will be isomorphic to $G_{k[\![t]\!]}$. There is a generic automorphism of $\xi$ given by the formula $(\alpha,(1-\alpha) / t)$, where $\alpha$ is a non-identity $n^{th}$ root of unity. This automorphism does not extend to $0$, and the generic automorphism group is abelian and hence acts trivially on itself by conjugation. It follows that no extension and modification of $\xi$ will allow this generic automorphism to extend either.
	\end{rem}

\subsection{Proof of \texorpdfstring{\Cref{prop:crit-unpunctured}}{Proposition}}

\begin{proof}[Proof of \Cref{prop:crit-unpunctured}]
To show the implication \eqref{I:valuative_A} $\Rightarrow$ \eqref{I:valuative_B}, i.e., that components of automorphism groups extend after passing to a nearby modification, if the same holds for automorphisms of finite order,  it suffices to show that every connected component of the $K$-group $\Aut_{\cX}(\xi_K)$ contains a finite type point of finite order. Let $g \in \Aut_{\cX}(\xi_K)$ be a finite type point. After a finite field extension we can decompose $g = g_s g_u$ under the Jordan decomposition, where $g_s$ is semisimple and $g_u$ is unipotent. Now consider the reduced Zariski closed $K$-subgroup $H \subset \Aut_{\cX}(\xi_K)$ generated by $g_s$. Because $g_s$ is semisimple, $H$ is a diagonalizable $K$-group and hence every component of $H$ contains an element of finite order. We may thus replace $g_s$ with a finite order element in the same connected component of $\Aut_{\cX}(\xi_K)$ which still commutes with $g_u$. If $\op{char}(K)>0$, then $g_u$ has finite order and we are finished. If $\op{char}(K)=0$, then $g_u$ lies in the identity component of $G$, so $g$ lies on the same component as the finite order element $g_s$.

Let us now prove the implication \eqref{I:valuative_C} $\Rightarrow$ \eqref{I:unpunctured}, i.e., that the specialization of components of automorphism groups after specializing modifications implies unpunctured inertia. Let $x \in |\cX|$ be a closed point, $p\co (U,u) \to (\cX,x)$ be a {formally} smooth local morphism {where $U$ is a localization of a smooth neighborhood of $x$} and $H \subset \Aut_{\cX}(p)$ a connected component. The image of the projection $H \to U$ is a constructible set whose closure contains $u$. It follows that we can find an essentially finite type DVR $R$ and a map $\Spec(R) \to U$ whose special point maps to $u$ and whose generic point lies in the image of $H \to U$. After an extension of the DVR $R$, we may assume that the generic point $\Spec(K) \to U$ lifts to $H$, and that the connected component $H' \subset H|_{\Spec(K)}$ containing this lift is geometrically connected. The hypotheses {\eqref{I:valuative_C} implies} that, after possibly further extending $R$, there exists a modification $\xi' \co \Spec(R) \to \cX$ of $\xi$ such that the closure of $H'$ in $\Aut_{\cX}(\xi)$ meets the fiber over $0 \in \Spec(R)$ and $0 \in \Spec(R)$ still maps to $u$. By construction $H'$ maps to $H$, which implies that $H \subset \Aut_{\cX}(p)$ meets the fiber over $u$.

If $\cX$ is locally reductive and $\Theta$-reductive, the implication \eqref{I:valuative_B} $\Rightarrow$ \eqref{I:valuative_C} follows from the following lemma.
\end{proof}

\begin{lem} \label{L:specializing_modification}  Let $\cX$ be a locally reductive and $\Theta$-reductive algebraic stack locally of finite type over a quasi-separated and locally noetherian algebraic space $S$.
	Let $R$ be a DVR and let $\xi \co\Spec(R) \to \cX$ be a morphism. If $\xi'$ is a nearby modification of $\xi$ and $g \in \Aut_\cX(\xi')$, then there is a finite extension of DVRs $R' / R$ with fraction field $K'$ and a modification $\xi'' \co \Spec(R') \to \cX$ of $\xi$ such that $g|_{K'}$ extends to an automorphism of $\xi''$ and $\xi''(0)$ is a specialization of $\xi(0)$.
\end{lem}

\begin{proof}
Let us first reduce to the situation that $\cX=[\Spec(A)/\GL_N]$ is a quotient stack:   Let $\kappa$ be the residue field of $R$ and let $\cZ \subset \cX_{\kappa} = \cX \times_S \Spec(\kappa)$ be the closure of the point $p:= \xi'(0) \in \cX_{\kappa}$. By \Cref{L:unique_closed_point} we know that $\cZ$ has a unique closed point $z \in |\cZ|$, and in particular $z$ is a specialization of both $p$ and $\xi(0)$ because $\xi'$ is a nearby modification of $\xi$. If necessary we pass to a finite extension of $R$ so that we may assume that $z \in \cZ(\kappa)$ as well. As $\cX$ is locally reductive and $\Theta$-reductive \Cref{P:GMS_one_point} implies that $\cZ \simeq [\Spec(A) / G_{z}]$ for some affine $G_{z}$-scheme $\Spec(A)$. Embedding $G_z \subset \GL_{N,\kappa}$ for some $N$, we may replace $G_z$ with $\GL_N$ and $\Spec(A)$ with the affine scheme $\GL_N \times_{G_z} \Spec(A)$.

Kempf's theorem \cite[Theorem 4.2]{kempf1978instability} implies that after passing to a finite purely inseparable extension of $\kappa$, which can be induced by a suitable finite extension of DVRs, 
there is a canonical filtration $f\co \Theta_{\kappa} \to [\Spec(A)/\GL_N]$ with an isomorphism $f(1) \simeq p$ such that $f(0) = z$. The fact that $f$ is canonical implies that any automorphism of $p = f(1)$ extends to an automorphism of the map $f$. In particular the restriction of $g \in \Aut_\cX(\xi)$ to $p = \xi'(0)$ extends uniquely to an automorphism of $f$ which we also denote $g$.

We now apply the strange gluing lemma (\Cref{C:strange_gluing}), which states that after composing $f$ with a suitable ramified cover $(-)^n \co \Theta_\kappa \to \Theta_\kappa$, the data of the map $\xi' \co \Spec(R) \to \cX$ and the filtration $f \co \Theta_\kappa \to \cX$, comes from a \emph{unique} map $\gamma\co \oST_{R} \to \cX$, where $f$ is the restriction of $\gamma$ to the locus $\{s=0\}$ and $\xi'$ is the restriction of $\gamma$ to the locus $\{t \neq 0\}$. The uniqueness of this extension guarantees that the automorphism $g$ of $\xi'$ and $f$ extends uniquely to an automorphism of $\gamma$, which we again denote $g$. Finally we construct our modification $\xi''$ as the composition
\[
\xi'' \co  \Spec(R[\sqrt{\pi}]) \to \oST_{R} \xrightarrow{\gamma} \cX,
\]
where the first map is given in $(s,t,\pi)$ coordinates by $(\sqrt{\pi},\sqrt{\pi},\pi)$, which maps the special point of $\Spec(R[\sqrt{\pi}])$ to the point $\{s=t=\pi=0\}$ of $\oST_R$. By construction the automorphism $g$ of $\gamma$ restricts to an automorphism of $\xi''$ extending $g|_{K[\sqrt{\pi}]}$, and the special point $\xi''(0)$ maps to the closed point $z$ of $\cZ$, which is a specialization of $\xi(0)$.
\end{proof}

\subsection{Proof of \texorpdfstring{\Cref{T:main_unpunctured}\eqref{T:main_unpunctured1}}{Theorem}}
It will be convenient to introduce a variant of the valuative criterion \eqref{I:valuative_A} in \Cref{prop:crit-unpunctured} for an algebraic stack $\cX$:
\begin{enumerate}[label=(\arabic*$'$)]
\item For any DVR $R$ with fraction field $K$, any morphism $\xi \co \Spec(R) \to \cX$, and any $g \in \Aut_\cX(\xi_K)$ of finite order, there is a modification $\xi' \co \Spec(R') \to \cX$ of $\xi$ such that $g|_{K'}$ extends to an automorphism of $\xi'$. \label{I:valuative_A_prime}
\end{enumerate}
Unlike in the valuative criterion \eqref{I:valuative_A}, criterion \ref{I:valuative_A_prime} does not require that the modification $\xi'$ is a nearby modification. Criterion \ref{I:valuative_A_prime} is not a sufficient condition for $\cX$ to have unpunctured inertia even when $\cX$ is locally reductive and $\Theta$-reductive, but it has useful formal properties. Note in particular that in an  {\textsf{S}}-complete stack, any modification is an elementary modification and in particular a nearby modification so condition \ref{I:valuative_A_prime} is equivalent to condition \eqref{I:valuative_A} of \Cref{prop:crit-unpunctured}.

\begin{prop} \label{P:proper-over-affine-mod-G}
Let $G$ be a geometrically reductive and \'etale-locally embeddable group scheme (e.g. reductive) over an algebraic space $S$ and let $W \to S$ be an affine morphism of finite type with an action of $G$. Then $[W/G]$ satisfies the criterion \ref{I:valuative_A_prime}. In particular, if in addition $S$ is separated, then $[W/G]$ satisfies condition \eqref{I:valuative_A} of \Cref{prop:crit-unpunctured}.
\end{prop}

We will prove \Cref{P:proper-over-affine-mod-G} at the end of this subsection, after establishing some preliminary results.

\begin{lem} \label{L:BGL_val_criteria}
The stack $B_{\bZ}\GL_N$ satisfies the condition \eqref{I:valuative_A}.
\end{lem}
\begin{proof}
In this case, because every vector bundle on $\Spec(R)$ is trivializable, the condition is equivalent to the claim that every element $g \in \GL_N(K)$ is conjugate to an element of $\GL_N(R)$ after passing to an extension of the DVR $R$. After an extension of $R$ we may conjugate $g$ to its Jordan canonical form. The fact that $g$ has finite order implies that the diagonal entries of the resulting matrix are roots of unity. Because the group of $k^{th}$ roots of unity is a finite group scheme, the entries of the Jordan canonical form must lie in $R$.
\end{proof}

\begin{lem} \label{L:proper_val_criteria}
Let $p\co\cX \to \cY$ be a proper representable morphism of noetherian stacks. If $\cY$ satisfies the valuative criterion \ref{I:valuative_A_prime}, then so does $\cX$.
\end{lem}
\begin{proof}
Since $p$ is representable and separated, for any morphism $\xi \co \Spec(R) \to \cX$ from a DVR, we have a closed immersion $\Aut_{\cX}(\xi) \hookrightarrow \Aut_{\cY}(p \circ \xi)$ of group schemes over $\Spec(R)$. Furthermore, because $p$ is proper, any modification of $p \circ \xi$ lifts uniquely to a modification of $\xi$. Therefore, given a generic automorphism of $\xi$, we may pass to an extension $R'/R$ and modify $p \circ \xi|_{R'}$ so that this generic automorphism extends, and then this lifts uniquely to a modification of $\xi|_{R'}$ such that the given generic automorphism extends.
\end{proof}

\begin{proof}[Proof of \Cref{P:proper-over-affine-mod-G}]
It suffices to show that $[\Spec(A) / G]$ satisfies the criterion \ref{I:valuative_A_prime}, where $G$ and $\Spec(A)$ are defined over a DVR R and with $A$ finitely generated over $R$. After passing to a finite extension of $K = \Frac(R)$ we may assume that $G$ embeds as a closed subgroup $G \hookrightarrow \GL_{N,R}$ for some $N$. We may then replace $G$ with $\GL_{N,R}$ and replace $\Spec(A)$ with $\GL_{N,R} \times^G \Spec(A)$, which will again be affine because $G$ is geometrically reductive. Furthermore we can assume that $A$ is reduced, because we are only considering maps from reduced schemes. So it suffices to prove the claim for $[\Spec(A)/\GL_{N,R}]$ for a reduced $R$-algebra $A$ of finite type.

Now consider a $\GL_N$-scheme $X$ which is reduced and projective over $\Spec(A^{\GL_N})$ such that $X$ contains $\Spec(A)$ as a dense $\GL_N$-equivariant open subscheme and the complement $X \smallsetminus \Spec(A)$ is the support of an ample $\GL_N$-invariant Cartier divisor $E$. The construction in \cite[Lem.~6.1]{Teleman} for smooth schemes in characteristic $0$ works here as well: We simply choose a closed $G$-equivariant embedding $\Spec(A) \hookrightarrow \bA_{A^{\GL_N}}(\cE)$ for some locally free $\GL_n$-module over $A^{\GL_N}$, and then let $X$ be the closure of $\Spec(A)$ in $\bP_{A^{\GL_N}}(\cE \oplus \cO)$. Thus we have a diagram:
\[
\xymatrix{\Spec(A) \ar@{^{(}->}[r] \ar[dr] & X \ar[d] \\ & \Spec(A^{\GL_N})}
\]

We claim that $\Spec(A)$ is precisely the semistable locus of $X$ with respect to $\cO_X(E)$ in the sense of \cite{Seshadri}. Indeed the tautological invariant section $s \co \cO_X \to \cO_X(E)$ which restricts to an isomorphism over $\Spec(A)$ shows that $\Spec(A) \subset X^{\ss}$. Conversely $s^n$ gives an isomorphism of $A^{\GL_N}$-modules for all $n>0$,$$\Gamma(\Spec(A),\cO_X(nE))^{\GL_N} \simeq A^{\GL_N}.$$Under this isomorphism any invariant global section $f \in \Gamma(X,\cO_X(nE))^{\GL_N}$, after restriction to the dense open subset $\Spec(A)$, agrees with a section of the form $g s^n$, where $g$ is the pullback of a function under the map $X \to \Spec(A^{\GL_N})$. It follows that $f = g \cdot s^n$ because $X$ is reduced. This shows that $X^{\ss} \subset \Spec(A)$.

Now \Cref{L:BGL_val_criteria} implies that the criterion \ref{I:valuative_A_prime} holds for $\Spec(A^{\GL_N}) \times B\GL_N$, and hence \Cref{L:proper_val_criteria} implies that the criterion holds for $[X/\GL_N]$. So in order to establish the criterion for $[\Spec(A)/\GL_N]$, it suffices to show that given a point $\xi \co \Spec(R) \to [X/\GL_N]$ along with a finite order automorphism $g$ of $\xi$, if $\xi_K$ lies in the open substack $[\Spec(A)/\GL_N]$, then after passing to an extension of $R$ one can modify the pair $(\xi,g)$ at the special point of $\Spec(R)$ so that the image of $\xi$ lies in $[\Spec(A)/\GL_N]$. Note that the stabilizer group $X$-scheme $\Stab_{\GL_N}(X) \subset X \times \GL_N$ is equivariant for the $\GL_N$ action which acts by the given action on $X$ and by conjugation on the $\GL_N$ factor. It suffices to show that given an $R$-point of $\Stab_{\GL_N}(X)$ whose generic point lies in $\Spec(A) \times \GL_N$, after passing to an extension of $R$ there is a modification of the resulting map $\xi \co \Spec(R) \to [\Stab_{\GL_N}(X) / \GL_N]$ whose image lies in $[\Stab_{\GL_N}(\Spec(A)) / \GL_N] = [ \Stab_{\GL_N}(X) \cap (\Spec(A) \times \GL_N) / \GL_N]$.

Note that $\Stab_{\GL_N}(X)$ is projective over $\Spec(A^{\GL_N}) \times \GL_N$. We claim that the semistable locus of $\Stab_{\GL_N}(X)$ for the action of $\GL_N$ with respect to the pullback of $\cO_X(E)$ is precisely $\Stab_{\GL_N}(\Spec(A))$. Indeed, this follows from the Hilbert--Mumford criterion \cite{Seshadri}. Any destabilizing cocharacter for $(x,g) \in \Stab_{\GL_N}(X)$ is destabilizing for $x \in X$. Conversely, for every point $(x,g) \in \Stab_{\GL_N}(X)$ whose underlying point $x \in X$ is unstable, Kempf's theorem on the existence of canonical destabilizing flags \cite{kempf1978instability} implies that there is a destabilizing cocharacter $\lambda$ for $x$ which commutes with $\Stab_{\GL_N}(x)$, and this $\lambda$ defines a destabilizing cocharacter for the point $(x,g)$. Applying GIT in the projective over affine situation to the projective morphism $\Stab_{\GL_N}(X) \to\Spec(A^{\GL_N}) \times \GL_N$ we therefore find a proper adequate moduli space $\Stab_{\GL_N}(X)^{\ss}/\!/\GL_N \to \Spec(A^{\GL_N}) \times \GL_N/\!/\GL_N$. Thus any morphism $\Spec(R) \to [\Stab_{\GL_N}(X) / \GL_N]  \to \Spec(A^{\GL_N}) \times \GL_N/\!/\GL_N$ for which the generic point maps into the semistable part admits a modification that maps to $[\Stab_{\GL_N}(\Spec(A)) / \GL_N]$.
\end{proof}

\begin{remark}
By appealing to \Cref{T:valuative-criterion-finite}, the proof in fact shows that a slightly stronger version of \ref{I:valuative_A_prime} holds in which the extension $K'/K$ of fraction fields is {\it finite}.  It follows that the statements of \Cref{prop:crit-unpunctured} and \Cref{T:main_unpunctured} remain true after replacing the valuative criteria \eqref{I:valuative_A} and \eqref{I:valuative_B} with the stronger condition where the extension $K'/K$ is required to be finite.
\end{remark}

\begin{proof}[Proof of \Cref{T:main_unpunctured}\eqref{T:main_unpunctured1}]
If $\cX$ has an adequate moduli space, then $\cX$ is necessarily $\Theta$-reductive (\Cref{P:theta-reductive}).  \Cref{prop:crit-unpunctured} establishes the implications \eqref{I:valuative_A} $\Rightarrow$ \eqref{I:valuative_B} $\Rightarrow$  \eqref{I:valuative_C} $\Rightarrow$  \eqref{I:unpunctured} between the valuative criteria.  To verify \eqref{I:valuative_A}, let $\xi \co \Spec(R) \to \cX$ be a morphism, and let $g$ be an automorphism of $\xi_K$ of finite order. Then we may choose an \'etale map $U \to X$ whose image contains the image of $\Spec(R)$ and such that $\cU := \cX \times_X U \simeq [\Spec(A)/\GL_N]$ (\Cref{P:local-on-ams}).  After replacing $R$ with an extension of DVRs we may assume that $\xi'$ lifts to a map $\xi' \co \Spec(R') \to \cU$. Furthermore the map $\cU \to \cX$ is inertia preserving in the sense that $I_\cU \simeq I_\cX \times_\cX \cU$, which implies that $g$ lifts to a finite order generic automorphism $g'$ of $\xi'_{K'}$. By \Cref{P:proper-over-affine-mod-G} the stack $\cU$ satisfies condition \eqref{I:valuative_A} of \Cref{prop:crit-unpunctured}. This provides a nearby modification of the map $\xi'$ for which $g'$ extends, and we can compose this with the map $\cU \to \cX$ to get a nearby modification of the original map for which $g|_{K'}$ extends.
\end{proof}

\subsection{The proof of \texorpdfstring{\Cref{T:main_unpunctured}\eqref{T:main_unpunctured2}}{Theorem}}
Let $\cX$ be a noetherian algebraic stack with affine stabilizers, let $p \co \cY \to \cX$ be an \'etale map with $\cY \cong [\Spec(A)/G]$, and let $R$ be a complete DVR with fraction field $K$ and residue field $\kappa$. Let $x \in |\cX|$ be a closed point such that $p$ induces an isomorphism $p^{-1}(\cG_x) \simeq \cG_x$, where $\cG_x$ denotes the residual gerbe of $x$.

\begin{lem}
The functor $\uMap(\Spec(R),\cY) \to \uMap(\Spec(R),\cX)$ defined by composition with $p$ induces an equivalence between the full subgroupoids of maps taking the special point of $\Spec(R)$ to $p^{-1}(x)$ and $x$ respectively. The same is true for the functor $\uMap(\oST_R,\cY) \to \uMap(\oST_R,\cX)$ and the subgroupoid taking the point $0 \in \oST_R(\kappa)$ to $p^{-1}(x)$ and $x$ respectively.
\end{lem}

\begin{proof}
The map $p$ is \'etale and induces an equivalence between the residual gerbe of $x \in |\cX|$ and $p^{-1}(x) \in |\cY|$. It therefore induces an equivalence between the $n^{th}$ order neighborhoods of these residual gerbes, so any map $\Spec(R) \to \cX$ mapping $0$ to $x$ lifts uniquely along $p$ over any nilpotent thickening of $0 \in \Spec(R)$. The result then follows from Tannaka duality and the fact that $\Spec(R)$ is coherently complete along its special point, so a compatible family of lifts over nilpotent thickenings of $0$ corresponds to a unique lift of the map $\Spec(R) \to \cX$ along $p$. The same argument applies to $\oST_R$, which is coherently complete along the inclusion $(B\bG_m)_\kappa \hookrightarrow \oST_R$ at the point $0$.
\end{proof}

Now let $\xi' \co \Spec(R) \to \cY$ be an $R$-point mapping the special point to $p^{-1}(x)$, and let $\xi = p \circ \xi'$.

\begin{lem} \label{L:stabilizer_preserving}
If $\cX$ and $\cY$ are {\textsf{S}}-complete, then $p$ induces an isomorphism
\[
\Aut_\cY(\xi'_K) \xrightarrow{\cong} \Aut_{\cX}(\xi_K).
\]
\end{lem}

\begin{proof}{For any stack $\cX$ and $\xi\in \cX(R)$ we can express $\Aut_{\cX}(\xi_K)$ in terms of $\oST_R \smallsetminus 0=\Spec(R) \cup_{\Spec(K)} \Spec(R)$, because the gluing condition for stacks identifies $\cX(\oST_R \smallsetminus 0)$ with the category of triples $(\xi_1,\xi_2,\phi_K)$ where $\xi_1,\xi_2\in \cX(R)$ and $\phi_K\in \Isom_{\cX}(\xi_{1,K},\xi_{2,K})$. Thus we have an isomorphism:
		$$
	\Aut_{\cX}(\xi_K) \simeq \{ (f,\phi_1,\phi_2) \,|\, f=(\xi_1,\xi_2,\phi_K) \in \cX(\oST_R \smallsetminus 0), \phi_1\colon \xi_1 \map{\cong}\xi, \phi_2\colon \xi_2 \map{\cong} \xi_i \}.
		$$ }

Because both $\cY$ and $\cX$ are  {\textsf{S}}-complete, restriction gives an equivalence of groupoids $\Map(\oST_R, - ) \to \Map(\oST_R \smallsetminus 0,-)$ for both stacks. It follows that
\begin{equation} \label{E:automorphisms}
\Aut_\cX(\xi_K) \simeq \{ \text{maps } \oST_R \to \cX + \text{ equivalences } f_1 \simeq \xi \simeq f_2\}
\end{equation}
and likewise for $\cY$. If $\xi$ maps the special point of $\Spec(R)$ to $x \in \cX$, then any map $f \co \oST_R \to \cX$ which admits an isomorphism $f_1 \simeq \xi$ must also map $(0,0)$ to $x$, because $x$ is closed. Now the previous lemma implies that $p$ induces a bijection of sets on the right hand side of \eqref{E:automorphisms} for $\cX$ and $\cY$, and thus also on the left hand side.
\end{proof}

\begin{proof} [Proof of  \Cref{T:main_unpunctured}\eqref{T:main_unpunctured2}]
We first verify criterion \eqref{I:valuative_A} of \Cref{prop:crit-unpunctured}. 
Consider a map from a DVR $\xi \co \Spec(R) \to \cX$. Let $x \in |\cX|$ be a closed point in the closure of $\xi(0)$, and let $p : [\Spec(A)/\GL_N] \to \cX$ be an \'etale quotient presentation around $x$. Then after an extension of DVRs $R'/R$ we may lift $\xi$ to a map $\xi' : \Spec(R') \to [\Spec(A)/\GL_N]$. Now  \Cref{L:specializing_modification} allows one to construct a modification of $\xi'$, after replacing $R'$ with a further finite extension, which maps the special point to the closed point $p^{-1}(x) \in [\Spec(A)/\GL_N]$. It follows from \Cref{L:stabilizer_preserving} that the map
\[
\Aut_{[\Spec(A)/\GL_N]} (\xi'_{K'}) \to \Aut_{\cX}(p\circ \xi'_{K'})
\]
is an isomorphism of $K'$-groups. In particular given a finite order element $g \in \Aut_{\cX}(\xi_K)$, one may lift this to $\xi'$ after replacing $R'$ with a further extension. We know that the criterion \eqref{I:valuative_A} of \Cref{T:main_unpunctured} holds for $[\Spec(A)/\GL_N]$ by \Cref{P:proper-over-affine-mod-G}, and after replacing $R'$ with a further extension this produces a nearby modification for which $g|_{K'}$ extends. Composing with $p$ gives a nearby modification of the original map $\xi$ for which $g$ extends.  

The same argument shows that $\cX$ satisfies the criterion \eqref{I:valuative_C}. It follows from \Cref{prop:crit-unpunctured} that \eqref{I:valuative_A} $\Rightarrow$ \eqref{I:valuative_B} and \eqref{I:valuative_C} $\Rightarrow$ \eqref{I:unpunctured}.
\end{proof}


\section{Semistable reduction and \texorpdfstring{$\Theta$}{Theta}-stability}\label{S:semistable_reduction}

In this section we explain how completeness properties of stacks induce similar properties of the substack of semistable objects, if these are defined using the theory of $\Theta$-stability. Our key result is \Cref{T:Langton-Algorithm} {which} is inspired by Langton's algorithm for semistable reduction for families of torsion-free sheaves on a projective variety. Recall from \Cref{R:langton} that this  algorithm starts with a family of bundles parametrized by a DVR $R$ such that the generic fiber is semistable and the special fiber is unstable, and then applies elementary modifications to arrive at a semistable family. Surprisingly, it turns out that his construction admits an analog that relies only on the geometry of the algebraic stack representing the moduli problem, not on the particular type of objects classified by the moduli problem. The structure we will need is that of a $\Theta$-stratification from \cite[Def.~2.1.1]{hlinstability} {which} formalizes the notion of canonical filtrations in geometric terms. 

\begin{defn}\label{D:theta-stratum}
	Let $\cX$ be an algebraic stack locally of finite type over a noetherian algebraic space $S$. 
\begin{enumerate}
	\item 	A \emph{$\Theta$-stratum in $\cX$} consists of a union of connected components $\cZ^+ \subset \Filt(\cX)=\uMap_S(\Theta,\cX)$ such that $\ev_1 \co \cZ^+ \to \cX$ is a closed immersion.	
	\item  A \emph{$\Theta$-stratification of $\cX$} indexed by a totally ordered set $\Gamma$ is a cover of $\cX$ by open substacks $\cX_{\leq c}$ for $c \in \Gamma$ such that $\cX_{\leq c} \subset \cX_{\leq c'}$ for $c < c'$, along with a $\Theta$-stratum $\cZ^+_c \subset \Filt(\cX_{\leq c})=\uMap_S(\Theta, \cX_{\leq c})$ in each $\cX_{\leq c}$ whose complement is $\bigcup_{c'<c} \cX_{\leq c'} \subset \cX_{\leq c}$. We require that $\forall x \in |\cX|$ the subset $\{c \in \Gamma| x \in \cX_{\leq c}\}$ has a minimal element. We assume for convenience that $\Gamma$ has a minimal element $0 \in \Gamma$.
	\item  We say that a $\Theta$-stratification is \emph{well-ordered} if for any point $x \in |\cX|$, the totally ordered set $\{c \in \Gamma | \ev_1(\cS_c) \cap \overline{\{x\}} \neq \emptyset\}$ is well-ordered.
\end{enumerate}		
\end{defn}	 
\begin{rem}
It will be convenient for us to identify a $\Theta$-stratum $\cZ^+$ with the closed substack it defines on $\cX$, i.e., we will sometimes say that a closed substack $\cZ^+\subset \cX$ is a $\Theta$-stratum, if there exist a union of  connected components $\cZ^{\prime+} \subset \uMap_S(\Theta,\cX)$ such that $\ev_1 \co \cZ^{\prime+} \to \cZ^+ \subset \cX$ is an isomorphism.
\end{rem}


Given a $\Theta$-stratification, we {refer to the open substack} $\cX^{\ss} := \cX_{\leq 0}$ as the semistable locus. For any unstable point $x \in \cX(k) \smallsetminus \cX^{\ss}(k)$, the $\Theta$-stratification determines a canonical filtration $f\co \Theta_{k} \to \cX$ with $f(1) \simeq x$, which we refer to as the \emph{HN filtration}.

Restricting a map $f\co \Theta \to \cX$ to $B \bG_m \hookrightarrow \Theta$ defines a map $$\gr \co \Filt(\cX)=\uMap_S(\Theta,\cX) \to \uMap_S(B\bG_m,\cX) =\Grad(\cX)$$ which corresponds to ``passing to the associated graded object'' of the filtration $f$. Composition with the projection $\Theta \to B\bG_m$ defines a section $\sigma \co \uMap_S(B\bG_m,\cX) \to \uMap_S(\Theta,\cX)$ of the map $\gr$ which corresponds to the ``canonical split filtration of a graded object.'' These maps define a canonical $\bA^1$ deformation retract of $\uMap_S(\Theta,\cX)$ onto $\uMap_S(B\bG_m,\cX)$, and in particular induce bijections on connected components \cite[Lem.~1.3.8]{hlinstability}. We refer to the union of connected components $\cZ \subset \uMap_S(B\bG_m,\cX)$ corresponding to $\cZ^{+}$ as the \emph{center} of the $\Theta$-stratum $\cZ^{+}$. The result is a diagram
\[
\xymatrix{ \cZ \ar@{^(->}[r]^\sigma & \cZ^{+}\ar@/^/[l]^{\gr} \ar@{^(->}[r]^{\ev_1} & \cX.}
\]

\subsection{The semistable reduction theorem}

\begin{thm}[Langton's algorithm]\label{T:Langton-Algorithm}
Let $\cX$ be an algebraic stack locally of finite type and quasi-separated, with affine {stabilizers}, over a noetherian algebraic space $S$, and let $\cZ^{+}\hookrightarrow \cX$ be a $\Theta$-stratum. Let $R$ be a DVR with fraction field $K$ and residue field $\kappa$. 
Let $\xi_{R} \co \Spec(R) \to \cX$ be an $R$-point such that the generic point $\xi_K$ is not mapped to $\cZ^{+}$, but the special point $\xi_k$ is mapped to $\cZ^{+}$: 
	$$\xymatrix{
		\Spec(K) \ar@{^(->}[r]\ar[d]^{\xi_K} & \Spec(R) \ar[d]^{\xi_R} & \Spec(\kappa)\ar@{_(->}[l]\ar[d]^{\xi_{\kappa}} \\
		\cX\smallsetminus \cZ^{+} \ar@{^(->}[r]^{j}& \cX & \ar@{_(->}[l]_{\iota} \cZ^{+}.}$$
Then there exists an extension $R \to R'$ of DVRs with $K \to K'=\Frac(R')$ finite and an elementary modification $\xi'_{R^\prime}$ of $\xi_{R^\prime}$ such that $\xi_{R^\prime}^\prime \co \Spec(R') \to \cX$ lands in $\cX\smallsetminus \cZ^+$. 
\end{thm}

\begin{rem}
In the proof of the above result we will apply the non-local slice theorem (\Cref{T:relativeslice}) for algebraic stacks. As the proof of this result has not appeared, we give an alternative argument using \cite[Thm.~1.2]{ahr}, which requires the additional hypothesis that $S$ is the spectrum of an algebraically closed field and that for any $x \in \cX(k)$, the automorphism group $G_x$ is smooth -- this suffices, in particular, for stacks over a field of characteristic $0$. 
\end{rem}

This theorem is stated for a single stratum, but it immediately implies a version for a stack with a $\Theta$-stratification:

\begin{thm}[Semistable reduction] \label{T:semistable_reduction}
Let $\cX$ be a quasi-separated algebraic stack with affine {stabilizers} that is locally finite type over a noetherian algebraic space $S$, with a well-ordered $\Theta$-stratification. Then for any morphism $\Spec(R) \to \cX$, after an extension $R \to R^\prime$ of DVRs with $K \to K'=\Frac(R')$ finite there is a modification $\Spec(R') \to \cX$, obtained by a finite sequence of elementary modifications, whose image lies in a single stratum of $\cX$.
\end{thm}
\begin{proof}
Beginning with a map $\xi_{R} \co \Spec(R) \to \cX$ such that $\xi_{K} \in \cZ_c^+$ and $\xi_{\kappa} \in \cZ^+_{c_0}$ for $c_0>c$, we may apply \Cref{T:Langton-Algorithm} iteratively to obtain a sequence of finite extensions of $R$ and elementary modifications of $\xi$ with special point in $\cZ^+_{c_i}$ for $c_0 > c_1 > \cdots$. Each $\cZ^+_{c_i}$ meets $\overline{\xi_{K}}$, so the well-orderedness condition guarantees that this procedure terminates, and it can only terminate when $c_i = c$.
\end{proof}

\begin{rem}
In the relative situation when $\cX$ is defined over a base algebraic stack $S$, one can base change the structure of a $\Theta$-stratification along a smooth map $S' \to S$, so both \Cref{T:Langton-Algorithm} and \Cref{T:semistable_reduction} extend immediately to the case of a quasi-separated and locally noetherian base stack $S$. 
\end{rem}

\subsubsection{Langton's algorithm in the basic situation}
The main idea of the proof is to reduce to the situation where $\cX=[\Spec(A)/\bG_m]$ is  the quotient of an affine scheme by an action of $\bG_m$, $\cZ=[(\Spec A)^{\bG_m}/\bG_m]$ is the substack defined by the fixed point locus of the action and $\cZ^+=[\Spec (A/I_{+}) /\bG_m]$ is the attracting substack, where
\[
I_+ := (\bigoplus_{n>0} A_n)
\]
is the graded ideal generated by elements of positive weight {for the $\bG_m$-action}. In this basic situation the theorem will then follow from an elementary calculation. We will first explain the proof of this special case and then show how to reduce to the basic situation.

\begin{lem}\label{L:graded_affine_reduction}
In the setting of \Cref{T:Langton-Algorithm}, suppose in addition that $\cX=[\Spec(A)/\bG_m]$ and that $\cZ^+=[\Spec(A/I_+)/\bG_m]$ {is the subscheme defined by the elements of positive weight in $A$}. Then the conclusion of \Cref{T:Langton-Algorithm} holds. 
\end{lem}

\begin{proof}
	Let us denote $\widetilde{X}:=\Spec(A)$ and $\widetilde{Z}^+:=\Spec(A/I_+)$. As $\tX\to \cX$ is a $\bG_m$-torsor, we can lift $\xi$ to a map $\xi'_{R} \co \Spec(R) \to \Spec(A)$, obtaining a diagram
	$$\xymatrix{
		\Spec(K) \ar@{^(->}[r]\ar[d]^{\xi'_K} & \Spec(R)\ar[d]^{\xi'_{R}} & \Spec(\kappa)\ar@{_(->}[l]\ar[d]^{\xi'_{\kappa}} \\
		\tX\smallsetminus\tZ^+ \ar@{^(->}[r]^-{j}& \Spec(A) & \ar@{_(->}[l]_{\iota}\Spec(A/I_+).}$$
	As $\xi'_{\kappa}\in \tZ^+=\Spec(A/I_+)$ and $A/I_+$ is generated by elements of non-positive weight, the $\bG_m$-orbit of $\xi'_{\kappa}$, corresponding to a map of graded algebras $A/I_+ \to k[t^{\pm 1}]$ where $t$ has weight $-1$, extends to an equivariant morphism $\bA^1_{\kappa} \to \tZ^+$. Thus the $\bG_m$-orbits of the points $\xi'_K,\xi'_{R},\xi'_{\kappa}$ define a diagram:
	$$\xymatrix@R=10pt{
		\bG_{m,K} \ar@{^(->}[r]\ar[dd]^{f_K} & \bG_{m,R} \ar[dd]^{f_R} & \bG_{m,{\kappa}}\ar@{_(->}[l]\ar@{_(->}[d] \\
		& & \bA^1_{\kappa} \ar[d]^{\overline{f}_{\kappa}} \\
		\tX\smallsetminus \tZ^+ \ar@{^(->}[r]^-{j}& \Spec(A) & \ar@{_(->}[l]_{\iota}\Spec(A/I_+).}$$
	We know that $f_R^\#(I_+) \in \pi(R[t^{\pm 1}])$ since $f_{\kappa}^\#$ factors through $A/I_+$, and we have $K[t^{\pm 1}]\cdot f_R^\#(I_+)=K[t^{\pm 1}]$ since the image of $f_K$ does not intersect $\tZ^+$.  	
	
	Let $a_i\in I_{d_i}$ be homogeneous generators of $I_+$. Then for all $i$ we have $f_R^{\#}(a_i) =\epsilon_i \pi^{n_i} t^{-d_i}$ for some $n_i > 0$ and $\epsilon_i\in R^\times \cup \{0\}$.  As $f_R^{\#}(I_+)$ is not $0$ we can define
\[
\frac{m}{d}:=\min_i \left\{ \frac{n_i}{d_i} \left| f_R^{\#}(a_i)\neq 0 \right.\right\}
\]
and let $R^\prime:=R[\pi^{\frac{1}{d}}]$. Since $n_i-\frac{d_i m}{d} \geq 0$ for each $i$, we can write 
\[
f_R^{\#}(a_i)=\epsilon_i (\pi^{n_i-\frac{d_i m}{d}}) (\pi^{\frac{m}{d}}t^{-1})^{d_i}= \epsilon_i (\pi^{n_i-\frac{d_i m}{d}}) s^{d_i}.
\]
Since $f^\#_R$ maps elements of negative weight to $R[t]$, we have a homomorphism of graded rings
	$$f_{R'}^{\prime\#}\co A \to R^\prime[s,t]/(st-\pi^{\frac{m}{d}}) = R^\prime[t,\pi^{\frac{m}{d}}t^{-1}]\subset R^\prime[t,t^{-1}]$$ 
Furthermore, composing with the map setting $s=1$ at least one $f_{R'}^{\prime\#}(a_i)$ is not mapped to $0$ mod $\pi^{\frac{1}{d}}$, i.e., $f_{R'}^\prime|_{\{s=1\}} \co \Spec (R^\prime) \mapsto \Spec(A_{a_{i}}) \subset \tX-\tZ^+$.
The graded homomorphism $f^{\prime\#}_{R^\prime}$ defines a morphism 
	$$ [\Spec \left( R^\prime[s,t]/(st-\pi^{\frac{m}{d}}) \right)/\bG_m] \to \cX=[\tX/\bG_m].$$

As $\pi^{\frac{m}{d}}$ is not a uniformizer for $R'$, this is not quite an elementary modification. However, we can embed $R'[s,t]/(st-\pi^{m/d}) \subset R'[s^{1/m},t^{1/m}] / (s^{1/m} t^{1/m} - \pi^{1/d})$. If we regard $s^{1/m}$ and $t^{1/m}$ as having weight $1$ and $-1$ respectively, the map $\Spec(R'[s^{1/m},t^{1/m}] / (s^{1/m} t^{1/m} - \pi^{1/d})) \to \Spec(R'[s,t]/(st-\pi^{m/d}))$ is equivariant with respect to the group homomorphism $\bG_m \to \bG_m$ given in coordinates by $z\mapsto z^m$. The resulting composition
\[
\oST_{R'} \to  [\Spec \left( R^\prime[s,t]/(st-\pi^{\frac{m}{d}}) \right)/\bG_m] \to \cX
\]
is the desired modification of $\xi_{R}$.
\end{proof}	

\subsubsection{Reduction to quasi-compact stacks}

We first show that by replacing $\cX$ by a suitable open substack we may assume that $\cX$ is quasi-compact.

\begin{lem} \label{L:quasi-compact_stratum}
In the setting of \Cref{T:Langton-Algorithm}, let $\sigma \co \cZ \to \cZ^+$ be the center of the $\Theta$-stratum $
\ev_1 \co \cZ^+ \hookarr \cX$. Then for any point $x \in |\cZ|$ and any open substack $\cU \subset \cX$ containing $\sigma(x)$, there is another open substack with $\sigma(x) \in \cV \subset \cU$ such that $\cZ^+ \cap \cV$ is a   $\Theta$-stratum in $\cV$.
\end{lem}
\begin{proof} We only need to find a substack $\cV \subset \cX$ containing $\sigma(x)$ such that for any $f\co \Theta_k \to \cX$, where $k$ is a field, with $f \in \cZ^+$ and $f(1) \in \cV$, we have $f(0) \in \cV$ as well.
Let $\cU' = (\ev_1 \circ \sigma)^{-1}(\cU) \subset \cZ$, and let $\cZ' = \cZ \smallsetminus \cU'$ be its complement. Then the open substack
\[
\cV:=\cU \smallsetminus (\cU \cap \ev_1(\gr^{-1}(\cZ'))) \subset \cX
\]
satisfies the condition.
\end{proof}

\subsubsection{Reminder on the normal cone to a $\Theta$-stratum}

The main problem in finding a presentation of the form $[\Spec(A/I_+)/\bG_m] \subset [\Spec(A)/\bG_m]$ is that for an arbitrary morphism $[\Spec(A)/\bG_m] \to \cX$ the preimage of the $\Theta$-stratum need not be defined by the ideal generated by the elements of positive weight. To find presentations for which this happens, we need to recall that the weights of the $\bG_m$-action of the restriction of the conormal bundle of a $\Theta$-stratum to its center $\cZ$ are automatically positive. This property was already important in the work of Atiyah--Bott \cite{AtiyahBott} and it appears in the language of spectral stacks in \cite[\S 1.2]{hlinstability}. For completeness we provide a classical argument:

\begin{lem}\label{L:stratum_weights}
	In the setting of \Cref{T:Langton-Algorithm}, let $\sigma \co \cZ \to \cZ^+$ be the center of the $\Theta$-stratum $
\ev_1 \co \cZ^+ \hookarr \cX$ and $x \in \cZ(k)$ be a $k$-point. By abuse of notation we will also denote $\sigma(x) \in \cX(k)$ by $x$.
	\begin{enumerate}
		\item Let $T_{\cX,x}=\bigoplus_{n\in \bZ} T_{\cX,x,n}$ be the decomposition of the tangent space at $x$ into weight spaces with respect to the $\bG_m$-action induced from the canonical cocharacter $\lambda_x\co \bG_m \to \Aut_{\cX}(x)$. Then we have $T_{\cZ^+,x}=\bigoplus_{n\geq 0}T_{\cX,x,n}$.
		\item $\bG_m$ acts with non-negative weights on $\Lie(Aut_{\cX}(x))$.
	\end{enumerate}
\end{lem}

\begin{proof}	
	Let us first show that $\bigoplus_{n\geq 0}T_{\cX,x,n}\subseteq T_{\cZ^+,x}$.
	Let $t\in \cX(k[\epsilon]/\epsilon^2)$ be a tangent vector in $T_{\cX,x,n}$ for some $n\geq 0$, i.e., $t$ comes equipped with an isomorphism $t \mod \epsilon \cong x$. 
	
	This means that we have a commutative diagram 
	$$\xymatrix@C=5em{[\Spec(k[\epsilon]/\epsilon^2)/\bG_m]\ar[r]^-t& \cX,\\ [\Spec(k)/\bG_m] \ar[ur]_-{(x,\lambda_x)}\ar@{^(->}[u] & }$$ where $\bG_m$ acts on $\Spec(k[\epsilon]/\epsilon^2)$ via $(\lambda,\epsilon) \mapsto \lambda^n\epsilon$. In other words, we have a commutative diagram
	$$\xymatrix@C=5em{
		\bG_m \times \Spec(k[\epsilon]/\epsilon^2) \ar[r]^-{(\lambda,\epsilon \mapsto \lambda^n \epsilon)} \ar[dr] & \Spec(k[\epsilon]/\epsilon^2) \ar[d]^t \\
		& \cX.
	}$$
	If $n\geq 0$ then the horizontal map extends to {a} {$\bG_m$-equivariant map on} $\bA^1$, i.e., we get an extension 
	$$\xymatrix@C=5em{
		\bA^1 \times \Spec(k[\epsilon]/\epsilon^2) \ar[r]^-{(\lambda,\epsilon \mapsto \lambda^n \epsilon)} \ar[dr] & \Spec(k[\epsilon]/\epsilon^2) \ar[d]^t \\
		& \cX
	}$$
	and this defines an extension of $t$ to a $k[\epsilon]/\epsilon^2$-valued point of $\uMap_S(\Theta,\cX)$.
	
	Conversely, an extension of the constant map $[\bA^1/\bG_m] \to [\Spec(k)/\bG_m] \to \cX$ to $[\bA^1 \times \Spec( k[\epsilon]/(\epsilon^2)) /\bG_m] \to \cX$ automatically factors through the first infinitesimal neighborhood of $x\in \cX$. On a versal first order deformation this corresponds to a homomorphism of graded algebras $k[\epsilon_1,\dots,\epsilon_d]/(\epsilon_i)_{i=1,\dots d}^2 \to k[\lambda,\epsilon]/(\epsilon^2)$, where we can choose $\epsilon_i$ to be homogeneous for the $\bG_m$-action defined by $\lambda_x$.  
	This has to vanish on those tangent directions $\epsilon_i$ on which $\lambda_x$ acts with negative weights. This shows (1).
	
	Similarly for (2), when we regard $x$ as a $k$ point of $\cZ \hookrightarrow \cZ^+ \subset \uMap_S(\Theta,\cX)$, it corresponds to a map which factors as $\Theta_k \to B_k G_x\hookrightarrow \cX$, where we abbreviated $G_x=\Aut_{\cX}(x)$. We know that $\Aut_{{\cZ^+}}(x) \to \Aut_{\cX}(x)$ is an equivalence, so by the classification of $G_x$-bundles on $[\bA^1/\bG_m]$ (see \cite[Lem.~1.7]{heinloth-stability} or \cite[Prop.~A.0.1]{hlinstability}) 
	this implies that for the canonical cocharacter $\lambda_x \co \bG_m \to G_x$ we have $G_x=P(\lambda_x)$ as an algebraic group. In particular this means that $\bG_m$ acts with non-negative weights on the Lie algebra of $G_x=P(\lambda_x)$.
\end{proof}

\subsubsection{Reduction to the basic situation - Case of smooth stabilizers over a field}

\begin{lem} \label{L:stratum_chart_2} 
	Let $\cX$ be an algebraic stack  of finite type with affine diagonal  over an algebraically closed field $k$. 
	Let $\cZ^+\subset \cX$ be a $\Theta$-stratum with center $\sigma \co \cZ \to \cZ^+$, and let $x_0 \in \cZ(k)$ be a point such that $x:=\sigma(x_0)$ has a smooth automorphism group. Then there is a smooth representable morphism $p \co [\Spec(A)/\bG_m] \to \cX$ whose image contains $x$ and such that 
	\[
	 p^{-1}(\cZ^+) = [\Spec(A/I_+) / \bG_m] \hookrightarrow [\Spec(A)/ \bG_m].
	\]
\end{lem}
\begin{proof}
	The point $x_0$ has a canonical non-constant homomorphism $\bG_m \to \Aut_{\cZ}(x_0)$, which induces a canonical homomorphism $\lambda: \bG_m \to G_{x} := \Aut_{\cX}(x)$. We may replace $\bG_m$ with its image in $G_{x}$ and thus assume that $\lambda$ is injective. As we assumed that $G_x$ is smooth the quotient $G_{x} / \lambda(\bG_m)$ is smooth, so we may apply \cite[Thm.~1.2]{ahr}  to obtain a smooth representable morphism 
	$$ p \co [\Spec(A)/\bG_m] \to \cX$$
	together with a point $w\in \Spec(A)(k)$ in $p^{-1}(x)$ which is fixed by $\bG_m$ and such that $p^{-1}(B_k G_{x})\cong B_k \bG_m$. The isomorphism $p^{-1}(B_k G_{x})\cong B_k \bG_m$ implies that the relative tangent space to $\tilde{p} \co \Spec(A) \to \cX$ at $w$ is naturally identified with $\Lie(G_x) / \Lie(\bG_m)$ on which $\bG_m$ acts with non-negative weights by part (2) of \Cref{L:stratum_weights}.
	
	Note that connected components of $\Spec(A)^{\bG_m}$ can be separated by invariant functions, so we may replace $\Spec(A)$ with a $\bG_m$-equivariant affine open neighborhood of $w$ so that $\Spec(A)^{\bG_m}$ is connected. It follows that $\Spec(A/I_+)$ is connected as well.

	This implies that $\cZ^+_A:=[\Spec(A/I_+)/\bG_m]\subseteq [\Spec(A)/\bG_m]$ is isomorphic to a connected component of $\uMap(\Theta,[\Spec(A)/\bG_m])$ and $\cZ_A:=[\Spec(A^{\bG_m})/\bG_m]\subset \cZ^+_A$ is the center of $\cZ^+_A$. As $p(x)\in \cZ$ connectedness now implies that $p(\cZ_A)\subset \cZ_0$ and therefore we also have  $[\Spec(A/I_+)/\bG_m] \subset p^{-1}(\cZ^+)$.
	
	To conclude that $\cZ^+_A \cong p^{-1}(\cZ^+)$ after possibly shrinking $A$, it suffices to check that the inclusion $[\Spec(A / I_+)/\bG_m] \subseteq p^{-1}(\ev_1(\cZ^+))$ of closed substacks of $[\Spec(A)/\bG_m]$ is an isomorphism locally at $w$. Consider the pull-back:
	
	$$\xymatrix{
		p^{-1}(\ev_1(\cZ^+))=\Spec(B) \ar[d] \ar@{^(->}[r]& \Spec(A)\ar[d]^{p}\\
		\cZ^+\ar@{^(->}[r] & \cX
	}$$
	Then $B$ is a graded ring and we still have an exact sequence $$T_{p,w} \to T_{\Spec(B),w} \to \cT_{\cZ^+,x}.$$ As $\bG_m$ acts with non-negative weight on the relative tangent bundle at $w$ and also on $\cT_{\cZ^+,x}$ by  \Cref{L:stratum_weights}, this shows that $\bG_m$ acts with non-negative weights on $T_{\Spec(B),w}$. In particular the maximal ideal $\mathfrak{m}_{w} \subset B$ of $w$ is generated by elements of non-positive weight locally at $w$. 
	
	Therefore, after possibly shrinking $A$ we may assume that $B =\bigoplus_{n \leq 0} B_n$ is non-positively graded.  
	As $\Spec(A/I_+) \subset \Spec(A)$ was the contracting subscheme for $\bG_m$ we find that locally around $w$ we thus have  $p^{-1}(\cZ^+) \subset \Spec(A / I_+)$ locally around $w$. This proves our claim.
\end{proof}
\subsubsection{Reduction to the basic situation - general case}
\begin{lem} \label{L:stratum_chart}
In the setting of \Cref{T:Langton-Algorithm} where $\cZ^+ \hookarr \cX$ is a $\Theta$-stratum and $\cX$ is quasi-compact, there is a smooth representable morphism $p : [\Spec(A)/\bG_m] \to \cX$ such that $p^{-1}(\cZ^+)$ is the $\Theta$-stratum
\[
 p^{-1}(\cZ^+) = [\Spec(A/I_+) / \bG_m] \hookrightarrow [\Spec(A)/ \bG_m],
\]
and $\cZ^+$ is contained in the image of $p$.
\end{lem}
\begin{proof}
Because $\cX$ is finite type over the base space $S$, we may apply \Cref{L:filtration_atlas} to obtain a smooth, surjective and representable map $p \co [\Spec(A)/\bG_m^n] \to \cX$ such that $\uMap(\Theta,[\Spec(A)/\bG_m^n]) \to \uMap(\Theta,\cX)$ is also smooth, surjective and representable. From \Cref{P:theta-description}, we know that $\uMap(\Theta,[\Spec(A)/\bG_m^n])$ is the disjoint union indexed by cocharacters $\bG_m \to \bG_m^n$ of stacks of the form $[\Spec(A/I_+) / \bG_m^n]$, where $I_+$ is the ideal generated by positive weight elements with respect to a given cocharacter. Choosing different connected components if necessary and forgetting all but the relevant cocharacter in each component, we can construct a non-positively graded algebra $C = \bigoplus_{n \leq 0} C_n$ along with a smooth surjective representable map $[\Spec(C)/\bG_m] \to \cZ^+$.

We now discard the previously constructed $\Spec(A)$ and apply the relative slice theorem (\Cref{T:relativeslice}) to the smooth surjective map $[\Spec(C)/\bG_m] \to \cZ^+$, where we regard $\cZ^+$ as a closed substack of $\cX$. This provides a map $p : [\Spec(A')/\bG_m] \to \cX$ along with an isomorphism $C \cong A' / I_{p^{-1}(\cZ^+)}$, where $I_{p^{-1}(\cZ^+)}$ is the ideal corresponding to $p^{-1}(\cZ^+)$. By construction $C$ has no positive weight elements, so the ideal $I_+$ generated by positive weight elements of $A'$ is contained in $I_{p^{-1}(\cZ^+)}$.

Because $p$ is smooth, the relative cotangent complex of $\Spec(C) \hookrightarrow \Spec(A')$ is $p^\ast(\bL_{\cZ^+/\cX})$. In particular, the fiber of the conormal bundle of $\Spec(C) \hookrightarrow \Spec(A')$ has positive weights at every point of $\Spec(C)^{\bG_m}$ by \Cref{L:stratum_weights}. One may therefore find a collection of positive weight elements of $I_{p^{-1}(\cZ^+)}$ that generate the fiber of $I_{p^{-1}(\cZ^+)}$ at every closed point of $\Spec(C)^{\bG_m}$.

Moreover, as $C$ is non-positively graded, the orbit closure of every point in $\Spec(C)$ meets the fixed locus $\Spec(C)^{\bG_m}$. So by Nakayama's lemma we can actually find a collection of homogeneous elements of $I_+$ that generate the fiber of $I_{p^{-1}(\cZ^+)}$ at every point of $\Spec(C)$ and hence in a $\bG_m$-equivariant open neighborhood of $\Spec(C) \hookrightarrow \Spec(A')$. We may thus invert a weight $0$ element $a \in A'$ so that these elements of $I_+$ generate $(I_{p^{-1}(\cZ^+)})_a \subset A'_a$ and $C = A' /I_{p^{-1}(\cZ^+)} = A'_a / (I_{p^{-1}(\cZ^+)})_a$ is unaffected.

In particular we have shown that after inverting a weight $0$ element of $A'$, we have a smooth map $p \co [\Spec(A')/\bG_m] \to \cX$ such that $[\Spec(A'/I_+)/\bG_m] = p^{-1}(\cZ^+)$ and the map $\Spec(A'/I_+) \to \cZ^+$ is surjective.
\end{proof}

We can now prove the semistable reduction theorem:

\begin{proof}[Proof of \Cref{T:Langton-Algorithm}]
Consider a map $\xi \co \Spec(R) \to \cX$ as in the statement of the theorem. Observe that for any smooth map $p \co \cY \to \cX$ such that {$\cZ^+$} induces a $\Theta$-stratum $p^{-1}(\cZ^+)$ in $\cY$ and the image of $p$ contains the image of $\xi$, if we know the conclusion of the theorem holds for $\cY$ then the conclusion holds for $\cX$ as well: indeed after an extension of $R$ we may lift $\xi$ to a map $\xi' \co \Spec(R') \to \cY$, construct an elementary modification in $\cY$ such that the new map $\xi'' \co \Spec(R') \to \cY$ lies in $\cY \smallsetminus p^{-1}(\cZ^+)$, and observe that the composition of this elementary modification with $p$ gives an elementary modification of $\xi$ such that the new map $p \circ \xi'' \co \Spec(R') \to \cX$ lies in $\cX \smallsetminus \cZ^+$.

Using this observation and the fact that $\xi_k$ lies in $\cZ^+$, we may use \Cref{L:quasi-compact_stratum} to replace $\cX$ with a quasi-compact open substack, then use \Cref{L:stratum_chart} to construct a smooth map $p \co [\Spec(A)/\bG_m] \to \cX$ whose image contains the image of $\xi$ and for which $\cZ^+$ induces a $\Theta$-stratum. Then we are finished by \Cref{L:graded_affine_reduction}.
\end{proof}

\subsection{Comparison between a stack and its semistable locus}

As an immediate consequence of the semistable reduction theorem, we have the following:
\begin{cor} \label{C:valuative_criterion_semistable}
Let $\cX$ be a quasi-separated algebraic stack with affine {stabilizers} that is locally finite type over a noetherian algebraic space $S$. Let $\cX = \bigcup_{c \in \Gamma} \cX_{\leq c}$ be a well-ordered $\Theta$-stratification of $\cX$. If $\cX \to S$ satisfies the existence part of the valuative criterion for properness with respect to DVRs, then so does $\cX_{\leq c} \to S$ for every $c \in \Gamma$.  In particular, if the semistable locus $\cX^{\ss}:=\cX_{\le 0}$ is quasi-compact, then $\cX^{\ss} \to S$ is universally closed.
\end{cor}
\begin{proof} 
Consider a DVR $R$ and a map $\Spec(R) \to S$ along with a lift $\Spec(K) \to \cX_{\le c}$. If $\cX \to S$ satisfies the existence part of the valuative criterion, then after an extension of $R$ one can extend this lift to a lift $\Spec(R') \to \cX$ of $\Spec(R) \to S$. By hypothesis the generic point lies in $\cX_{\le c}$, so by \Cref{T:semistable_reduction} after passing to a further extension of $R$ there is a sequence of elementary modifications resulting in a modification $\Spec(R') \to \cX_{\le c}$. Note that because $\Spec(R)$ is the good moduli space of $\oST_R$, and good moduli spaces are universal for maps to an algebraic space \cite[Thm.~6.6]{alper-good}, any elementary modification of a map $\Spec(R) \to S$ is trivial. It follows that our modified map $\Spec(R') \to \cX_{\le c}$ is a lift of the original map $\Spec(R) \to S$. By \cite[Lem.~2.4.6]{hlp}, one can check that a finitely presented morphism of noetherian algebraic stacks $\cX \to \cY$ is universally closed by checking that $\cX \times \bA^n \to \cY \times \bA^n$ is closed for all $n\geq 0$. From this it follows that the existence part of the valuative criteria with respect to DVRs implies universal closedness.
\end{proof}

Next let us briefly recall the notion of $\Theta$-stability from \cite[Def.~4.1.1 \& 4.1.3]{hlinstability} and \cite[Def.~1.2]{heinloth-stability}.

\begin{defn} \label{D:semistable-locus}
Given a cohomology class $\ell \in H^2(\cX;\bR)$, we say that a point $p \in |\cX|$ is \emph{unstable} with respect to $\ell$ if there is a filtration $f \co \Theta_k \to \cX$ with $f(1) = p \in |\cX|$ and such that $f^\ast(\ell) \in H^2(\Theta_k;\bR) \simeq \bR$ is positive. The \emph{$\Theta$-semistable locus} $\cX^{\ss}$ is the set of points which are not unstable. 
\end{defn}
The above definition is simply an intrinsic formulation of the Hilbert--Mumford criterion for semistability in geometric invariant theory.
We are somewhat flexible with what type of cohomology theory we use: if $\cX$ is locally finite type over $\bC$ one may use the Betti cohomology of the analytification of $\cX$; if $\cX$ is locally finite type over another field $k$, one can use Chow cohomology; and in general one may use the Neron--Severi group $NS(\cX)_\bR$ for $H^2(\cX;\bR)$. In \cite[\S 4.1.2]{hlinstability}, the properties of the cohomology theory needed for the theory of $\Theta$-stability are axiomatized.

\begin{prop}   \label{P:semistable_locus}
Let $\cX$ be an algebraic stack locally of finite type with affine diagonal over a noetherian algebraic space $S$, and let $\cX^{\ss}$ be the $\Theta$-semistable points with respect to a class $\ell \in H^2(\cX;\bR)$.  Suppose that either   
\begin{enumerate}[label=(\alph*)]
\item $\cX^{\ss}$ is the open part of a $\Theta$-stratification of $\cX$, i.e., $\cX^{\ss} = \cX_{\leq 0}$, such that for each HN filtration $g \co \Theta_k \to \cX$ of an unstable point one has $g^\ast(\ell) >0$ in $H^2(\Theta_k;\bR)$, or
\item $\cX^{\ss} \subset \cX$ is open and $\cX \to S$ is $\Theta$-reductive.
\end{enumerate}
Then
\begin{enumerate}
\item if $\cX \to S$ is  {\textsf{S}}-complete, then so is $\cX^{\ss} \to S$, and
\item if $\cX \to S$ is $\Theta$-reductive, then so is $\cX^{\ss} \to S$.
\end{enumerate}
\end{prop}

In the proof, we will need the following:

\begin{lem}\label{L:semistable_filtrations}
Under the hypotheses of \Cref{P:semistable_locus}, given a filtration $f \co \Theta_k \to \cX$ such that $f(1)$ is semistable with respect to $\ell$, then $f^\ast(\ell) = 0$ if and only if $f(0)$ is semistable as well.
\end{lem}
\begin{proof}
The proof is a geometric reformulation of the corresponding argument for semistability for vector bundles. 
For any semistable point $x \in \cX(k)$ and any cocharacter $\lambda \co \bG_m \to G_x$, the restriction of $\ell$ to $H^2([\Spec(k)/\bG_m];\bR)\simeq \bR$ along the resulting map $f_\lambda \co \Theta_k \to [\Spec(k)/\bG_m] \to \cX$ must vanish, because the invariants for $\lambda$ and $\lambda^{-1}$ differ by sign and are both non-positive. 

For the converse suppose that $f^\ast(\ell)=0$ and $f(0) \notin \cX^{\ss}$.  
We claim that there is a filtration $g \co \Theta_k \to \cX$ of $f(0)$ with $g^*(l) > 0$ which is invariant under the action of $\bG_m$ on $f(0)$ induced from the filtration $f \co \Theta_k \to \cX$.  This is automatic in case (a) as HN filtrations are canonical.  For case (b), since $\cX \to S$ is $\Theta$-reductive,  the representable map $\uMap(\Theta,\cX)\to \cX$ satisfies the valuative criterion for properness, so the fiber of this map over $f(0) \in \cX(k)$, which is denoted $\Flag(f(0))$ is an algebraic space {locally} of finite type over $k$ which satisfies the valuative criterion for properness. The action of $\bG_m$ by automorphisms of $f(0)$ gives a $\bG_m$-action on $\Flag(f(0))$. Given some point $g_1 \in \Flag(f(0))(k)$ for which $g_1^\ast(\ell)>0$, we can consider the orbit $\bG_m \to \Flag(f(0))$ of $g_1$. Because $\Flag(f(0))$ satisfies the valuative criterion for properness, this map extends to an equivariant map $\bA^1_k \to \Flag(f(0))$. This map sends $0 \in \bA^1_k$ to a fixed point for the action of $\bG_m$ on $\Flag(f(0))$, which corresponds to a $\bG_m$-invariant filtration $g$ of $f(0)$, and $g$ is on the same connected component of $\Flag(f(0))$ as $g_1$, so $g^\ast(\ell) = g_1^\ast(\ell) > 0$.

Denote by $R=k[\![\pi]\!]$ the completion of the local ring of the affine line with coordinate $\pi$ at $0$. Then the map $f_R \co \Spec(R) \to [\Spec(k[\pi])/\bG_m] = \Theta_k \map{f} \cX$ and $g \co \Theta_k \to \cX$ define the datum needed to apply the gluing lemma \Cref{C:strange_gluing}, which says that after restricting $f_R$ to $R' = R[\pi^{1/n}]$ for $n\gg 0$ there is a unique extension $ F_{R'} \co \overline{ST}_{R'} \to \cX$ such that $F|_{t\neq 0} \cong f_{R'}$ and $F|_{s=0}\cong g$. Let $\pi^\prime = \pi^{1/n}$ denote the uniformizer in $R'$.

As $f_{R'}$ was the restriction of a map $f_{\bA^1} \co \bA_{k}^1 \to \Theta_{k} \to \cX$ we find that this morphism extends canonically to
$$ F \co [\Spec(k[\pi^\prime,s,t]/(st-\pi^\prime))/\bG_m] = [\Spec(k[s,t]) / \bG_m] \to \cX.$$
By uniqueness of the extension $F_{R'}$ and the fact that $g$ is fixed by the $\bG_m$-action on $f(0)$ induced by $f$, this morphism comes equipped with a descent datum for the standard $\bG_m$-action on $\bA^1=\Spec(k[\pi^\prime])$. We therefore obtain
$$ \overline{F} \co [\Spec(k[\pi^\prime,s,t]/(st-\pi^\prime))/\bG_m^2] = [\bA^2_k / \bG_m^2] \to \cX.$$
where the action of the second copy of $\bG_m$ is with weight $-1$ on $s$ and trivial on $t$. Choosing a different basis for the cocharacter lattice of $\bG_m^2$, we see that this is equivalent to the usual action of $\bG_m^2$ on $\bA_k^2$. In particular, every cocharacter $\lambda\co \bG_m \to \bG_m^2$ has the form $\lambda(t) = (t^a,t^b)$ for some pair $\langle a,b\rangle$.

If $a,b \geq 0$, then the point $(1,1) \in \bA^2_k$ has a limit under $\lambda(t)$ as $t\to 0$, and restricting $\overline{F}$ to the corresponding line we obtain a filtration $f_{\langle a,b\rangle}\co \Theta_k \to \cX$, with $f(1) \simeq f_{\langle a,b\rangle}(1)$. By construction the original filtration $f$ corresponds to $f_{\langle1,0\rangle}$.

Note that $\overline{F}^\ast(\ell)|_{[0/\bG_m^2]}$ defines a character $\chi$ of $\bG_m^2$ such that for $a,b>0$ the weight $f^*_{\langle a,b\rangle}(\ell)$ is given by the canonical pairing $(\langle a,b\rangle, \chi)$ between characters and cocharacters. As $F|_{\{s=0\}}=g$ was assumed to be destabilizing and $F$ was defined by the subgroup $\langle -1,1\rangle$ we have $(\langle -1,1 \rangle, \chi) = g^*(l) >0$. We also have $(\langle1,0\rangle,\chi) = f_{\langle 1,0\rangle}^\ast(\ell) = 0$ by hypothesis. Thus for $a\gg b > 0$ we have $f^*_{\langle a,b\rangle}(\ell) >0$, contradicting the assumption that $f(1)$ was semistable.
\end{proof}

\begin{rem} \label{R:more_general_stability}
In \Cref{T:abelian_semistable_GMS} below we will use a slightly more general notion of stability: we replace the weight of $f^\ast(\ell)$ with any function $\ell$ from the set of filtrations in $\cX$ to a totally ordered real vector space $V$, and define $x \in \cX$ to be semistable if $\ell(f) \leq 0$ for any filtration $f$ with $f(1) = x$. Then the proof above applies verbatim, provided that 1) $\ell(f)$ is locally constant in algebraic families of filtrations, and 2) for any map $F \co [\bA^2_k / \bG_m^2] \to \cX$ the function on the cocharacter lattice $\lambda \mapsto \ell(f_\lambda) \in V$ is linear, where $f_\lambda$ denotes the filtration associated to the cocharacter $\lambda$ as in the proof of \Cref{L:semistable_filtrations}. The second condition is equivalent to requiring the function $\lambda \mapsto \ell(g_\lambda)$ is linear, where $g_\lambda$ denotes the constant filtration of $F(0,0)$ induced from the $\bZ$-grading of $F(0,0)$ associated to the cocharacter $\lambda$.
\end{rem}

\begin{rem}
The proof of \Cref{L:semistable_filtrations} is a special case of the technique used to prove the perturbation theorem on filtrations \cite[Thm.~3.5.3]{hlinstability}. This theorem constructs a bijection between filtrations of a point $x \in \cX$ that are ``close'' to a given filtration $f$ and filtrations of the associated graded object $f(0) \in \Map_k(B \bG_m,\cX)$ that are ``close'' to the canonical filtration defined by the action of $\bG_m$ on $f(0)$. In this language, the proof that $f(0)$ is semistable if $f(1)$ is semistable and $f^\ast(\ell)=0$ amounts to the observation that if $f(0)$ had a destabilizing filtration as a graded object, then because the canonical filtration of $f(0)$ has weight $0$, one can find destabilizing filtrations of $f(0)$ which are arbitrarily close to the canonical filtration, and then one can identify these with destabilizing filtrations of $f(1)$ using the perturbation theorem.
\end{rem}

\begin{proof}[Proof of \Cref{P:semistable_locus}]
Consider a DVR $R$ and a diagram
\[
\xymatrix{ \Spec(R) \cup_{\Spec(K)} \Spec(R) \ar[r] \ar[d] & \cX^{\ss} \ar[d] \\ \oST_R \ar@{-->}[ur] \ar[r] & S }.
\]
By hypothesis we can fill the dotted arrow uniquely to a map $\oST_R \to \cX$. We claim that in fact the map $\oST_R \to \cX$ factors through $\cX^{\ss}$. Because $\cX^{\ss}$ is open, it suffices to check that the unique closed point maps to $\cX^{\ss}$. By hypothesis the point $(\pi,s,t) = (0,1,0)$ and the point $(\pi,s,t)=(0,0,1)$ map to $\cX^{\ss}$. Restricting the map $\oST_R \to \cX$ to the locus $\Theta_k \simeq \{s=0\}$ and $\Theta_k \simeq \{t=0\}$ give filtrations $f_1$ and $f_2$ in $\cX$ of points in $\cX^{\ss}$, and if one has $f_{1}^\ast(\ell)<0$ then the other has $f_{2}^\ast(\ell)>0$, which would contradict the fact that $f(1) \in \cX^{\ss}$. Therefore $f^\ast(\ell) = 0$ for both filtrations, and it follows from \Cref{L:semistable_filtrations} that $f(0) \in \cX^{\ss}$ as well.

For the corresponding claim for $\Theta$-reductivity is proved similarly. For the analogous filling diagram, we start with a map $f \co \Theta_R \smallsetminus \{(0,0)\} \to \cX^{\ss}$ and fill it to a map $\tilde{f} \co \Theta_R \to \cX$. We claim that $(0,0)$ maps to $\cX^{\ss}$ as well, and hence because $\cX^{\ss} \subset \cX$ is open it follows that $\tilde{f}$ lands in $\cX^{\ss}$. Because the restriction $f_K$ of $f$ to $\Theta_K \subset \Theta_R \smallsetminus \{(0,0)\}$ maps to $\cX^{\ss}$, we know from \Cref{L:semistable_filtrations} that $f_K^\ast(\ell) = 0$. The function $f \mapsto f^\ast(\ell) \in \bR$, regarded as a function on $\Map(\Theta,\cX)$, is locally constant. It therefore follows that the restriction $\tilde{f}_k \co \Theta_k \to \cX$ of $\tilde{f}$ also has $\tilde{f}_k^\ast(\ell) = 0$. It follows that $\tilde{f}_k(0) \in \cX^{\ss}$.
\end{proof}

\begin{cor} \label{P:moduli_space_semistable_locus}
Let $\cX$ be an algebraic stack locally of finite type with affine diagonal over a noetherian algebraic space $S$ defined over $\bQ$. Assume that $\cX \to S$ is  {\textsf{S}}-complete and $\Theta$-reductive. Let $\cX^{\ss} \subset \cX$ be the $\Theta$-semistable locus with respect to some class $\ell \in H^2(\cX;\bR)$. If $\cX^{\ss} \subset \cX$ is a quasi-compact open substack, then $\cX^{\ss}$ admits a good moduli space which is separated over $S$. Furthermore if in addition $\cX \to S$ satisfies the existence part of the valuative criterion for properness and $\cX^{\ss}$ is the open part of a well-ordered $\Theta$-stratification of $\cX$, then the good moduli space for $\cX$ is proper over $S$.
\end{cor}
\begin{proof}
The map $\cX^{\ss} \to S$ is  {\textsf{S}}-complete and $\Theta$-reductive by \Cref{P:semistable_locus}.  
By \Cref{T:A}, $\cX^{\ss}$ admits a separated good moduli space $\cX \to X$. Under the additional hypotheses, $\cX^{\ss} \to S$ satisfies the existence part of the valuative criterion for properness by \Cref{C:valuative_criterion_semistable}, and hence $X \to S$ is proper by \Cref{P:gms-S-complete}.
\end{proof}

\subsection{Application: Properness of the Hitchin fibration}

Let us illustrate how the semistable reduction theorem (\Cref{T:semistable_reduction}) can be used to simplify and extend classical semistable reduction theorems for principal bundles and Higgs bundles on curves. 

The setup for these results is the following (see e.g., \cite[\S 2]{NgoHitchin}). Let $C$ be a smooth projective, geometrically connected curve over a field $k$ and $G$ a reductive algebraic group. As the notions are slightly easier to formulate over algebraically closed fields and the valuative criteria allow for extensions of the ground field, we will assume that $k$ is algebraically closed in this section. 

We denote by $\Bun_G$ the stack of principal $G$-bundles on $C$, i.e., for a $k$-scheme $S$ we have that $\Bun_G(S)$ is the groupoid of principal $G$-bundles on $C\times S$. Fix a line bundle $\cL$ on $C$. A $G$-Higgs bundle with coefficients in $\cL$ on $C$ is a pair $(\cP,\phi)$ where $\cP$ is a $G$-bundle on $C$ and $\phi\in H^0(C,(\cP\times^G \Lie(G))\tensor \cL)$. We denote by $\Higgs_G$ the stack of $G$-Higgs bundles with coefficients in $\cL$.

The stack $\Higgs_G$ comes equipped with the forgetful morphism $\Higgs_G \to \Bun_G$ and the Hitchin morphism $h \co \Higgs_G \to \cA_G$. Here $\cA_G\cong \bigoplus_{i=1}^r H^0(C,\cL^{d_i})$ {is the Hitchin base, i.e., the affine space over $k$ defined by the vector space of global sections of the line bundles $\cL^{d_i}$}, where $d_1,\dots,d_r$ are the degrees of homogeneous generators of $k[\Lie(G)^\ast]^G$ and $h$ is defined by mapping $(\cP,\phi)$ to the characteristic polynomial of $\phi$.  

On both $\Bun_G$ and $\Higgs_G$ there is a classical notion of stability, which is defined in terms of reductions to parabolic subgroups. 

Let us recall how this notion is related to $\Theta$-stability. For vector bundles there is an equivalence (\Cref{P:theta-description}, \cite[Lem.~1.10]{heinloth-stability})
$$ \uMap(\Theta,\Bun_{\GL_n}) \cong  \left\langle (\cE,\cE^i)_{i\in \bZ} \left| \begin{array}{l} \cE \in \Bun_{\GL_n}, \cE^i\subseteq \cE^{i+1} \subseteq \cE \text{ subbundles} \\
\cE^i = \cE \text{ for } i\gg 0, \cE^i = 0 \text{ for } i\ll 0
\end{array}  \right.\right\rangle$$
which is given by assigning to a weighted filtration of a vector bundle $\cE$ the canonical $\bG_m$-equivariant degeneration of $\cE$ to the associated graded bundle. 

This construction has an analog for principal bundles. To state this we fix (as in \Cref{P:theta-description}) a complete set of conjugacy classes of cocharacters $\Lambda \subset \Hom(\bG_m,G)$. As in \S \ref{S:drinfeld} we denote by $P_\lambda^+\subseteq G$ the parabolic subgroup defined by $\lambda$ and by $L_\lambda \subset P_\lambda^+$ the Levi subgroup defined by $\lambda$ which is isomorphic to the quotient of $P_\lambda^+$ by its unipotent radical $U_\lambda^+ \subset P_\lambda^+$. Then there is an equivalence (see e.g., \cite[Lem.~1.13]{heinloth-stability}) 
$$ \uMap(\Theta, \Bun_G)  \cong  \coprod_{\lambda \in \Lambda} \Bun_{P_\lambda^+}.$$

For Higgs bundles note that the forgetful map $\Higgs_G\to \Bun_G$ is representable and therefore $\uMap(\Theta,\Higgs_G) \subset \uMap(\Theta,\Bun_G)\times_{\Bun_G}\Higgs_G$, i.e., a filtration of a Higgs bundle is the same as a filtration of the underlying principal bundle that preserves the Higgs field $\phi$. 

Recall that a $G$-bundle $\cE$ is called semistable, if for all $\lambda$ and all $\cE_\lambda \in  \Bun_{P_\lambda^+}$ with $\cE_\lambda\times^{P_\lambda^+}G \cong \cE$ we have $\deg(\cP_\lambda\times^{P_\lambda^+} \Lie(P_\lambda^+)) \leq 0$.  Similarly a Higgs bundle is called semistable if the same condition holds for all reductions that respect the Higgs field $\phi$. This stability notion can be viewed as $\Theta$-stability induced from the so called determinant line bundle $\cL_{\det}$ on the stack $\Bun_G$ whose fiber at a point $\cP$ (resp. a point $(\cP,\phi)$)  is given by the one dimensional vector space $\det(H^1(C,\cP\times^G \Lie(G))) \tensor \det(H^0(C,\cP\times^G \Lie(G)))^{-1}$ (see \cite[\S 1.F]{heinloth-stability}, \cite{hlinstability}).

As usual we denote by $\Bun_G^{\ss} \subseteq \Bun_G$ the open substack of semistable bundles and for by $\Bun_{P_\lambda^+}^{\ss} \subseteq \Bun_{P_\lambda^+}$ the open substack of bundles such that the associated $L_\lambda$-bundle is semistable.

Finally let us recall how the notion of Harder--Narasimhan reduction can be used to equip the stacks $\Bun_G$ and $\Higgs_G$ with a (well-ordered) $\Theta$-stratification if the characteristic of $k$ is not too small, i.e., such that Behrend's conjecture holds for $G$ (see \cite[Thm.~1]{heinloth-bounds} for explicit bounds depending on $G$; note that char $2$ has to be excluded for groups of type $B_n,D_n$ as well).

For any unstable $G$-bundle $\cP$ there exists a canonical Harder--Narasimhan reduction $\cP^{HN}$ to a parabolic subgroup $P_\lambda^+$, where $\lambda$ is uniquely determined up to a positive integral multiple. We denote by 
\begin{align*}
\underline{d}:=\underline{\deg}(\cP^{HN})\co \Hom(P_\lambda,\bG_m) &\to \bZ\\
\chi \mapsto \deg(\cP_\lambda^{HN} \times^\chi \bG_m)
\end{align*}
the degree of $\cP^{HN}$ and by $\Bun_{P_\lambda^+}^{\underline{d},\ss}\subset \Bun_{P_\lambda^+}^{\ss}$ the connected component defined by $\underline{d}$. The instability degree of $\cP_\lambda$ is defined as 
$$ \ideg(\cP) := \deg(\cP^{HN}\times^{P_\lambda^+}\Lie(P)).$$

Behrend showed that the morphism $\Bun_{P_\lambda^+}^{\underline{d},ss} \to \Bun_G$ defined by the inclusion $P_\lambda^+\subset G$ is radicial if the degree $\underline{d}$ is the degree of a canonical reduction \cite{behrend-thesis} and the map is an embedding if Behrend's conjecture holds for $G$ \cite[Lem.~2.3]{heinloth-semistable} this condition is satisfied if the characteristic of $k$ is not too small with respect to $G$ (e.g., $>31$).

Moreover the instability degree $\ideg$ is upper semicontinuous in families and if this invariant is constant on a family, then the family admits a global Harder--Narasimhan reduction \cite[Prop.~7.1.3]{behrend-thesis}\cite[Prop.~2.2]{heinloth-semistable}. Thus the Harder--Narasimhan reduction of bundles defines a $\Theta$-stratification on $\Bun_G$ if the characteristic of $k$ is not too small. The same arguments apply for Higgs bundles and this shows the following lemma.

\begin{lem}
	If the characteristic of $k$ is large enough so that Behrend's conjecture holds for $G$, then the Harder--Narasimhan stratifications of $\Bun_G$ and $\Higgs_G$ form a well-ordered $\Theta$-stratification.
\end{lem}	

To apply the semistable reduction theorem to $\Higgs_G \to \cA_G$ we need to show that this morphism satisfies the existence part of the valuative criterion for properness. The existence result is probably well known (see e.g. \cite[\S 8.4]{ChaudouardLaumonLFP1} for an argument over the regular locus) but we could not find a general reference. 
\begin{lem}\label{L:valuativeHiggs}
	Suppose that the characteristic of $k$ is not a torsion prime for  $G$ and very good for $G$.
	Let $R$ be a DVR with fraction field $K$, and let $(\cE_K,\phi_K)\in \Higgs_G(K)$ be a Higgs bundle such that $h(\cE_K,\phi_K)\in  \cA_G(K)$ extends to $\cA_G(R)$. Then there exists an extension $R \to R'$ of DVRs with $K \to K'=\Frac(R')$ finite	
	and a point $(\cE_R^\prime,\phi_R^\prime)\in \Higgs_G(R^\prime)$ extending $(\cE_K,\phi_K)$. 
\end{lem}
\begin{proof}
	First let us assume that the derived group of $G$ is simply connected. The generic point of $C$ will be denoted by $\eta$, $\cg=\Lie(G)$ and $\car:=\cg/\!/G$ is the space of characteristic polynomials of elements of $\cg$.
	
	Let $(\cE_K,\phi_K)\in \Higgs_G(K)$ be a Higgs bundle such that $h(\cE_K,\phi_K)\in \cA_G(R)\subset \cA_G(K)$. 
	
	We argue as in \cite[\S 8.4]{ChaudouardLaumonLFP1}. After a finite extension of $K$ we may assume that $\cE_K$ is trivial at the generic point $K(\eta)$ of $C_K$. Choosing trivializations of $\cE|_{K(\eta)}$ and $\cL|_\eta$ identifies $\phi_K$ with an element in $X_K\in \cg(K(\eta))$. To conclude the argument as in loc.cit., it is sufficient to show that after passing to a finite extension of $K$ we can conjugate $X_K$ to an element of $\cg(R(\eta))$,  because this allows one to extend $\phi_K$ to the trivial bundle over $R(\eta)$ and as sections of affine bundles extend canonically in codimension 2, the Higgs field $\phi_K$ will then  define a Higgs field for any extension $\cE_R$ of $\cE_K$ that is trivial over $R(\eta)$.
	
	We denote by $X_K=X_K^s + X_K^n$ the Jordan decomposition of $X_K$ into the semisimple and nilpotent part. 
	
	As $h(\cE,\phi)$ extends to $R$ we know that the image of $X_K$ in $\car=\cg/\!/G$ defines an $R(\eta)$-valued point. We can use the Kostant section $\car \to \cg$ to obtain $Y_R\in \cg(R(\eta))$ with $h(Y_R)=h(X_K)$. 
	
	We claim that we can modify $Y_R$ such that its generic fiber $Y_K$ is semisimple. To see this let us consider the Jordan decomposition $Y_K=Y_K^s + Y_K^n$. By our assumptions on the characteristic of $K$ the main result of \cite{McNinchSL2} shows that there exists a parabolic subgroup $P_{\lambda}^+\subset G$ defined by a cocharacter $\lambda \co \bG_m \to G_K$ such that $Y_K^n$ is contained in the Lie algebra of the unipotent radical of $P^+_{\lambda}$ and as $Y_K^s$ is in the centralizer of $Y_K^n$  the element $Y_K$ also lies in $\Lie(P_{\lambda}^+)$. As parabolic subgroups extend over valuation rings, we find that $Y_R$ is contained in a parabolic subgroup $P_R\subset G_R$ and we can choose $\lambda$ to be a cocharacter defined over $R$ as well. 
	
	As $P(\lambda)$ is defined to be the set of points {$p \in G$} such that $\lim_{t \to 0} \lambda(t)p\lambda(t)^{-1}$ exists, the limit $\lim_{t \to 0}\lambda(t) \cdot Y_R$ will be an $R$-valued point $Y^\prime_R$ such that $Y^\prime_K=Y_K^s$ is semisimple.
	
	As the semi-simple part of $X_K$ is the unique closed orbit in the conjugacy class of $X_K$ we know that $X_K^s$ and $Y_K^s$ lie in the same closed orbit. As we assumed that the derived group of $G$ is simply connected and that $p$ is not a torsion prime for $G$ the centralizer $Z_G(X_s)$ is a connected reductive group \cite[Thm.~0.1]{SteinbergTorsion}. By Steinberg's theorem, any $Z_G(X_s)$ torsor over $K(\eta)$ splits after a finite extension of $K$, so after possibly extending $K$ the elements $X_K^s$ and $Y_K^s$ are conjugate. Thus after conjugating $X_K$ we may assume that $X_K^s=Y_K^s$, i.e., we may assume that the semisimple part of $X_K$ extends to $R$.
	
	Now we can apply the previous argument to $X_K$, namely the element $X_K$ is contained in a parabolic subalgebra defined by a cocharacter $\lambda$, such that $X_K^n$ is contained in its unipotent radical, so that for some $a\in K$ the element $\lambda(a)\cdot X_K^n$ will extend to $R$ as well.
	
	Finally for any group $G$ we can consider a $z$-extension $$0 \to Z \to \widetilde{G} \to G \to 1$$  where $Z$ is a central torus and the derived group $G^\prime$ of $\widetilde{G}$ is simply connected. Then the map $\Bun_{\widetilde{G}} \to \Bun_G$ is a smooth surjection. Moreover the covering $G^\prime \to G$ is separable, because we assumed that the fundamental group of $G$ has no $p$-torsion. Therefore $\Lie(\widetilde{G})  \cong \Lie(Z) \oplus \Lie(G)$ and therefore the map $\Higgs_{\widetilde{G}}\to \Higgs_G$ also admits local sections. Thus it suffices to prove the result for $\widetilde{G}$.
\end{proof}

The semistable reduction theorem (\Cref{T:semistable_reduction}) now allows us to deduce:

\begin{cor} \label{C:higgs_properness}
	Suppose that the characteristic $p$ of $k$ is large enough such that Behrend's conjecture holds for $G$, such that $p$ is not a torsion prime for $G$ and such that $p$ is very good for $G$, then the Hitchin morphism
	$$ h \co \Higgs_G^{\ss} \to \cA_G$$
	satisfies the existence part of the valuative criterion for properness, i.e., if $R$ is a DVR with fraction field $K$ and 
	$x_K\co \Spec(K) \to \Higgs_G^{\ss}$ is a morphism such that $h(x_K)\co \Spec(K) \to \cA_G$ extends to $R$, then there exists an extension $R \to R'$ of DVRs with $K \to K'=\Frac(R')$ finite	
	and a morphism $x_{R^\prime} \co \Spec(R') \to \Higgs^{\ss}_G(R^\prime)$ extending $x_K$. 
\end{cor}
Note that in characteristic $0$ this result is due to Faltings \cite[Thm.~II.4]{FaltingsStableG} and for the regular part of the Hitchin fibration this is due to Chaudouard--Laumon \cite[Thm.~8.1.1]{ChaudouardLaumonLFP1}. Over the complex numbers, the result can also be deduced from results of Simpson as explained in \cite{deCataldo}.

\begin{remark}
Since $\Higgs_G^{\ss}$ is quasi-compact, the conclusion is equivalent to saying that the Hitchin morphism is universally closed. In particular the induced morphism on an adequate moduli space will be proper. 
\end{remark}

\begin{proof}
	By \Cref{L:valuativeHiggs} we can find an extension of $x_K$ to $x_R$. As $\Higgs_G$ admits a well ordered $\Theta$-stratification by Harder--Narasimhan reductions we can therefore apply the semistable reduction \Cref{T:semistable_reduction} to conclude.
\end{proof}


\newcommand{\base}{k}

\section{Good moduli spaces for objects in abelian categories}

In this section we study the moduli functor for objects in a $k$-linear abelian category $\cA$, following the foundational work of Artin and Zhang \cite{artin-zhang}, who explained that many of the results known for categories of quasi-coherent sheaves on a scheme can be carried out in an abstract setting.  This construction has been studied more recently in \cite{gaitsgory2005notion, abramovich2006sheaves, calabrese2013moduli}. The general setup is very useful as it for example gives an easy way to formulate moduli problems for objects in the heart of different t-structures in the derived category of coherent sheaves on a scheme. It turns out that this setup leads to moduli problems in which the conditions of $\Theta$-reductivity,  {\textsf{S}}-completeness, and unpunctured inertia can be checked rather easily. Following the convention of \cite{artin-zhang}, throughout this chapter we exceptionally fix a base ring $k$, which is allowed to be any commutative ring.

\subsection{Formulation of the moduli problem}
Let us start by recalling the setup of \cite{artin-zhang}. Let $\cA$ be a $k$-linear abelian category that is assumed to be cocomplete, i.e., arbitrary small colimits exist in $\cA$. To formulate a reasonable moduli problem we first need to recall some finiteness conditions on objects.
Recall that an object $E \in \cA$ is
\begin{itemize}
	\item \emph{finitely presentable} (also known as \emph{compact}) if the canonical map 
	\begin{equation} \label{E:colim-map}
	\colim_{\alpha \in I} \Hom(E,F_\alpha) \to \Hom(E,\colim_{\alpha \in I} F_\alpha)
	\end{equation}
	is an isomorphism for any (small) filtered system $\{F_{\alpha}\}_{\alpha \in I}$ in $\cA$;  
	\item \emph{finitely generated} if \eqref{E:colim-map}  is an isomorphism for any filtered system of monomorphisms in $\cA$, or equivalently, if $E = \bigcup_\alpha E_\alpha$ for a filtered system of subobjects, then $E = E_\alpha$ for some $\alpha$ \cite[Prop.~3.5.6]{popescu}; and
	\item \emph{noetherian} if every ascending chain of subobjects of $E$ terminates, or equivalently, if every subobject of $E$ is finitely generated.
\end{itemize}
We denote by $\cA^{\rm{fp}}$ the full subcategory of $\cA$ consisting of finitely presentable objects.

\begin{example}  If $\cA = \Mod_R$ for a (possibly non-commutative) ring $R$, an object $E \in \cA$ is finitely presentable (resp. finitely generated, resp. noetherian) if and only if the corresponding module is.  The analogous statement holds if $\cA = \QCoh(X)$ for a quasi-compact, quasi-separated scheme $X$. 
\end{example} 
We say that $\cA$ is
\begin{itemize}
	\item \emph{locally of finite type} if every object in $\cA$ is the union of its finitely generated subobjects; 
	\item \emph{locally finitely presented} if every object in $\cA$ can be written as the filtered colimit of finitely presentable objects, and $\cA^{\rm{fp}}$ is essentially small; 
	\item \emph{locally noetherian} if it has a set of noetherian generators.
\end{itemize}
If $\cA$ is locally noetherian, then finitely generated, finitely presentable, and noetherian objects coincide \cite[Prop.~B1.3]{artin-zhang}, and the category $\cA^{\rm{fp}}$ is closed under kernels and hence abelian. Our main results will assume that $\cA$ is locally noetherian.

The next ingredient to the formulation of moduli problems is the observation that the existence of colimits allows one to define a tensor product, which in turn provides a notion of base change.

More precisely, there is a canonical ${\base}$-bilinear functor 
\begin{equation} \label{E:tensor}
(-) \tensor_{\base} (-) \co \Mod_{\base} \times \cA \to \cA
\end{equation}
which is characterized by the formula $$\Hom_{\cA}(M \tensor_{\base} E,F) = \Hom_{\Mod_{\base}}(M,\Hom_{\cA}(E,F))$$
for objects $E,F \in \cA$ and a $\base$-module $M$.  Explicitly, if one presents $M$ as the cokernel of a morphism $\base^I \to \base^J$ for index sets $I$ and $J$, then $M \tensor_{\base} E$ can be computed as $\coker(E^I \to E^J)$ where the morphism $E^I \to E^J$ is induced by the matrix defining $\base^I \to \base^J$.

The functor $(-)\tensor_{\base} (-)$ commutes with filtered colimits and is right exact in each variable. If $M \in \Mod_{\base}$ is flat and $\cA$ is locally noetherian then $M \tensor_{\base}(-)$ is exact \cite[Lem.~C1.1]{artin-zhang}.

\begin{defn} \label{D:flatness} \cite[\S C1]{artin-zhang}
	We say that an object $E \in \cA$ is \emph{flat} if $(-) \tensor_{\base} E \co \Mod_{\base} \to \cA$ is exact.
\end{defn}

This tensor product leads to a base change formalism as follows.

\begin{defn}[Base change categories] \cite[\S B2]{artin-zhang}
For a commutative ${\base}$-algebra $R$, let $\cA_R$ denote the category of $R$-module objects in $\cA$, i.e., pairs $(E, \xi_E)$ where $E \in \cA$ and $\xi_E \co R \to \End_{\cA}(E)$ is a morphism of $\base$-algebras, and a morphism $(E, \xi_E) \to (E', \xi_{E'})$ in $\cA_R$ is a morphism $E \to E'$ in $\cA$ compatible with the actions of $\xi_E$ and $\xi_{E'}$.
\end{defn}

For a commutative ${\base}$-algebra $R$, $\cA_R$ is an $R$-linear abelian category \cite[Prop.~B2.2]{artin-zhang}, and $\cA_{\base} = \cA$. Given a homomorphism of commutative rings $\phi \co R_1 \to R_2$, the forgetful functor 
$$\phi_\ast \co \cA_{R_2} \to \cA_{R_1}$$ 
is faithfully exact, commutes with filtered colimits and faithful, and $\phi_\ast$ is fully faithful if $\phi$ is surjective \cite[Prop.~B2.3]{artin-zhang}. Moreover, $\phi_\ast$ admits a left adjoint 
$$\phi^*=R_2 \otimes_{R_1} (-) \co \cA_{R_1} \to \cA_{R_2}$$
by \cite[Prop.~B3.16]{artin-zhang}.

\begin{remark}
	The property of being locally noetherian is not stable under base change, but if $\cA$ is locally noetherian and $R$ is an essentially finite type $k$ algebra, then $\cA_R$ is locally noetherian \cite[Cor.~B6.3]{artin-zhang}.
\end{remark}	

The above constructions allow one to prove descent if $\cA$ is locally noetherian \cite[Thm.~C8.6]{artin-zhang}, i.e., if $R\to S$ is a faithfully flat map of commutative $\base$-algebras then $\cA_R$ is equivalent to the category of objects in $\cA_S$ equipped with a descent datum.

\begin{remark}[$\cA_{\cX}$ for stacks 	over $\base$]\label{R:AforStacks}
	As the assignment $R \mapsto \cA_R$ satisfies descent if $\cA$ is locally noetherian, it defines a stack $\underline{\cA}$ in the fppf topology on $\base$-$\alg$. As in \cite{lmb} this  extends the category $\cA_R$ naturally not only to schemes but also to algebraic stacks, i.e., for any algebraic stack $\cX$ over $\base$ we can define 
	$$\cA_{\cX} := \Map_{\text{Fibered Cat}/\base\text{-}\alg}(\cX,\underline{\cA}).$$  
If $\cX$ is the quotient stack for a groupoid of affine schemes $\cX= [X_1 \rightrightarrows X_0]$ with $X_i=\Spec(R_i)$, then descent implies that the category $\cA_{\cX}$ is naturally equivalent to the category of objects of $\cA_{X_0}$ equipped with a descent datum. We will use this description for the stacks $\Theta$ and $\oST_R$.
\end{remark}

Faithfully flat descent also allows one to extend the functor $R_2 \otimes_{R_1} (-) \co \cA_{R_1} \to \cA_{R_2}$ above to a functor $f^\ast \co \cA_{\cY} \to \cA_{\cX}$ for any morphism of stacks $f \co \cX \to \cY$. To prove extension theorems, we will need the following construction.
\begin{lem}[Push-forward in $\cA$]\label{R:ChechPushForward} 	
Suppose that $\cA$ is locally noetherian. If $f\co \cX \to \cY$ is a quasi-compact morphism with affine diagonal of algebraic stacks then the restriction functor $f^*\co \cA_\cY\to \cA_\cX$ admits a right adjoint $f_\ast$ which commutes with filtered colimits and flat base change.
\end{lem}	
\begin{proof}
Let us first prove the claim when $\cY = \Spec(A)$ is affine. In this case we can choose a presentation of $\cX$ by a groupoid in affine schemes $\cX \simeq [\Spec(R_1) \rightrightarrows \Spec(R_0)]$, and an object $E \in \cA_{\cX}$ is described by an object $E_0 \in \cA_{R_0}$ along with a descent datum, i.e., denoting by $d_i\co R_0\to R_1$ the structure maps of the presentation this is an isomorphism $\phi \co R_1 \otimes_{d_0,R_0} E_0 \simeq R_1 \otimes_{d_1,R_0} E_0$, satisfying a cocycle condition.
	Then using the fact that homomorphisms in $\cA_{[\Spec(R_1) \rightrightarrows \Spec(R_0)]}$ are homomorphisms in $\cA_{R_0}$ which commute with the respective cocycles, one may verify directly that
	\begin{equation}\label{E:pushforward}
	f_\ast(E) = \ker \left( E_0 \xrightarrow{\phi \circ \eta_0 - \eta_1} B_1 \otimes_{d_1,B_0} E_0 \right) \in \cA_A,
	\end{equation}
	where $\eta_i \co E_0 \to B_1 \otimes_{d_i,B_0} E_0$ for $i=0,1$ is the unit of the adjunction between $B_1 \otimes_{d_i,B_0} (-)$ and the forgetful functor $(d_i)_\ast \co\cA_{B_1} \to \cA_{B_0}$. Both objects in \eqref{E:pushforward} are regarded as objects of $\cA_A$ via the canonical forgetful functor. {The fact that $f_\ast$ commutes with filtered colimits can be deduced from the formula \eqref{E:pushforward} and the fact that filtered colimits are exact in $\cA_A$\cite[Prop.~B2.2]{artin-zhang}. The fact that $f_\ast$ commutes with flat base change can be deduced from the fact that if $A \to A'$ is a flat ring map, then $A'\otimes_A(-)$ is exact and hence commutes with the formation of the kernel in \eqref{E:pushforward}}.
	
Now let $\cY$ be an algebraic stack with affine diagonal. For any $E \in \cA_\cX$ and any morphism $\Spec(R) \to \cY$ we have shown that the object $(f_R)_\ast(E|_{\cX_R}) \in \cA_R$ is defined and its formation commutes with flat base change. Faithfully flat descent implies that these objects descend to a unique object $f_\ast(E) \in \cA_{\cY}$. Faithfully flat descent also shows that the resulting functor $f_\ast \co\cA_\cX \to \cA_\cY$ is right adjoint to $f^\ast$.
\end{proof}	

\begin{rem} \label{R:Cech_pushforward}
	If $X$ is a separated scheme and $f\co X \to \Spec(R)$ a morphism then the above construction reproduces the usual Cech description of the push forward, i.e., given $E\in \cA_X$ choose a covering $X_0 = \bigsqcup_i U_i \to X$ by open affines then \eqref{E:pushforward} reduces to $$f_\ast(E) = \ker \bigg( \bigoplus_i (f_i)_\ast(E|_{U_i}) \to \bigoplus_{i<j} (f_{ij})_\ast(E|_{U_{ij}}) \bigg),$$ where $U_{ij}:= U_i \cap U_j$.
\end{rem}

\begin{defn}[Moduli functor]\label{D:abelian_moduli}
	Let ${\base}$ be a commutative ring and let $\cA$ be a locally noetherian, cocomplete, and ${\base}$-linear abelian category. Then we define the category $\cM_{\cA}$ fibered in groupoids over $\base$-$\text{alg}$ by assigning the groupoid
	\[
	\cM_\cA(R) := \langle \text{objects } E \in \cA_R \text{ which are flat and finitely presented} \rangle 
	\]
	for a $\base$-algebra $R$.
\end{defn}

\begin{lem}
	The category fibered in groupoids $\cM_\cA$ is a stack in the big fppf topology on $\base$-${\rm alg}$ and extends naturally to a stack on the big fppf topology on schemes over $\base$.
\end{lem}
\begin{proof}
	As we already quoted the result \cite[Thm.~C8.6]{artin-zhang} that the categories $\cA_R$ satisfy flat descent, we only need to check that the conditions of flatness and finite presentation are preserved by descent and pull-back. Given a ring map $R_1 \to R_2$, the pullback functor $R_2 \otimes_{R_1} (-)$ preserves flat objects by \cite[Lem.~C1.2]{artin-zhang}, and it preserves finitely presentable objects because its right adjoint  commutes with filtered colimits. If $R_1\to R_2$ is a faithfully flat map of $\base$-algebras we have $R_2 \otimes_{R_1}(M \otimes_{R_1} (-)) \simeq (R_2 \otimes_{R_1} M) \otimes_{R_2} (R_2 \otimes_{R_1}(-))$, and the exactness of a sequence in $\cA_{R_1}$ can be checked after applying the functor $R_2 \otimes_{R_1}(-)$, so $E \in \cA_{R_1}$ is flat if $R_2 \otimes_{R_1} E$ is. Also, one can directly verify from the description of $\Hom$ in the category of descent data for the map $R_1 \to R_2$ that any descent data for a finitely presentable object $E \in \cA_{R_2}$ is a finitely presentable object in the category of descent data.
\end{proof}

\begin{rem}
If $\cA$ is a locally noetherian $k$-linear abelian category and $R$ is a $k$-algebra, most finitely presented objects in $\cA_R$ will have an underlying object in $\cA$ that is not finitely presented. This is why it is essential that $\cA$ be cocomplete, and thus rather big, in our definition.
\end{rem}

\begin{warn}[Choice of cocompletion matters]
Different algebraic moduli functors can arise from embedding the same small category into different cocomplete categories. Let $X$ be a smooth projective scheme over $\bC$. The natural moduli functor for finite dimensional representations of the fundamental group $\pi_1(X)$ uses $R$-module objects in the category $\cA$ of \emph{all} representations of $\pi_1(X)$, i.e., modules over the group algebra $R[\pi_1(X)]$. Because $\cA^{\rm fp}$ is not Hom-finite, to get an algebraic moduli functor one must define a family over $\Spec(R)$ to be a flat and finitely presented $R[\pi_1(X)]$-module whose underlying $R$-module is also finitely presented. (If one naively applies \Cref{D:abelian_moduli} to the formal Ind-completion of the category of finite dimensional $\pi_1(X)$-representations, the resulting moduli functor is not algebraic.)

On the other hand, the category of \emph{finite dimensional} $\pi_1(X)$-modules is equivalent to the category of vector bundles with flat connection. The latter admits a different Ind-completion as quasi-coherent $\cD_X$-modules. The natural algebraic moduli functor from this perspective defines a family over $\Spec(R)$ to be a quasi-coherent $R \otimes \cD_X$-module on $X$ whose underlying quasi-coherent sheaf on $X_R$ is $R$-flat and coherent. Closed points of this stack are in bijection with those of the moduli of finite dimensional $\pi_1(X)$ representations. In fact, the non-abelian Hodge correspondence implies that these stacks have homeomorphic good moduli spaces \cite[Sect.~7]{simpson}, but they are not algebraically isomorphic.

\end{warn}

\subsection{Verification of the valuative criteria for the stack \texorpdfstring{$\cM_\cA$}{MA}}
To apply our existence results we will check $\Theta$-reductivity and  {\textsf{S}}-completeness for the stack $\cM_{\cA}$ with respect to DVRs which are essentially of finite type over $k$. The first step is to show that, as for module categories, $\bG_m$-equivariant objects can be interpreted as graded objects. Let us recall these notions.

Recall from \cite[\S B7]{artin-zhang} that the category of $\bZ$-graded objects $\cA^{\bZ}$ consists of functors $\bZ \to \cA$, where $\bZ$ is regarded as a category with only identity morphisms. Concretely, we think of objects in $\cA^{\bZ}$ as objects $E \in \cA$ equipped with the extra data of a direct sum decomposition $E \cong \bigoplus_{n \in \bZ} E_n$. Given a $\bZ$-graded ${\base}$-algebra $A$, \emph{a $\bZ$-graded $A$-module object} is an object of $\cA^{\bZ}$ whose underlying object $E = \bigoplus_{n \in \bZ} E_n \in \cA$ is equipped with an $A$-module structure such that  multiplication $A \tensor_\base E \to E$ maps $A_n \tensor_\base E_m$ to $E_{n+m}$.
The category $\cA_A^{\bZ}$ of $\bZ$-graded $A$-module objects is abelian and locally noetherian if $\cA_A$ is \cite[Prop.~B7.5]{artin-zhang}. 

Now we encode the $\bZ$-grading of a graded ${\base}$-algebra $A$ by a morphism of $\base$-algebras $\sigma_A : A \to A[t^{\pm 1}]$ that equips $A$ with the structure of a comodule over the coalgebra $k[t^{\pm 1}]$. Concretely, $\sigma_A(a)=\sum a_n t^n$, where $a_n$ is the $n^{th}$ graded piece of $a \in A$. Objects of $\cA_{[\Spec(A)/\bG_m]}$ are by definition  objects $E \in \cA_A$ together with a cocycle, which can be encoded by a coaction morphism \[\sigma : E \to E[t^{\pm1}] := A[t^{\pm 1}] \tensor_A E.\]
This is a morphism of $A$-module objects, where $E[t^{\pm1}]$ is equipped with the twisted $A$-module structure induced by $\sigma_A$, i.e., $a \in A$ acts by $\sum_n a_n t^n$. Note that the object in $\cA$ underlying $E[t^{\pm1}]$ is a direct sum of one copy of $E$ for each monomial $t^n$, so one can write $\sigma = \sum_n \sigma_n t^n$, where $\sigma_n : E \to E$ are morphisms in $\cA = \cA_k$. Rather than being $A$-linear, one has $\sigma_n(a \cdot -) = \sum_{i+j=n} a_i \cdot \sigma_j(-)$ for any $a \in A$ with homogeneous pieces $a_i$.

The cocycle condition on $\sigma$ amounts to the condition that the following diagrams in $\cA$ must commute:
\begin{equation} \label{E:coaction}
\xymatrix{
	E \ar[r]^-{\sum \sigma_n t^n} \ar[d]^{\sum \sigma_n t^n} & E[t^{\pm 1}] \ar[d]^{t \mapsto t t'} \\
	E[t^{\pm1}] \ar[r]^-{\sum \sigma_n (t')^n} & E[t^{\pm1},(t')^{\pm1}]
}
\qquad
\xymatrix{
	E \ar[r]^-{\sum \sigma_n t^n} \ar[rd]_{\id}					& E[t^{\pm1}] \ar[d]^{t \mapsto 1}\\
	&  E
}
\end{equation}
One can show that the morphisms in this diagram are in fact $A$-linear, where $a \in A$ acts on $E[t^{\pm1},(t')^{\pm1}]$ by $\sum_n a_n (tt')^n$.

\begin{prop} \label{P:abelian_graded_alg}
Let ${\base}$ be a commutative ring and let $\cA$ be a locally noetherian ${\base}$-linear abelian category.  Let $A$ be a $\bZ$-graded ${\base}$-algebra. Then there is a natural equivalence $\cA^\bZ_A \to \cA_{[\Spec(A)/\bG_m]}$ that maps $E \in \cA_A^\bZ$ to the object of $\cA_{[\Spec(A)/\bG_m]}$ defined by the coaction morphism $\sigma = \sum \sigma_n t^n :  E \to E[t^{\pm1}]$, where $\sigma_n : E \to E$ is the $k$-linear projection onto the weight $n$ summand of $E$. This restricts to an equivalence between $\cM_{\cA}([\Spec(A)/\bG_m])$ and the groupoid of objects in $\cA_A^{\bZ}$ whose underlying non-graded $A$-module object is flat and finitely presented.
\end{prop}
\begin{proof}
Given a coaction $\sigma : E \to E[t^{\pm1}]$ as discussed above, the commutativity of \eqref{E:coaction} is equivalent to the identities $\sum_{m,n} \sigma_m \sigma_n t^m (t')^n = \sum_n \sigma_n (tt')^n$ and $\sum_n \sigma_n = \id$. As in the case of graded $\base$-modules, this implies that $\sigma_n$ are a collection of mutually orthogonal idempotent endomorphisms of $E$ that induce a direct sum decomposition $E = \bigoplus_n E_n$ in $\cA$, where $E_n$ is the image of $\sigma_n$. Conversely for any grading $E \cong \bigoplus_n E_n$, one can define a coaction by letting $\sigma_n$ be the projection onto $E_n \subset E$. As discussed above, the identity $\sigma_n a = \sum_{i+j=n} a_i \sigma_j$ is equivalent to the fact that the coaction $\sigma$ is $A$-linear and thus that it defines a descent datum for $\cA_{[\Spec(A)/\bG_m]}$. We leave it to the reader to verify that this identity is also equivalent to the compatibility of the grading of $E$ and $A$ in the sense that $A_n \otimes E_m$ maps to $E_{m+n}$.

The claim that $E \in \cA_A^\bZ$ corresponds to an element of $\cM_\cA([\Spec(A)/\bG_m])$ if and only if the underlying object in $\cA_A$ is flat and finitely presented follows from the fact that $\cM_\cA$ is a stack for the fppf topology.
\end{proof}

We can use the above result to describe objects of $\cA$ over $\Theta_R = [\Spec(R[x])/\bG_m]$ for any $\base$-algebra $R$ and over $\oST_R := [\Spec \big(R[s,t] / (st-\pi)\big)  / \bG_m]$ for any DVR $R$ over $\base$ with uniformizing parameter $\pi$.  Both descriptions are in terms of filtered objects.

\begin{cor}\label{C:reese-ThetaR} Suppose that $\cA$ is locally noetherian. Let $R$ be a $\base$-algebra then the category $\cA_{\Theta_R}$ is equivalent to the category of sequences of morphisms
	$$E \co \quad \cdots \to E_{n+1} \xrightarrow{x} E_{n} \to \cdots$$
	in $\cA_R$, such that the restriction of $E$
	\begin{itemize}
		\item along $\Spec(R) \hookarr \Theta_R$ is $\colim E_i$, and
		\item along $B \bG_{m,R} \hookarr \Theta_R$ is $\bigoplus_{n \in \bZ}  E_{n}/x(E_{n+1})$.
	\end{itemize}
	This equivalence restricts to an equivalence between $\cM_{\cA}(\Theta_R)$ and the groupoid of $\bZ$-weighted filtrations $\cdots \subset E_{n+1} \subset E_{n} \subset \cdots$ of an object $E_{\infty}$ in $\cA_R$ such that $E_{n}/E_{n+1} \in \cA_R$ is flat and finitely presented, $E_n = E_{\infty}$ for $n \ll 0$ and $E_n = 0 $ for $n \gg 0$.
\end{cor}	

\begin{proof}  
	By \Cref{P:abelian_graded_alg}, an object $E$ in $\cA_{\Theta_R}$ is a $\bZ$-graded $R[x]$-module object of $\cA_R$.  This corresponds to a $\bZ$-graded object $E = \bigoplus_{n \in \bZ} E_n$ in $\cA_R$ together with a multiplication $x \co E \to E$ mapping $E_{n+1}$ to $E_{n}$.
	The restriction of $E$ along the open immersion $\Spec(R) \hookarr \Theta_R$ is the $\bG_m$-invariants (i.e. the degree $0$ component) of the $\bZ$-graded $R[x^{\pm1}]$-module object $E \tensor_{R[x]} R[x^{\pm1}] \in \cA_{R[x^{\pm1}]}$.   We compute that
	$$\begin{aligned}
	E \tensor_{R[x]} R[x^{\pm1}]  & = E \tensor_{R[x]} (\colim( \cdots \xrightarrow{x} R{[x]} \xrightarrow{x} R{[x]} \xrightarrow{x} \cdots) \\
	& = \colim ( \cdots \xrightarrow{x} E \xrightarrow{x} E \xrightarrow{x} \cdots) \\
	& = \bigoplus_{n \in \bZ} \colim (\cdots \xrightarrow{x} E_{n+1} \xrightarrow{x} E_{n} \xrightarrow{x} \cdots)
	\end{aligned}$$
	whose $\bG_m$-invariants is $\colim E_n$. The restriction of $E$ along the closed immersion $B \bG_{m,R} \hookarr \Theta_R$ is the object $E \tensor_{R[x]} (R[x]/x) \cong E/xE$ in $\cA_R$ with the $\bZ$-grading $E/xE = \bigoplus_{n \in \bZ} E_{n} / x(E_{n+1})$. This proves the first claim.
	
	If $E  \in \cA_{\Theta_R}$ is flat and finitely presented, then so is the corresponding object in $\cA_{R[x]}$ (also denoted by $E$).  Since $E \in \cA_{R[x]}$ is flat, it is torsion free and thus $x \co E \to E$ is injective.  The base change $E \tensor_{R[x]} (R[x]/x)= \bigoplus E_n / E_{n+1}$ is a flat and finitely presented object in $\cA_R$ which implies that the sum is finite and the summands are finitely presentable and flat. This implies that $E_n$ stabilizes for $n \ll 0$. Also, if $E \in \cA_{R[x]}$ is finitely presentable, it must admit a surjection $R[x] \otimes_R F \to E$ for some $F \in \cA_R^{\rm{fp}}$, corresponding to a map $F \to E$ in $\cA_R$.  Because $F$ is finitely presentable it must factor through $\bigoplus_{n \leq N} E_n \subset E$ for some $N$, and hence the image of $R[x] \otimes_R F \to E$ lies in  $\bigoplus_{n \leq N} E_n$ as well, which implies that $E_n = 0$ for $n \gg 0$.
	
	Conversely, if $\cdots \subset E_{n+1} \subset E_n \subset \cdots$ satisfies the conditions above, then each $E_n$ is constructed as a finite sequence of extensions of flat and finitely presentable objects in $\cA_R$ and is thus flat and finitely presentable. It follows that the graded $R[x]$-module objects $R[x] \otimes_R E_n$ are flat and finitely presentable as objects of $\cA_{R[x]}$. Furthermore, the object $E \in \cA_{R[x]}^\bZ$ corresponding to this $\bZ$-weighted filtration 
	can be constructed as a finite sequence of extensions of objects of the form $R[x] \otimes_R E_n\langle-n\rangle$, where the $\langle -n \rangle$ denotes a grading shift so that the resulting object is homogeneous of degree $n$. Hence $E \in \cA_{\Theta_R}$ is finitely presentable and flat.
\end{proof}

\begin{cor} \label{C:reese}  Suppose that $\cA$ is locally noetherian.  Let $R$ be a DVR over $\base$ with uniformizing parameter $\pi$ and residue field $\kappa$.
	The category $\cA_{\oST_R}$ is equivalent to the category of diagrams in $\cA_R$
	\begin{equation} \label{E:abelian_str}
	\xymatrix{E \co & \cdots \ar@{}|-{}[r] \ar@/^/[r]^s & E_{n-1} \ar@{}|-{}[r] \ar@/^/[r]^s \ar@/^/[l]^t & E_{n} \ar@{}|-{}[r] \ar@/^/[r]^s \ar@/^/[l]^t & E_{n+1}  \ar@/^/[r]^s \ar@/^/[l]^t \ar@{}|-{}[r] & \ar@/^/[l]^t \cdots },
	\end{equation}
	satisfying $st=ts=\pi$. Under this equivalence the restriction of $E$
	\begin{itemize}
		\item along $\Spec(R) \xhookrightarrow{s \neq 0} \oST_R$ is
		$\colim (\cdots \xrightarrow{s} E_{n-1} \xrightarrow{s} E_{n} \xrightarrow{s} \cdots)$,
		\item along $\Spec(R) \xhookrightarrow{t \neq 0} \oST_R$ is
		$\colim (\cdots \xleftarrow{t} E_{n-1} \xleftarrow{t} E_{n} \xleftarrow{t} \cdots)$,
		\item along $\Theta_{\kappa} \xhookrightarrow{s = 0} \oST_R$ is the object corresponding to the sequence $(\cdots \xleftarrow{t} E_n/sE_{n-1} \xleftarrow{t} E_{n+1}/sE_n \xleftarrow{t} \cdots)$, and
		\item along $\Theta_{\kappa} \xhookrightarrow{t = 0} \oST_R$ is $(\cdots \xrightarrow{s} E_{n-1}/tE_{n} \xrightarrow{s} E_{n}/tE_{n+1} \xrightarrow{s} \cdots)$.
	\end{itemize}
	This equivalence restricts to an equivalence between $\cM_{\cA}(\oST_R)$ and the groupoid consisting of objects $E$ such that:  (a) $s$ and $t$ are injective, (b) $s \co E_{n-1}/tE_n \to E_{n}/t E_{n+1}$ is injective for all $n$, (c) each $E_n \in \cA_R$ is finitely presentable, (d) $s \co E_{n-1} \to E_n$ is an isomorphism for $n \gg 0$, and (e) $t \co E_{n} \to E_{n-1}$ is an isomorphism for $n \ll 0$.
\end{cor}

\begin{proof}  
	The equivalence between $\cA_{\oST_R}$ and the category of diagrams as in \eqref{E:abelian_str} is argued as before (\Cref{C:reese-ThetaR}). 
	Let us first show that flatness is characterized by conditions (a) and (b). 
	Suppose $E \in \cA_{\oST_R}$ is flat, then the pullbacks of $E$ to $R[s,t]/(st-\pi)$ and $\kappa[s]$ (by setting $t=0$) are both flat and in particular torsion free which gives conditions (a) and (b). Conversely if $E$ is given by the diagram \eqref{E:abelian_str} flatness is a local condition and (a) implies that the restriction of $E$ to $s\neq 0$ (or $t\neq 0$) is torsion free and thus flat. Condition (b) implies that the restriction to $\kappa[s]$ is $s-$torsion free and so $E$ is also flat at the origin $s=0=t$ by applying  \cite[Lem.~C1.12]{artin-zhang}.
	
	To check that a finitely presentable, flat object $E$ satisfies conditions (c)--(e) note that these are closed under cokernels, so we only need to check these for a generating class of objects.
	Now if $M \in \cA_R$ is an $R$-module, then
	\[
	R[s,t]/(st-\pi) \otimes_R M \simeq \bigoplus_{n<0} M \cdot t^{-n} \oplus M \oplus \bigoplus_{n>0} M \cdot s^n,
	\]
	is a $\bZ$-graded $R[s,t]/(st-\pi)$-module for which $t$ is an isomorphism on negatively graded pieces, and $s$ is an isomorphism on positively graded pieces. Therefore for any finitely presentable $M \in \cA_R$, $R[s,t]/(st-\pi) \otimes_R M\in \cA_R^{\bZ}$ satisfies the conditions (c)--(e) of the lemma, and the same holds if $M \in \cA_R^\bZ$ is graded object with finitely many non-trivial graded pieces which are each finitely presentable. {As any finitely generated object of $\cA^\bZ_{R[s,t]/(st-\pi)}$ admits a surjection from an object of this form,} any $E \in (\cA^\bZ_{R[s,t]/(st-\pi)})^{\rm{fp}}$ admits a presentation of the form $E \simeq \coker(R[s,t]/(st-\pi) \otimes M_1 \to R[s,t]/(st-\pi) \otimes M_0)$ for some $M_0,M_1 \in (\cA^\bZ_R)^{\rm{fp}}$, which proves (c)--(e) for $E$.

Conversely, suppose that the diagram \eqref{E:abelian_str} satisfies the conditions (c)--(e) of the lemma. Then (d) and (e) imply that for $N \gg 0$, the canonical homomorphism of graded $R[s,t]/(st-\pi)$-modules $R[s,t] / (st-\pi) \otimes_R \bigoplus_{-N\leq n \leq N} E_n \to E$ is surjective. Let $K = \bigoplus_n K_n$ be the kernel of this homomorphism. Because $R[s,t] / (st-\pi) \otimes_R \bigoplus_{-N\leq n \leq N} E_n$ satisfies conditions (d) and (e), it follows that $K$ satisfies these conditions as well. Therefore $R[s,t]/(st-\pi) \otimes_R \bigoplus_{-M \leq n \leq M} K_n \to K$ is surjective as well for $M \gg 0$. By (c) each $K_n$ is the kernel of a surjection of finitely presented objects, and it is thus finitely generated. We have thus expressed $E$ as the cokernel
\[
E = \coker \big(R[s,t]/(st-\pi) \otimes_R \bigoplus_{|n| \leq M} K_n \to R[s,t] / (st-\pi) \otimes_R \bigoplus_{|n| \leq N} E_n \big)
\]
of a homomorphism from a finitely generated object to a finitely presented object, so $E$ is finitely presented.
\end{proof}

Now our extension results will follow from a basic result about extensions in codimension 2 that again carries over from quasi-coherent modules.

\begin{lem} \label{L:abelian_flatness}
	Let $j \co U \to X$ be an open subscheme of a regular noetherian scheme of dimension $2$ whose complement is $0$-dimensional. Then $j_\ast \co\cA_{U} \to \cA_X$ maps flat objects to flat objects, and induces an equivalence between the full subcategory of flat objects over $X$ and over $U$, with inverse given by $j^\ast \co \cA_X \to \cA_U$.
\end{lem}
\begin{proof}
	It suffices to show that $j_\ast$ preserves flat objects, and that both the unit and counit of the adjunction between $j^\ast$ and $j_\ast$ are equivalences on flat objects. The property of being an isomorphism is local by descent, so we may assume that $X = \Spec(R)$ is affine and $U$ is the complement of a single closed point. Localizing further it suffices to consider the case of $X = \Spec(R)$ for a regular ring $R$ of dimension $2$ and $U \subset X$ the complement of the closed point $p$ whose maximal ideal is generated by a regular sequence $x,y \in R$. In particular $U=\Spec(R_{x}) \cup \Spec(R_{y})$.
	
	By construction (\Cref{R:ChechPushForward}) we then have
	$$ j_*(E) = \ker( E|_{R_{x}} \oplus E|_{R_{y}}\to E|_{R_{xy}}).$$
	So the natural map $j^\ast(j_\ast(E)) \to E$ is an equivalence for any $E$ as $\cA_{U}$ was defined by descent from affine schemes. 
	
	Conversely if $E\in \cA_R$ is flat we can tensor $E$ with the left exact sequence 
	$$ 0 \to R \to R_{x}\oplus R_{y} \to R_{xy}$$ 
	to find that the sequence 
	$$ 0 \to E \to  E|_{R_{x}} \oplus E|_{R_{y}}\to E|_{R_{xy}}$$ is still exact, so $E \cong j_*(E|_U)$.
	
	Finally we must show that $j_\ast$ preserves flat objects. By \cite[Lem.~C1.12]{artin-zhang} we only need to show that $$\Tor_1^R(R/\cp,j_*E)=0$$ for all prime ideals $\cp$ of $R$. (Here $\Tor$ is defined as usual by choosing projective resolutions of $R$-modules.) 
	For prime ideals in $U$ this follows from flatness of $E$, so it suffices to show that $\Tor_1(\kappa, j_\ast(E)) = 0$ for any flat $E \in \cA_U$, where $\kappa=R/(x,y)$ is the residue field of the missing point, i.e., we need to show that tensoring $j_*E$ with the Koszul complex
	$$ 0 \to R \to R\oplus R \to R \to \kappa \to 0$$
	gives an exact sequence
	\[
	0 \to j_\ast(E) \xrightarrow{x \oplus (-y) } j_\ast(E) \oplus j_\ast(E) \xrightarrow{y \oplus x} j_\ast(E).
	\]
	This is the pushforward of the tensor product of $E$ with the short exact sequence of flat objects on $U$, $0 \to \cO_U \to \cO_U \oplus \cO_U \to \cO_U \to 0$ and therefore exact, because $j_\ast$ is left exact.
\end{proof}

This construction now allows us to check our conditions for the existence of good moduli spaces for $\cM_\cA$. In the following we will assume that $\cA$ is locally noetherian and use the result \cite[Cor.~B6.3]{artin-zhang} stating that then for any $\base$-algebra that is essentially of finite type (i.e., a localization of a finitely generated $\base$-algebra) the category $\cA_R$ is again locally noetherian and thus finitely presentable and noetherian are equivalent notions for objects of $\cA_R$.

We start with {\textsf{S}}-completeness.

\begin{lem} \label{L:abelian_S_complete}
	If $\cA$ is locally noetherian, then $\cM_\cA$ is  {\textsf{S}}-complete with respect to any DVR $R$ that is essentially of finite type over ${\base}$.
\end{lem}
\begin{proof}
Let us denote by $j\co \oST_R \smallsetminus 0 \subset \oST_R$ the inclusion and take any $E\in \cM_{\cA}(\oST_R \smallsetminus 0)$. By \Cref{L:abelian_flatness}, $j_\ast(E)$ is flat, so we only need to show that it is finitely presentable, i.e., we have to check conditions (c)--(e) of \Cref{C:reese}.
	
	We shall compute $j_\ast(E)$ explicitly. Let us denote by $K$ the fraction field of $R$. As $E$ is flat, it is defined by an object $F \in \cA_{K}$ and two $R$-module subobjects $E_1,E_2 \subset F$ such that $K\otimes_R E_i \to F$ is an isomorphism. 
	
	Let $j_i \co \Spec(R) \to \oST_R$ for $i=1,2$ denote the two open immersions and $j_{12} \co \Spec(K) \to \oST_R$ their intersection, then  
	\[
	j_\ast(E) = \ker\left( (j_1)_\ast (E_1) \oplus (j_2)_\ast(E_2) \to (j_{12})_\ast(F) \right).
	\]
	Let us describe this kernel as a graded $R[s,t]/(st-\pi)$-module. In terms of  graded modules, $j_1$ is given by the graded inclusion $R[s,t]/(st-\pi) \subset R[t^{\pm 1}]$ and $j_2$ by the graded inclusion $R[s,t]/(st-\pi) \subset R[s^{\pm 1}]$.
	
	As the maps $E_i \to F$ are injective, we may thus identify $j_\ast(E) \subset (j_{12})_\ast(F)$ with the intersection of two subobjects $(j_1)_*(E_1)$ and $(j_2)_*(E_2)$ under the equivalence of \Cref{P:abelian_graded_alg}, we compute
	\begin{gather*}
	(j_{12})_\ast(F) = R[t^{\pm 1}] \otimes_R F = \bigoplus_n F t^n, \\
	(j_1)_\ast(E_1) = E_1 \otimes_R R[t^{\pm 1}]= \bigoplus_n E_1 t^n, \\
	(j_2)_\ast(E_2) = E_2 \otimes_R R[s^{\pm 1}] \simeq \bigoplus_{n \in \bZ} (\pi^{-n} \cdot E_2) t^{n} \subset (j_{12})_\ast(F)
	\end{gather*}
	where in the third line we have used the identification $s = t^{-1} \pi$. We compute that:
	$$j_\ast(E) \simeq \bigoplus_{n \in \bZ} \big( E_1 \cap (\pi^{-n} \cdot E_2) \big) t^n \subset \bigoplus_{n\in \bZ} F t^n.$$
	
	Now each graded piece of $j_\ast(E)$ is finitely presentable because they are subobjects of the noetherian object $E_1$ (using that $\cA_R$ is locally noetherian). The union of the ascending sequence $\cdots \subset E_1 \cap (\pi^{-n} \cdot E_2) \subset E_1 \cap (\pi^{-n-1} \cdot E_2) \subset \cdots $ is $E_1$ because $K \otimes_R E_2 \simeq F$, and because $E_1$ is finitely generated, this union must stabilize.  By symmetry the same argument applies to $E_2$ and thus $j_\ast(E)$ is finitely presentable.
\end{proof}

Next we show $\Theta$-reductivity.
\begin{lem} \label{L:abelian_Theta_reductive}
	If $\cA$ is locally noetherian, then $\cM_\cA$ is $\Theta$-reductive with respect to any DVR $R$ that is essentially of finite type over ${\base}$.
\end{lem}
\begin{proof}
	As in the proof of \Cref{L:abelian_S_complete} let us denote by  $K$ the fraction field of $R$, $j \co \cU \hookarr \Theta_R$
 	 the complement of the closed point and we take $E \in \cM_{\cA}(\cU)$. Then by \Cref{L:abelian_flatness}, $j_\ast(E)$ is flat and so we need to show that it is finitely presentable.
	
	We pass to the presentation $\bA^1_R\to \Theta_R$, where the open subset $U \subset \bA^1_R$ corresponding to $\cU$ is covered by the two affine subschemes defined by $R[x] \subset K[x]$ and $R[x] \subset R[x^{\pm1}]$. Now $E \in \cA_\cU$ corresponds to an object $F \in \cA(K)$, a $R$-submodule object $E_1 \subset F$ such that $K \otimes_R E_1 \to F$ is an isomorphism, and a weighted descending filtration $\cdots F_{n+1} \subset F_n \subset \cdots \subset F$ satisfying the hypotheses of \Cref{C:reese-ThetaR}, then $j_\ast(E)$ corresponds to the graded $R[x]$-module object
	\[
	j_\ast(E) = \bigoplus_{n \in \bZ} \big( F_{n} \cap E_1 \big) x^{-n} \subset \bigoplus_{n \in \bZ} F x^{-n} = (\Spec(K) \to \Theta_R)_\ast(F). 
	\]
	Because $\cA_R$ is locally noetherian each graded piece $G_n := F_n \cap E_1$ of $j_\ast(E)$ will be finitely presentable, the maps $G_{n+1} \to G_n$ are injective, $G_n = 0$ for $n \gg0$, and $G_n = E_1$ for $n \ll 0$. Thus $j_\ast(E)$ is finitely presentable by \Cref{C:reese}.
\end{proof}

The valuative criteria for universal closedness turn out to be satisfied as well:
\begin{lem} \label{L:abelian_properness}
	If $\cA$ is a locally noetherian abelian category, the stack $\cM_\cA$ satisfies the valuative criterion for universal closedness, i.e., the existence part of the valuative criterion for properness, with respect to DVRs which are essentially of finite type over $\base$.
\end{lem}
\begin{proof}
If $R$ is a DVR, the statement that an object $E\in \cA_R$ is flat if and only if it is torsion free follows from  \cite[Lem.~C1.12]{artin-zhang} as the condition is equivalent to the vanishing of $\Tor_1$. If $j \co \Spec(K) \to \Spec(R)$ denotes the inclusion of the generic point, then for any $E \in \cA_{K}$, we can write $j_\ast(E) = \bigcup_\alpha F_\alpha$ as a directed union of finitely generated (hence finitely presentable) subobjects which must be torsion free. If $E$ is finitely generated then $E = \bigcup_\alpha K \otimes_R F_\alpha$ must stabilize, so there is some flat and finitely presentable object $F_\alpha$ extending $E$.
\end{proof}

We will also use the following below:

\begin{lem} \label{L:abelian_semi_simple}
Suppose that $\cA$ is locally noetherian.  If $\cM_\cA$ is an algebraic stack with affine stabilizers, $\kappa$ is a field over ${\base}$, and $E \in \cM_{\cA}(\kappa)$ represents a closed point, then $E$ is a semisimple object in $\cA_{\kappa}$. 
\end{lem}
\begin{proof}
	Because $E$ is finitely presented, it can not be expressed as an infinite sum of non-zero objects. Therefore, we only have to show that every finite filtration of $E$ splits. Now by \Cref{C:reese} any finite filtration of $E$ corresponds to a map $\Theta_{\kappa} \to \cM_{\cA}$ mapping $1 \mapsto E$. Because $E$ is a closed point, the resulting map must factor through a map $\Theta_{\kappa} \to B_{\kappa}\Aut_{\cM_{\cA}} (E)$. 
	We know from the classification of torsors on $\Theta_{\kappa}$ (\cite[Lem.~1.7]{heinloth-stability} or \cite[Prop.~A.0.1]{hlinstability}) that any such map factors through the projection $\Theta_{\kappa} \to B_{\kappa}\bG_m$, and thus the corresponding filtration of $E$ is split.
\end{proof}

\begin{lem} \label{L:abelian_diagonal}
If $\cM_\cA$ is an algebraic stack locally of finite presentation over $\base$, then the diagonal of $\cM_\cA$ is affine.
\end{lem}
\begin{proof}
If $R$ is a valuation ring over $\base$ with fraction field $K$ and $E,F \in \cM_\cA(R)$, then $F \to K \otimes_R F$ is injective and hence so is the restriction map $\Hom_R(E,F) \to \Hom_K(K \otimes_R E , K \otimes_R F) \simeq \Hom_R(E,K \otimes_R F)$. This implies the valuative criterion for separatedness of the diagonal of $\cM_\cA$.

For any ring $R$ over $\base$ and $E,F \in \cM_\cA(R)$, we claim that the functor $R'/R \mapsto \Hom_{R'}(R'\otimes_R E,R'\otimes_R F)$ is represented by a separated algebraic space $\uHom_R(E,F)$ locally of finite presentation over $R$. First, observe that the subfunctor $P \subset \Aut_R(E\oplus F)$ classifying automorphisms of the form $\left[ \begin{smallmatrix} A & 0 \\ C & D \end{smallmatrix} \right]$ is representable by a closed subspace, because it is the preimage of the (closed) identity section under the map of separated $R$-spaces $\Aut_R(E\oplus F) \to \Aut_R(E\oplus F)$ given by
\[
\begin{bmatrix} A & B \\ C & D \end{bmatrix} \mapsto \begin{bmatrix} 1 & B \\ 0 & 1 \end{bmatrix}.
\]
Next observe that we have a group homomorphism $P \to \Aut_R(E) \times \Aut_R(F)$ over $R$ given by
\[
\begin{bmatrix} A & 0 \\ C & D \end{bmatrix} \mapsto (A,D).
\]
The preimage of the (closed) identity section is the subgroup classifying automorphisms of the form $\left[ \begin{smallmatrix} 1 & 0 \\ C & 1 \end{smallmatrix} \right]$, which is canonically identified with the functor $\uHom(E,F)$. Thus $\uHom(E,F)$ is a closed subgroup of $\Aut_R(E\oplus F)$, and it is representable, separated, and locally finitely presented over $R$.

From the functor of points definition of the algebraic space $X := \uHom_R(E,F)$, there is a natural action of $\bG_m$ which scales the homomorphism. Furthermore, the resulting map $\bG_m \times X \to X$ extends (uniquely) to $\bA^1 \times X$, i.e., $X = X^+$ in the terminology of \Cref{S:drinfeld}. Under these hypotheses, the morphism $X^+ \to X^{\bG_m}$ is affine \cite{ahr2}. Also, the fixed locus $X^{\bG_m}$ is the zero section $X^{\bG_m} \simeq \Spec(R) \hookrightarrow X$, and hence $X = X^+ = \uHom_R(E,F)$ is affine as well.

The algebraic $R$-space $\op{Isom}_R(E,F)$ is the closed subspace of $\uHom_R(E,F) \times \uHom_R(F,E)$ obtained as the preimage of the identity section under the map of separated $R$-schemes $\uHom_R(E,F) \times_R \uHom_R(F,E) \to \uHom_{R}(E,E) \times \uHom_R(F,F)$. Hence $\op{Isom}_R(E,F)$ is affine, i.e., $\cM_\cA$ has affine diagonal.
\end{proof}

\begin{rem}
If $\cM_\cA$ is an algebraic stack locally of finite presentation over $k$, then $\cA$ must be Hom-finite. For instance, if $\kappa$ is a field over $k$, then in the proof above we show that for any $E,F \in \cM_\cA(\kappa)$ the algebraic space $\underline{\Hom}_\kappa(E,F)$ is of finite presentation over $\Spec(\kappa)$. In particular this implies that the $\kappa$-module $\Hom_\kappa(E,F)$ is finite dimensional.
\end{rem}

\subsection{Construction of good moduli spaces}

We now apply the previous discussion to construct moduli good moduli spaces for objects in a $\base$-linear abelian category $\cA$. As our results require linear reductivity we will now need to assume that ${\base}$ is a noetherian commutative ring over a field of characteristic $0$.

Furthermore we assume that $\cM_{\cA}$, the moduli functor of flat families of finitely presentable objects of \Cref{D:abelian_moduli}, is an algebraic stack locally of finite type and with affine diagonal over $\base$.

\begin{ex} \label{EX:bridgeland_1}
Let $X$ be a projective scheme over an algebraically closed field $k$, and consider the heart of a $t$-structure $\cC \subset \D^b(X)$ that is noetherian, bounded with respect to the usual $t$-structure, and satisfies the ``generic flatness property'' \cite[Prop.~3.5.1]{abramovich2006sheaves}. Then if we consider the ind-completion $\cA := \Ind(\cC)$, it is shown in \cite[Prop.~6.2.7]{hlinstability} that $\cM_\cA$ is an open sub-functor of the moduli functor $\cD^b_{pug}(X)$ of universally gluable relatively perfect complexes on $X$ and is hence an algebraic stack locally of finite type with affine diagonal over $k$ \cite{lieblich}, \cite[\spref{0DPW}]{stacks-project}. Furthermore, \cite[Prop.~6.2.7]{hlinstability} shows that $\cM_\cA$ agrees with the most commonly studied moduli functor for flat families of objects in the heart of a $t$-structure, introduced in \cite[Defn.~3.3.1]{abramovich2006sheaves} (see also \cite[Prop.~3.3.7]{polishchuk2007constant}).
\end{ex}

The first result concerns a situation in which no additional stability condition is required.

\begin{thm} \label{T:abelian_moduli_space_1}
Let $k$ be an excellent ring of characteristic $0$ and $\cA$ be a cocomplete $\base$-linear abelian category that is locally noetherian. Assume that $\cM_\cA$ is an algebraic stack locally of finite type over $\base$. Then any quasi-compact closed substack $\cX \subset \cM_{\cA}$ admits a proper good moduli space, and in this case points of $\cX$ must parameterize objects of $\cA$ of finite length.
\end{thm}
\begin{proof}
The stack $\cX$ is $\Theta$-reductive (\Cref{L:abelian_Theta_reductive}) and  {\textsf{S}}-complete (\Cref{L:abelian_S_complete}) with respect to essentially finite type DVRs because both properties pass to closed substacks, and it has affine diagonal by \Cref{L:abelian_diagonal}. \Cref{P:SCompleteReductive} implies that closed points of $\cX$ have linearly reductive stabilizers, and thus $\cX$ is locally reductive. \Cref{P:theta-reductive-essentially-finite-type-DVRs} and \Cref{P:S-complete-essentially-finite-type-DVRs} then imply that $\cX$ is $\Theta$-reductive and {\textsf{S}}-complete, so \Cref{T:A} implies the existence of a separated good moduli space $\cX \to X$. Since $\cX$ satisfies the existence part of the valuative criterion for properness with respect to essentially finite type DVRs (\Cref{L:abelian_S_complete}), $X$ is proper by \Cref{P:gms-S-complete}. The fact that closed points of $\cX$ are represented by semisimple objects in $\cA_\kappa$ for fields $\kappa$ of finite type over ${\base}$ is \Cref{L:abelian_semi_simple}.
\end{proof}

In general we will need a notion of ``semistable'' objects in $\cA$. As in applications different notions of stability  are used, we will use an abstract setup that includes many of these. We will illustrate how classical notions of stability fit into this context below.

{We denote by $\pi_0(\cM_\cA)$ the set of connected components of the stack $\cM_\cA$}. For any  $\nu \in \pi_0(\cM_\cA)$, we let $\cM_{\cA}^\nu \subset \cM$ be the corresponding open and closed substack. Our notion of semistability on $\cM_\cA^\nu$ will be encoded by a locally constant function on $|\cM_\cA|$
\[
\p_\nu \co {|\cM_\cA|\to} \pi_0(\cM_\cA) \to V
\]
where $V$ is a totally ordered abelian group, {such that} $\p_\nu(E) = 0$ for any $E \in \cM_\cA^\nu$, and $\p_{\nu}$ is additive in the sense that $\p_\nu(E \oplus F) = \p_{\nu}(E) + \p_\nu(F)$. We will say that a point of $\cM_\cA^\nu$ represented by $E \in \cA_\kappa$ for some algebraically closed field $\kappa$ over ${\base}$, is \emph{$p_\nu$-semistable} if for any subobject $F \subset E$, $\p_\nu(F) \leq 0$ and $p_\nu$-unstable otherwise.\footnote{Note that because the flag space $\uMap(\Theta, \cM_{\cA}^\nu) \times_{\ev_1, \cM_{\cA}^\nu, [E]} \Spec(\kappa)$ of $[E] \co \Spec(\kappa) \to \cM_{\cA}^\nu$  is an algebraic space locally of finite type over $\kappa$, if there is a destabilizing subobject of $E$ after base change to an arbitrary field extension $\kappa \subset \kappa'$, then there is a destabilizing subobject for $E$ over $\kappa$, so this definition does not depend on the choice of representative.} Note that this definition is unaffected by embedding $V$ in a larger totally ordered group, so we may assume that $V$ is a totally ordered vector space over $\bR$ by the Hahn embedding theorem.

Using \Cref{C:reese-ThetaR} to identify maps $f \co \Theta_\kappa \to \cM_\cA$ with $\bZ$-weighted descending filtrations $\cdots \subset E_{w+1} \subset E_w \subset\cdots$ in $\cA_{\kappa}$, we define a locally constant function $\ell \co |\uMap_{\base}(\Theta,\cM_\cA^\nu)| \to V$ by the formula
\[
\ell(\cdots \subset E_{w+1} \subset E_w \subset \cdots) := \sum_w w \p_\nu(E_w / E_{w+1}).
\]
\begin{lem}
	A point $x \in |\cM_\cA^\nu|$ is $p_\nu$-unstable if and only if there is some $f \in |\uMap_{\base}(\Theta,\cM_{\cA}^\nu)|$ with $f(1) = x$ and $\ell(f)>0$.
\end{lem}
\begin{proof}
	If $F \subset E$ is a destabilizing subobject, then we can simply consider the filtration where $E_2=0$, $E_1=F$, and $E_w = E$ for $w \leq 0$. This filtration has $\ell(E_\bullet) = {p_\nu(F)}>0$. 
	
	The converse is a linear algebra statement: Given a $\bZ$-weighted filtration such that $\ell(\cdots \subset E_{w+1} \subset E_w \subset \cdots) := \sum_w w \p_\nu(E_w / E_{w+1})>0$ and $\p_v(E) = \sum_w \p_\nu(E_w / E_{w+1})=0$  it follows that for some index $i$ we have $$\p_v(E_i) = \sum_{w\geq i} \p_\nu(E_w / E_{w+1})>0$$
	so one of the filtration steps will be destabilizing.
\end{proof}

\begin{ex}[Gieseker stability]
Let $\cA = \QCoh(X)$ denote the category of quasi-coherent sheaves on a scheme $X$ which is projective over a field $k$ of characteristic 0. Then $\cM_\cA$ is the usual stack of flat families of coherent sheaves on $X$. Fix a numerical $K$-theory class $\gamma \in K_0^{\num}(X)$ corresponding to sheaves whose support has dimension $d$, fix an ample line bundle $\cO_X(1)$ on $X$, and let $\cM_\cA^\gamma$ denote the open and closed substack of objects of class $\gamma$. For any $E \in \Coh(X)$, let $P(E) = \alpha_n(E) t^n + \cdots + \alpha_0(E)$ denote its Hilbert polynomial with respect to $\cO_X(1)$.

We say that a coherent sheaf $F$ of class $\gamma$ is semistable with respect to $\cO_X(1)$ if for all nonzero subsheaves $E \subset F$, the asymptotic inequality
\[
p_\gamma(E) := \alpha_d(F) P(E) - \alpha_d(E) P(F) \le 0
\]
holds for all $t \gg 0$. If $E \subset F$ is a nonzero subsheaf with $\dim(\op{supp}(E)) < d$, then $\alpha_d(E) =0$ but $\alpha_d(F) >0$, so $E$ destabilizes $F$. Therefore a semistable sheaf must be pure, and in this case $P(E) / \alpha_d(E)$ and $P(F)/\alpha_d(F)$ are the reduced Hilbert polynomials, so our definition is equivalent to the classical definition of Gieseker semistability \cite[Def.~1.2.4]{huybrechtslehn}. Note also that $P(F) = P(\gamma)$ only depends on $\gamma$, so $p_\gamma$ defines a function $p_\gamma \co \cM_\cA \to V_d$, where $V_d$ denotes the vector space of polynomials of degree $\leq d$ totally ordered by asymptotic inequality as $t \to \infty$. The function $p_\gamma$ is locally constant in flat families and additive in short exact sequences.
\end{ex}

\begin{ex}[Bridgeland stability]
	Consider a projective scheme $X$ over an algebraically closed field ${\base}$ of characteristic $0$. A Bridgeland stability condition is determined by the heart of a $t$-structure $\cC \subset \D^b(X)$, and a \emph{central charge} homomorphism $Z \co K_0(\D^b(X)) \to \bC$, which we assume factors through the numerical $K$-theory $K_0^{\num}(\D^b(X))$, i.e., the quotient $K_0(\D^b(X))$ by the kernel of the Euler pairing $\chi(-,-)$. The central charge $Z$ is required to be a stability function on $\cC$ with the Harder--Narasimhan property \cite[Prop. 5.3]{bridgeland}.  The simplest example is  when $X$ is a projective curve, $\cC = \Coh(X)$ is the usual $t$-structure, and
	\[
	Z(E) = - \deg(E) + i \rank(E).
	\]
	For a stability condition $(\cC, Z)$ on a general $X$, we take the formula above as the definition of rank and degree of objects in $\cC$. As in \Cref{EX:bridgeland_1} we let $\cA := \Ind(\cC)$. Given a class $\gamma \in K_0^{\num}(\D^b(X))$, there is an open and closed substack $\cM_\cA^\gamma \subset \cM_\cA$ whose $k$-points classify objects $E \in \cA$ of numerical class $\gamma$.
Then Bridgeland semistability on the stack $\cM_\cA^\gamma$ is determined by the degree function
	\[
	\p_\gamma(E) := \deg(E) \rank(\gamma) - \deg(\gamma) \rank(E),
	\]
	so $E \in \cM_\cA^\gamma$ is unstable if and only if there is a subobject $F \subset E$ such that $\p_\gamma(F)>\p_\gamma(E) = 0$.
\end{ex}

\begin{thm} \label{T:abelian_semistable_GMS}
Let $k$ be an excellent ring of characteristic $0$ and $\cA$ be a locally noetherian, cocomplete and $\base$-linear abelian category. Assume that $\cM_\cA$ is an algebraic stack locally of finite type over $\base$.  Let $\nu \in \pi_0(\cM_\cA)$ be a connected component, and let $\p_\nu \co \pi_0(\cM_\cA) \to V$ be an additive function defining a notion of $p_\nu$-semistability on $\cM_\cA^\nu$, as above. 
	
	If the substack of $p_\nu$-semistable points $\cM_\cA^{\nu,\ss} \subset \cM_\cA^\nu$ is open and quasi-compact then $\cM_\cA^{\nu,\ss}$ admits a separated good moduli space. If in addition $\cM_\cA^{\nu,\ss}$ is the open piece of a $\Theta$-stratification of $\cM_\cA^\nu$, then $\cM_{\cA}^{\nu,\ss}$ admits a proper good moduli space.
\end{thm}

\begin{proof}
	We have already seen that $\cM_\cA^\nu$ has affine diagonal (\Cref{L:abelian_diagonal}), and with respect to essentially finite type DVRs $\cM_\cA^\nu$ is $\Theta$-reductive (\Cref{L:abelian_Theta_reductive}),  {\textsf{S}}-complete (\Cref{L:abelian_S_complete}), and satisfies the existence part of the valuative criterion for properness (\Cref{L:abelian_properness}).  It follows from \Cref{P:semistable_locus} and \Cref{R:more_general_stability} that $\cM_\cA^{\nu,\ss}$ is $\Theta$-reductive and  {\textsf{S}}-complete with respect to essentially finite type DVRs as well, so \Cref{P:SCompleteReductive} implies $\cM_\cA^{\nu,\ss}$ is locally reductive, and it is $\Theta$-reductive and \textsf{S}-complete by \Cref{P:theta-reductive-essentially-finite-type-DVRs} and \Cref{P:S-complete-essentially-finite-type-DVRs} respectively. Therefore, \Cref{T:A} implies that there is a separated good moduli space $\cM_\cA^{\nu,\ss} \to M$. For the final statement,  \Cref{L:abelian_properness} and \Cref{C:valuative_criterion_semistable} applied to the $\Theta$-stratification of $\cM_\cA^\nu$ imply that $\cM_\cA^{\nu,\ss}$ satisfies the existence part of the valuative criterion for properness with respect to essentially finite type DVRs and hence $M$ is proper over $\Spec(k)$. 
\end{proof}

\begin{ex}[Gieseker stability, continued] 
The Harder--Narasimhan stratification with respect to Gieseker semistability defines a $\Theta$-stratification of the stack $\uCoh(X)$ (this is essentially the content of \cite{nitsure2011schematic}, using different language), and the semistable locus is open and bounded.  Therefore \Cref{T:abelian_semistable_GMS} provides an alternate construction of a proper good moduli space for Gieseker semistable sheaves.
\end{ex}

\begin{ex} [Bridgeland stability, continued]
We consider Bridgeland stability conditions whose central charge factors through a fixed surjective homomorphism $\op{cl} \co K_0(\D^b(X)) \to K_0^{\num}(\D^b(X)) \to \Gamma$, where $\Gamma$ is a finitely generated free abelian group, and we let $\Stab_\Gamma(X)$ denote the space of Bridgeland stability conditions whose central charge factors through $\op{cl}$ and which satisfy the support property with respect to $\Gamma$. This is given the structure of a complex manifold in such a way that the map $\Stab_\Gamma(X) \to \Hom(\Gamma, \bC)$ which forgets the $t$-structure is a local isomorphism \cite[Thm. 1.2]{bridgeland}.


A stability condition $(\cC,Z)$ is algebraic if $Z(\Gamma) \subset \bQ + i \bQ$. If $(\cC,Z)$ is algebraic, then $\cC$ is noetherian \cite[Prop.~5.0.2]{abramovich2006sheaves}. Let $\Stab_\Gamma^\ast(X) \subset \Stab_\Gamma(X)$ be a connected component which contains an algebraic stability condition $\sigma_0 = (\cC_0, Z_0)$ for which:
\begin{itemize}
\item the heart $\cC_0$ satisfies the generic flatness condition, 
\item $\cC_0$ is bounded with respect to the usual $t$-structure, and
\item $\cM_{\Ind(\cC_0)}^{\gamma,\ss}$ is bounded for every $\gamma \in \Gamma$.
\end{itemize}
Then by \cite[Prop.~4.12]{toda2015moduli} the same is true for \emph{any} algebraic stability condition in the connected component $\Stab_\Gamma^\ast(X)$. Furthermore, if $\sigma = (\cC,Z) \in \Stab_\Gamma(X)$ is algebraic and satisfies the generic flatness and boundedness conditions, then the Harder--Narasimhan filtration defines a $\Theta$-stratification of $\cM_{\Ind(\cC)}^\gamma$ whose open piece is $\cM_{\Ind(\cC)}^{\gamma,\ss}$ \cite[Thm.~6.5.3, Rem.~6.5.6]{hlinstability}. Therefore \Cref{T:abelian_semistable_GMS} implies that $\cM_{\Ind(\cC)}^{\gamma,\ss}$ has a proper good moduli space for any algebraic stability condition in $\Stab^\ast(X)$.

Finally, we claim that for an \emph{arbitrary} stability condition $\sigma \in \Stab^\ast(X)$, one can find an algebraic stability condition $\sigma'$ which defines the same moduli functor $\cM_{\Ind(\cC)}^{\gamma,\ss}$, so $\cM_{\Ind(\cC)}^{\gamma,\ss}$ has a proper good moduli space for any stability condition in $\Stab^\ast(X)$.

To establish this claim, fix a class $\gamma \in \Gamma$. For any $\gamma' \in \Gamma$ which is linearly independent of $\gamma$ over $\bQ$, consider the real codimension $1$ subset
\[
\cW_{\gamma'} := \{\sigma = (\cC,Z) \in \Stab_\Gamma^\ast(X) \, | \, Z(\gamma') \in \bR_{>0} \cdot Z(\gamma) \}.
\]
If one restricts to a small compact neighborhood $\fB\subset \Stab_\Gamma^\ast(X)$ containing $\sigma$, then there is a finite subset $S \subset \Gamma$ such that for any $S' \subset S$ the moduli functor $\cM_{\Ind(\cC)}^{\gamma,\ss}$ is constant for all $\sigma \in \fC_{S'} \cap \fB$ \cite[Prop.~2.8]{toda}, where
\[
\fC_{S'} := \bigg(\bigcup_{\gamma' \in S'} \cW_{\gamma'} \bigg) \smallsetminus \bigcup_{\gamma' \notin S'} \cW_{\gamma'}.
\]
$\cM_{\Ind(\cC)}^{\gamma,\ss}$ is constant for $\sigma \in \fB \cap \fC_{S'}$ because $\cM_{\Ind(\cC)}^{\gamma,\ss}$ can only change if the set of classes $\gamma' \in \Gamma$ with $Z(\gamma') \in \bR \cdot Z(\gamma)$ changes, and the set of such $\gamma'$ is constant for $\sigma \in \fC_{S'} \cap \fB$.

The condition that $\sigma \in \bigcup_{\gamma' \in S'} \cW_{\gamma'}$ amounts to the claim that if $W \subset \Gamma_\bQ$ is the span of $\gamma$ and the $\gamma' \in S'$, then $\dim_\bQ(Z(W)) = 1$. We may write $Z = Z_1 \oplus Z_2$ under a choice of splitting $\Gamma_\bQ \simeq W \oplus U$, and the condition now amounts to $\rank (Z_1) = 1$. As rational points are dense in the space of rank 1 real matrices, we may find an arbitrarily small perturbation $Z'$ of $Z$ which is rational, and this perturbed central charge defines our new algebraic stability condition $\sigma' \in \fB \cap \fC_{S'}$.
\end{ex}


\section{Good moduli spaces for moduli of \texorpdfstring{$\cG$}{G}-torsors}

To illustrate our general theorems we now construct good moduli spaces for semistable torsors under Bruhat--Tits group schemes. This generalizes the results obtained by Balaji and Seshadri who constructed such moduli spaces for generically split groups over the complex numbers. In \cite{heinloth-stability}, the third author analyzed the coarse moduli space for the moduli of stable bundles, whose existence is guaranteed by the Keel--Mori theorem. Here we extend this analysis to include semistable bundles which are not stable. As in this article we are interested in existence theorems for good moduli spaces (instead of adequate moduli) we will have to assume that we work over a base field $k$ of characteristic 0 in this section.

Let us briefly introduce the setup from \cite{heinloth-stability}. Let $C$ be a smooth geometrically connected, projective curve over a field $k$ and $\cG/C$ a smooth Bruhat--Tits group scheme over $C$, i.e., $\cG$ is smooth a affine group scheme over $C$ that has geometrically connected fibers, such that over some dense open subset $U\subset C$ the group scheme is reductive and over all local rings at points $p$ in $\Ram(\cG) := C\smallsetminus U$ the group scheme $\cG|_{\Spec(\cO_{C,p})}$ is a connected parahoric Bruhat--Tits group. The simplest examples are of course reductive groups $G\times C$.

The stack of $\cG$-torsors is denoted by $\Bun_{\cG}$ and this is a smooth algebraic stack. To define stability one usually chooses a line bundle on $\Bun_{\cG}$. As explained in \cite[\S 3.B]{heinloth-stability} there are natural choices in our situation. First there is the determinant line bundle $\cL_{\det}$ given by the adjoint representation, i.e., the fiber at a bundle $\cE\in \Bun_{\cG_p}$ is $\cL_{\det,\cE} = \det(H^*(C,\ad(\cE)))^\vee$, where $\ad(\cE)=\cE\times^{\cG}\Lie(\cG_p/C)$ is the adjoint bundle of $\cE$.

Next any collection of characters $\underline{\chi}\in \prod_{p\in \Ram(\cG)} \Hom(\cG_p,\bG_m)$ defines line bundles on the classifying stacks $B\cG_p$ and one obtains  a line bundle $\cL_{\underline{\chi}}$ on $\Bun_{\cG}$, by pull back via the map $\Bun_{\cG} \to B\cG_p$ defined by restriction of $\cG$ torsors on $C$ to the point $p$. We will denote by $\cL_{\det,\underline{\chi}}:=\cL_{\det}\tensor \cL_{\underline{\chi}}$, call the corresponding notion of stability $\underline{\chi}$-stability and denote by $\Bun_{\cG}^{\underline{\chi}-\ss} \subset \Bun_{\cG}$ the substack of $\underline{\chi}$-semistable torsors.

Under explicit numerical conditions on $\underline{\chi}$ this satisfies the positivity assumption of loc.cit. \cite[Prop.~3.3]{heinloth-stability}, i.e., the restriction of $\cL_{\det}$ to the affine Grassmannian $\Gr_{\cG,p}$ classifying $\cG$ bundles together with a trivialization on $C\smallsetminus p$ is nef. 
The parameter $\underline{\chi}$ will be called positive if $\cL_{\det,\underline{\chi}}$ is ample on $\Gr_{\cG,p}$ for all $p$. It is called admissible if $\underline{\chi}$ furthermore satisfies the numerical condition of \cite[Sec. 3.F]{heinloth-stability}.

\begin{thm}[Good moduli for semistable $\cG$-torsors]\label{T:GoodModuliBunG}
	Assume $k$ is a field of characteristic $0$, $C$ is a smooth, projective, geometrically connected curve over $k$, $\cG$ is a parahoric Bruhat--Tits group scheme over $k$ and $\underline{\chi}$ is a admissible stability parameter. 
	Then $\Bun_{\cG}^{\underline{\chi}-\ss}$ admits a proper good moduli space $M_{\cG}$.
\end{thm}
As remarked before, in the case that $\cG$ is a generically split group scheme, the space $M_{\cG}$ was constructed by Balaji and Seshadri \cite{BalajiSeshadri}.

To prove the theorem we only need to check that $\Bun_{\cG}$ satisfies the assumptions of our main Theorems \ref{T:A}, \ref{T:B} and \ref{T:C}, i.e., we need to show that the line bundle $\cL_{\det,\chi}$ defines a well-ordered $\Theta$-stratification of $\Bun_{\cG}$, that this stack satisfies the existence part of the valuative criterion for properness and that $\Bun_{\cG}^{\underline{\chi}-\ss}$ is $\Theta$-reductive and  {\textsf{S}}-complete. This will be done in a series of Lemmas.

\begin{lem} \label{L:theta-strat-bun-g}
	The canonical reduction of $\cG$-torsors defines a $\Theta$-stratification on $\Bun_{\cG}$ with semistable locus $\Bun_{\cG}^{\underline{\chi}-\ss}$. This stratification admits a well-ordering.
\end{lem}

Let us briefly recall the context of $\underline{\chi}$-stability. To simplify the presentation will assume that our base field is algebraically closed. Then by \cite[Lem. 3.9]{heinloth-stability} any map $f\co\Theta\to \Bun_{\cG}$ arises as a Rees construction $\Rees(\cE_{\cP},\lambda)$ from a $\cG$-bundle $\cE$, together with a reduction $\cE_{\cP}$ to a parabolic subgroup $\cP\subset \cG$ and a generic cocharacter $\lambda\co \bG_{m,k(C)}\to \cG_{k(C)}$ that is dominant for $\cP_{\eta}$. The argument of the proof implies that all components of  $\uMap(\Theta,\Bun_\cG)$ can be identified with components of $\Bun_{\cP}$ for parabolics $\cP$ equipped with a dominant $\lambda$. As reductions to parabolics are determined by the induced filtration of the adjoint bundle, this implies in particular that the forgetful morphism from any component to $\Bun_\cG$ is quasi-compact. We denote by $\wt_{\cE_{\cP}}(\lambda)$ the weight of $f^*(\cL_{\chi})$.

\begin{proof}[Proof of \Cref{L:theta-strat-bun-g}]

In \cite[Thm.~2.2.2, Simplif.~2.2.3, Simplif.~2.2.4]{hlinstability} a list of necessary and sufficient conditions for a stratification to be a $\Theta$-stratification is given. We will first verify these criteria (which we number as in loc.cit.), then briefly recall how the general proof works in the specific context of $\cG$-bundles in order to illustrate the relation to classical arguments of Behrend. It suffices to assume the ground field $k$ is algebraically closed \cite[Lem.~2.3.3]{hlinstability}.

\medskip
\noindent{\textit{Existence and uniqueness of HN filtrations $(1)$:}}
\medskip

$HN$ filtrations are maps $f \co \Theta \to \Bun_\cG$ corresponding to the canonical parabolic reductions for unstable $\cG$-bundles constructed in \cite[Sec.\ 3.F]{heinloth-stability}. To define canonical reductions in this setup we fixed an invariant inner form $(,)$ on the generic cocharacters of $\cG$ defining the norm $||\cdot||$ we showed \cite[Prop.\ 3.6]{heinloth-stability} that for any admissible $\underline{\chi}$ and any $\cG$-torsor $\cE$ there exists a canonical reduction $\cE_{\cP}$ to a parabolic subgroup $\cP\subset \cG$ defined by a dominant cocharacter $\lambda\co \bG_{m,\eta} \to \cG_{\eta}$ maximizing $\mu_{\max}(\cE):= \max(\frac{\wt_{\cE_\cP}(\lambda)}{||\lambda||})$.

The uniqueness of the reduction followed by checking that for any choice of a generic maximal torus $\cT_{\eta}\subset \cG_{\eta}^{\cE} \cong \cG_{\eta}$ that contains a maximal split torus the map sending a Borel subgroup $\cB_{\eta}\subset \cG_{\eta}$ containing $\cT_{\eta}$ to the cocharacter defined by $-\wt_{\cE_{\cB}}(\cdot)$ defines a complementary polyhedron in the sense of Behrend \cite[Def.\ 2.1]{behrend-stability}.
	
\medskip
\noindent{\textit{Consistency of HN filtrations $(5)$:}}
\medskip

Any canonical reduction $f\co\Theta\to \Bun_{\cG}$ of $f(1)$ also defines the canonical reduction of the associated graded bundle $f(0)$. Indeed, the perspective of complementary polyhedra shows that a canonical reduction $\cE_{\cP}$ can also be characterized by the property that (1) $\wt_{\cE_{\cP}}(\lambda^\prime)>0$ for any non-zero $\lambda^\prime$ that is non-negative on any root of $\cP$, and (2) such that for any reductions $\cE_{\cB}$ of $\cE_{\cP}$ to a Borel subgroups $\cB\subset \cP$ we have $\wt_{\cE_{\cB}}(\check{\alpha}) < 0$ for all simple roots coroots $\check{\alpha}$ in the centralizer of $\lambda$ \cite[Sec. 3]{behrend-stability} \cite[Thm 4.3.2]{heinloth-schmitt}.

\medskip
\noindent{\textit{Specialization of HN filtrations $(2')$:}}
\medskip

For any family $\cE_R \in \Bun_{\cG}(R)$ defined over a DVR $R$ with fraction field $K$ and residue field $\kappa$ we have $\mu_{\max}(\cE_K)\leq \mu_{\max}(\cE_\kappa)$ and equality holds only if the canonical filtration over $K$ extends to the family, by \cite[Lem.~3.17]{heinloth-stability}.

\medskip
\noindent{\textit{Local finiteness $(4)$:}}
\medskip

For any family of $\cG$-bundles over a finite type scheme $U$, only finitely many connected components of $\uMap(\Theta,\Bun_\cG)$ are necessary to realize the HN filtration of every fiber $\cE_u$ for $u \in U$. This follows from the fact that there are only finitely many $\lambda$ such the canonical reduction of $\cE_u$ for some $u \in U$ has type $\cP_\lambda$, which is established in the proof of \cite[Prop.~3.18]{heinloth-stability}.

\medskip
\noindent{\textit{Completing the proof:}}
\medskip

The function $\frac{\wt_{\cE_\cP}(\lambda)}{||\lambda||}$ defines a locally constant real valued function $\mu$ on the components of $\uMap(\Theta,\Bun_\cG)$ containing a HN filtration. We have verified that $\mu$ satisfies the conditions of \cite[Thm.~2.2.2, Simplif.~2.2.3, Simplif.~2.2.4]{hlinstability} and thus defines a weak $\Theta$-stratification on $\Bun_\cG$ in which the HN filtrations of unstable points correspond to canonical reductions of $\cG$-bundles. In our context, the argument goes as follows:

``Local finiteness'' and quasi-compactness of $\Bun_{\cP_\lambda} \to \Bun_\cG$ implies that the stratification of $\Bun_{\cG}$ defined by $\mu_{\max}$ is constructible (See also the proof of \cite[Prop.~3.18]{heinloth-stability}). Since the invariant $\mu_{\max}$ is semicontinuous, this implies that for any constant $c$ the substacks $\Bun_{\cG}^{\mu_{\max}\leq c}$ defined by the condition $\mu_{\max}(\cE)\leq c$ are open. To show that this defines a $\Theta$-stratification we are therefore left to show that the closed substacks $\Bun_{\cG}^{\mu_{\max}\leq c}\smallsetminus \Bun_{\cG}^{\mu_{\max}<c}$ are unions of connected components of $\uMap(\Theta,\Bun_{\cG}^{\mu_{\max}\leq c})$.
	
Let us fix an unstable bundle $\cE$ with $\mu_{\max}(\cE)=c$ and canonical reduction given as a reduction to $\cP_\lambda$. Let us denote by $p_\lambda\co \Bun_{\cP_\lambda}^{\HN} \to \Bun_{\cG}^{\mu_{\max}\leq c}$ the restriction of the canonical map $\Bun_{\cP_\lambda} \to \Bun_{\cG}$. The ``consistency of HN filtrations,'' i.e., the fact that passing to the associated graded of the canonical filtration preserves the strata, shows that $\Bun_{\cP_\lambda}^{\HN}$ is indeed a component of $\uMap(\Theta,\Bun_{\cG}^{\mu_{\max}\leq c})$, and by uniqueness of the filtration and the ``specialization property,'' the map $p_\lambda$ is proper and universally injective. Hence the stratification induced by $\mu_{\max}$ is a weak $\Theta$-stratification.

Any weak $\Theta$-stratification is a $\Theta$-stratification in characteristic $0$ \cite[Cor.~2.6.1]{hlinstability}. Alternatively, it is not hard to check directly that $p_\lambda$ is injective on tangent spaces at any HN filtration, as in the case of Behrend's conjecture, so $p_\lambda$ is a closed immersion whose image is $\Bun_{\cG}^{\mu_{\max}\leq c}\smallsetminus \Bun_{\cG}^{\mu_{\max}<c}$. Finally, this $\Theta$-stratification that admits a well-ordering, because for any $c$ and any connected component of $\Bun_{\cG}$ the open substack $\Bun_{\cG}^{\mu_{\max}\leq c}$ are of finite type \cite[Prop.~3.18]{heinloth-stability}, so $\Bun_{\cG}^{\mu_{\max}\leq c}\smallsetminus \Bun_{\cG}^{\mu_{\max}< c}$ can only contribute finitely many strata on each component of $\Bun_{\cG}$.
\end{proof}

\begin{rem}
Assume for simplicity that the ground field $k$ is algebraically closed. The notion of $\underline{\chi}$-stability is controlled by the class $\ell = c_1(\cL_{\underline{\chi}}) \in H^2(\Bun_\cG;\bQ)$. If $p \in C$ is a regular point for $\cG$, then generic cocharacters of $\cG$ induce cocharacters in $\cG_p$, and under this map the norm $||\lambda||$ on generic cocharacters can be induced from a conjugation invariant norm on cocharacters of $\cG_p$. It follows that $||\lambda||$ is induced by a class in $H^4(\Bun_\cG;\bQ)$ defined via pullback
\[
(\op{Sym}^2(N^\ast_\bQ))^W \simeq H^4(B\cG_p;\bQ) \to H^4(\Bun_\cG;\bQ),
\]
along the restriction morphism $\Bun_\cG \to B \cG_p$, where $N$ denotes the coweight lattice of $\cG_p$ and $W$ denotes the Weyl group.

Therefore the function $\mu(\cE_\cP,\lambda) = \frac{\wt_{\cE_\cP}(\lambda)}{||\lambda||}$ is a standard numerical invariant in the sense of \cite[Def.~4.1.10]{hlinstability} satisfying condition (R) by \cite[Lem.~4.1.15, Lem.~4.1.16]{hlinstability}. It is also not difficult to check this directly. 
So in the proof of \Cref{L:theta-strat-bun-g} one could also apply \cite[Thm.~4.5.1]{hlinstability}, which provides a shorter list of criteria for a numerical invariant to define a weak $\Theta$-stratification. In particular this implies that condition (5) above is automatic here.
\end{rem}

\begin{lem}\label{L:BunGScomplete}
	The stack $\Bun_{\cG}^{\underline{\chi}-\ss}$ is  {\textsf{S}}-complete and every closed point has linearly reductive stabilizer. 
\end{lem}
\begin{proof}
	 {\textsf{S}}-completeness holds because of the existence of a blow up of $\oST_R$ to linking two specializations. $\cL_{\det,\underline{\chi}}$ is positive on the exceptional lines, so  if the blow-up was necessary, one of the bundles was unstable.
	
	In particular every closed substack of $\Bun_{\cG}^{\Bun_{\cG}^{\underline{\chi}-\ss}}$ is again  {\textsf{S}}-complete, so that by \Cref{P:SCompleteReductive} the automorphism groups of closed points are geometrically reductive. As we assumed our base field to be of characteristic $0$ in this section, these groups are linearly reductive.  
\end{proof}

\begin{lem}
The stack $\Bun_{\cG}^{\underline{\chi}-\ss}$ is $\Theta$-reductive.
\end{lem}
\begin{proof}
	This again follows from semi-continuity of the numerical invariant as in \cite[Lem. 3.17]{heinloth-stability}. We briefly recall the argument:
	Let $R$ be a DVR with fraction field $K$ and residue field $\kappa$. Let $f\co \Theta_R\smallsetminus 0 \to \Bun_{\cG}^{\underline{\chi}-\ss}$ be a morphism. The restriction of $f$ to $\Spec(R)$ defines a family $\cE_R$ of $\cG$-torsors over $C_R$.
	
	The restriction of $f|_{\Theta_K}$ defines a filtrations on $\cE_K$ and  
	after possibly passing to a finite extension of $R$ we may by \cite[Rem. 3.10]{heinloth-stability} assume that this filtration is given by a reduction $\cE^{\cP}_K$ of $\cE_K$ to a parabolic subgroup $\cP\subset \cG$ which is defined by a generic cocharacter $\lambda\co \bG_m \to \cG_{k(\eta)}$.
		
	As $\cG_{k(\eta)}/\cP_{k(\eta)}$ is projective, any reduction of $\cE_{K(\eta)}$ to $\cP$ extends to an open subset $U\subset C_R$ whose the complement consists of finitely many closed points of the special fiber. Also the reduction over $U_\kappa$ extends canonically to a reduction $\cE^{\cP^\prime}_{\kappa}$ to a parabolic subgroup $\cP^\prime\subset \cG_\kappa$ over $C_\kappa$. 
	Now, as in loc.\ cit.\  at any point $p$ the adjoint bundle $\ad(\cE^{\cP^\prime}_{\kappa}) \subset \ad(\cE_\kappa)$ is the saturation of the adjoint bundle at the generic fiber. As $\underline{\chi}$ is admissible this implies that either the reduction $\cE^{\cP}_K$ extends over $C_R$ or the weight of the reduction $\cE^{\cP^\prime}_{\kappa}$ is strictly larger than the weight of $f|_{\Theta_K}$, which is $0$ as $f$ is a reduction in $\Bun_{\cG}^{\underline{\chi}-\ss}$. As by assumption $\cE_\kappa$ is $\underline{\chi}$-semistable the weight cannot increase, so the reduction extends to $C_R$. Now we can apply \Cref{L:semistable_filtrations} to find that this filtration also lies in $\Bun_{\cG}^{\underline{\chi},ss}$ and thus defines an extension of $f$ to $\Theta_R$.
\end{proof}

\begin{proof}[Proof of \Cref{T:GoodModuliBunG}] 
	We just proved that $\Bun_{\cG}^{\chi-\ss}$ is  {\textsf{S}}-complete, $\Theta$-reductive, and that every closed point has a linearly reductive stabilizer.

	By \cite[Prop.~3.3]{heinloth-stability} the stack $\Bun_{\cG}$ satisfies the existence criterion for properness, i.e., if $R$ is a DVR with fraction field $K$ and $\cE_K\in \Bun_{\cG}(K)$ is a $\cG$-torsor over $C_K$ then there exists a finite extension $R^\prime$ of $R$ such that $\cE_K$ extends to a torsor over $C_R$.  Therefore we can apply \Cref{T:C} to deduce the existence of a proper good moduli space.
\end{proof}


\appendix

\section{Strange gluing lemma}

In this section, we establish two gluing results:  \Cref{T:strange-gluing} and \Cref{T:strange-gluing2}.  Additionally, we apply both results to give a refinement of the classical semistable reduction theorem in GIT (\Cref{T:valuative-criterion-finite}).

\subsection{Gluing results}

Let $R$ be a DVR with fraction field $K$ and residue field $\kappa$, and let $\pi \in R$ be a uniformizer parameter. For $n>0$ we will consider the following quotient stack
\[
\oST_R^{n,1} = [\Spec(R[s,t]/(s t^n - \pi)) / \bG_m]
\]
where the $\bG_m$-action is encoded by giving $s$ weight $n>0$ and giving $t$ weight $-1$. We have a closed immersion $\Theta_{\kappa} \hookrightarrow \oST^{n,1}_R$ defined by $s = 0$ and an open immersion $\Spec(R)  \hookrightarrow \oST_R^{n,1}$ defined by $t \neq 0$. 

We will denote $0 \in \Spec(R)$ as the closed point and $1 \in \Theta_k$ as the open point.  Observe that any morphism  $\oST^{n,1}_R \to \cX$ restricts to morphisms $f \co \Theta_\kappa \to \cX$ and $\xi \co \Spec(R) \to \cX$ along with an isomorphism $\phi \co \xi(0) \simeq f(1)$ in $\cX(\kappa)$. 

\begin{thm} \label{T:strange-gluing}
Let $\cX$ be an algebraic stack with affine diagonal and locally of finite presentation over an algebraic space $S$. Let $R$ be a DVR with residue field $\kappa$ and consider morphisms $f \co \Theta_\kappa \to \cX$ and $\xi \co \Spec(R) \to \cX$ over $S$ together with an isomorphism $\phi \co \xi(0) \simeq f(1)$. For all $n\gg 0$, there is a morphism $\oST^{n,1}_R \to \cX$ unique up to unique isomorphism extending the triple $(f, \xi, \phi)$. 
\end{thm}

This theorem is inspired by the perturbation theorem \cite[Thm.~3.5.3]{hlinstability}, which is an analogous result for constructing map $[\bA^2_\kappa/\bG_{m,\kappa}^2] \to \cX$ from maps from the loci $\{s=0\}$ and $\{t\neq 0\}$. 

\begin{cor} \label{C:strange_gluing}
In the context of \Cref{T:strange-gluing}, for $n\gg 0$ the data of the morphisms $f$ and $\xi|_{\Spec(R[\pi^{1/n}])}$ with isomorphism $\phi$ extends canonically to a morphism $\oST^{1,1}_{R[\pi^{1/n}]} \to \cX$.
\end{cor}
\begin{proof}
Compose the uniquely defined map $\oST^{n,1}_R \to \cX$ of \Cref{T:strange-gluing} with the canonical map $\oST^{1,1}_{R[\pi^{1/n}]} \to \cX$ induced by the map of graded algebras $R[s,t]/(st^n-\pi) \to R[\pi^{1/n}][s^{1/n},t]/(s^{1/n}t-\pi)$, where $s^{1/n}$ has weight $1$.
\end{proof}

In order to establish \Cref{T:strange-gluing}, we will need to recall the following fact concerning pushouts.

\begin{lem} \label{L:pushout1}
If $\Spec(A) \to \Spec(B)$ is a closed immersion and $\Spec(A) \to \Spec(C)$ is a morphism, then
$$\xymatrix{
\Spec(A) \ar[d] \ar[r]		& \Spec(C) \ar[d] \\
\Spec(B) \ar[r]			& \Spec(B \times_{A} C)
}$$
is a pushout diagram in the category of algebraic stacks with affine diagonal.
\end{lem}

\begin{proof}
Ferrand established that the diagram is a pushout in the category of ringed spaces \cite[Thm. 5.1]{ferrand}. Temkin and Tyomkin establish that it is a pushout in the category of algebraic spaces \cite[Thm. 4.5]{temkin-tyomkin}.  The lemma follows by applying \cite[Thm 4.7]{temkin-tyomkin} or \cite[Lem. A.4]{hall-coherent} to the pullback of the diagram under an affine presentation of an algebraic stack $\cX$ with affine diagonal.
\end{proof}

\begin{lem} \label{L:pushout2}
Let $C = R[t,\pi/t,\pi/t^2,\ldots] \subset R[t]_t$. The commutative diagram
\begin{equation} \label{D:pushout}
\begin{split}
\xymatrix{
\Spec(\kappa) \ar[r] \ar[d]		& \Theta_{\kappa} \ar[d] \\
\Spec(R) \ar[r]			& [\Spec(C)/\bG_m]
}
\end{split}
\end{equation}
is cartesian and a pushout diagram in the category of algebraic stacks with affine diagonal.
\end{lem}

\begin{proof}
\Cref{L:pushout1} implies that both diagrams
\begin{equation*}
\begin{split}
\xymatrix{
\Spec(\kappa[t]_t) \ar[r] \ar[d]		& \Spec(\kappa[t]) \ar[d] \\
\Spec(R[t]_t) \ar[r]			& \Spec(C) 
}
\qquad
\xymatrix{
\bG_m \times \Spec(\kappa[t]_t) \ar[r] \ar[d]		& \bG_m \times \Spec(\kappa[t]) \ar[d] \\
\bG_m \times \Spec(R[t]_t) \ar[r]			& \bG_m \times \Spec(C) 
}
\end{split}
\end{equation*}
are pushout diagrams in the category of algebraic stacks with affine diagonal.  Since \eqref{D:pushout} is the $\bG_m$-quotient of the left diagram above, the statement follows from descent.
\end{proof}

\begin{proof}[Proof of \Cref{T:strange-gluing}]
By \Cref{L:pushout2}, the triple $(f, \xi, \phi)$ glues to a morphism $[\Spec(C) / \bG_m] \to \cX$ unique up to unique isomorphism. Write $C$ as a union $C = \bigcup C_n$, where $C_n := R[t,\pi/t^n] \subset R[t]_t$. Note that  $C_n \cong R[s,t]/(st^n - \pi)$ so in particular $[\Spec(C_n)/\bG_m] \cong \oST_R^{n,1}$ and that $\Spec(C) \to \Spec(C_n)$ is $\bG_m$-equivariant.  As $\cX \to S$ is locally of finite presentation, for $n \gg 0$ the morphism $[\Spec(C) / \bG_m] \to \cX$ factors uniquely as $[\Spec(C) / \bG_m] \to [\Spec(C_n) / \bG_m]  \to \cX$.
\end{proof}

To setup the second gluing result, for $n > 0$ consider the subalgebra $R[t/\pi^n] \subset K[t]$ and the quotient stack 
\[
\Theta_{R,n} = [\Spec(R[t/\pi^n]) / \bG_m]
\]
where $t$ has weight $-1$. We have a closed immersion $B_R \bG_m \hookrightarrow \Theta_{R,n}$ defined by $t/\pi^n=0$ and an open immersion $\Theta_K \hookrightarrow \Theta_{R,n}$ defined by $\pi \neq 0$.  Observe that any morphism  $\Theta_{R,n} \to \cX$ restricts to morphisms $g \co B_R \bG_m \to \cX$ and $\lambda \co \Theta_K \to \cX$ along with an isomorphism $\phi \co g|_{B_K \bG_m} \simeq \lambda|_{B_K \bG_m}$.  

\begin{thm} \label{T:strange-gluing2}
Let $\cX$ be an algebraic stack with affine diagonal and locally of finite presentation over an algebraic space $S$. Let $R$ be a DVR with fraction field $K$ and consider morphisms $g \co B_R \bG_m \to \cX$ and $\lambda \co \Theta_K \to \cX$ over $S$ together with isomorphism $\phi \co g|_{B_K \bG_m} \simeq \lambda|_{B_K \bG_m}$. For all $n\gg 0$, there is a morphism $\Theta_{R,n} \to \cX$ unique up to unique isomorphism extending the triple $(g, \lambda, \phi)$. 
\end{thm}

\begin{rem} Observe that $\Theta_{R,n} \cong \Theta_R$ and that the above theorem states that any triple $(g, \lambda, \phi)$ extends uniquely to a morphism $\Theta_R \to \cX$ after precomposing $\lambda \co \Theta_K \to \cX$ with the isomorphism $\Theta_K \to \Theta_K$, defined by $t \mapsto \pi^n t$, for $n \gg 0$.
\end{rem}

We will prove this theorem by using the following pushout result.

\begin{lem} \label{L:pushout3}
Let $D = R[t, t/\pi,t/\pi^2,\ldots] \subset K[t]$. The commutative diagram
\begin{equation} \label{D:pushout2}
\begin{split}
\xymatrix{
B_K \bG_m \ar[r] \ar[d]		& B_R \bG_m \ar[d] \\
\Theta_K \ar[r]			& [\Spec(D)/\bG_m]
}
\end{split}
\end{equation}
is cartesian and a pushout diagram in the category of algebraic stacks with affine diagonal.
\end{lem}

\begin{proof}
The proof is identical to the proof of \Cref{L:pushout2} using that $D = R[t]_t \times_{K[t]_t} K[t]$.
\end{proof}

\begin{proof}[Proof of \Cref{T:strange-gluing2}]
By \Cref{L:pushout2}, the triple $(g, \lambda, \phi)$ glues to a morphism $[\Spec(D) / \bG_m] \to \cX$ unique up to unique isomorphism. Writing $D = \bigcup_n R[t/\pi^n]$ and using that $\cX \to S$ is locally of finite presentation we have that for $n \gg 0$ the map $[\Spec(D) / \bG_m] \to \cX$ factors uniquely through $[\Spec(R[t/\pi^n])/\bG_m]$ to yield a map  $[\Spec(R[t/\pi^n])/\bG_m] \to \cX$. 
\end{proof}

\subsection{Semistable reduction in GIT}

\begin{thm} \label{T:valuative-criterion-finite}
Let $\cX$ be a noetherian algebraic stack with affine diagonal.  Assume that either 1) there is a good moduli space $\pi \co \cX \to X$, or 2) $\cX \cong [\Spec(A)/\GL_n]$ and that the adequate moduli space $\pi \co \cX \to X = \Spec(A^{\GL_N})$ is of finite type.  Given a DVR $R$ with fraction field $K$ and a commutative diagram
$$\xymatrix{
\Spec(K) \ar[r] \ar[d]		& \cX \ar[d]^{\pi} \\
\Spec(R) \ar[r]			& X
}$$
there exists an extension of DVRs $R \to R'$ with $K \to K':=\Frac(R')$ finite together with a morphism $h \co \Spec(R') \to \cX$ fitting into a commutative diagram
$$\xymatrix{
\Spec(K') \ar[r] \ar[d]	& \Spec(K) \ar[r] \ar[d]		& \cX \ar[d]^{\pi} \\
\Spec(R') \ar[r]  \ar@{-->}[urr]^h		& \Spec(R) \ar[r]			& X
}$$
such that $h(0) \in | \cX \times_X \Spec(\kappa')|$ is a closed point, where $\kappa'$ is the residue field of $R'$.
\end{thm}

\begin{remark} If $R$ is universally Japanese (e.g. excellent), then it can be arranged that $R \to R'$ is finite.
\end{remark}

\begin{remark} It follows from the valuative criterion for universally closedness (\cite[Thm.~7.3]{lmb}) that there exists a lift $\Spec(R') \to \cX$ with $K \to K'$ a {\it finitely generated} field extension.  Even in the case that $\cX=[\Spec(A)/G] \to \Spec(A^G) = X$ with $G$ linearly reductive and $A$ finitely generated over a field, the above result does not seem to appear in the literature. 
\end{remark}

\begin{proof}[Proof of \Cref{T:valuative-criterion-finite}]
Base changing by $\Spec(R) \to X$, we reduce to the case that $X=\Spec(R)$.  We are given a $K$-point $x_K \in \cX(K)$ and after possibly a finite extension of $K$ (and a corresponding extension of $R$), the unique closed point in $\pi^{-1}(\pi(x_K))$ is represented by a $K$-point $x'_K$.  The closure $\overline{ \{x'_K\}}$ is flat over $\Spec(R)$ and it follows from \cite[IV.17.16.2]{EGA} that after a finite extension of $K$,  the morphism $\overline{ \{x'_K\}} \to \Spec(R)$ has a section $z_R$.
Note that $K$-points $z_K$  and $x'_K$ are isomorphic.  By the Hilbert--Mumford criterion (\Cref{L:isotrivial}), the specialization $x_K \rightsquigarrow z_K$ can be realized by a morphism $\lambda \co \Theta_K \to \cX$ after possibly a further finite extension of $K$.

Restricting $\lambda$ to $0$ yields a map $B_K \bG_m \to \cX$.    Since $\Psi \co \underline{\Map}_R(B_R \bG_m, \cX) \to \cX$, induced by precomposing with $\Spec(R) \to B_R \bG_m$, satisfies the valuative criteria for properness, $B_K \bG_m \to \cX$ extends to a map $g \co B_R \bG_m \to \cX$ such that $z_R$ is isomorphic to $\Psi(g)$.  Since $\cX \to \Spec(R)$ is finitely presented (\cite[Thm.~A.1]{ahr}), \Cref{T:strange-gluing2} implies that the maps $\lambda \co \Theta_K \to \cX$ and $g \co B_R \bG_m \to \cX$ glue  to  a map  $[\Spec(R[t/\pi^n])/\bG_m] \to \cX$.  Restricting this map to $t/\pi^n-\pi$ yields a lift $h \co \Spec(R) \to \cX$ of $x_K$.

Finally, if $h(0) \in |\cX_{\kappa}|$ is not closed, then after possibly a finite extension of the residue field $\kappa$, 
we may use the Hilbert--Mumford criterion to construct a map $\eta \co \Theta_{\kappa} \to \cX_{\kappa}$ which realizes the specialization of $h(0)$ to a closed point.  \Cref{C:strange_gluing} implies that after a finite extension of $R$, $h$ and $\eta$ extends to a map $\oST_{R} \to \cX$.  Restricting this map to $s=t=\sqrt{\pi}$ yields the desired extension $\Spec(R[\sqrt{\pi}])  \to \cX$.
\end{proof}

\subsection{Valuative criteria using DVRs essentially of finite type}

In this section we prove a general lemma showing that over an excellent base it often suffices to check valuative criteria using DVRs essentially of finite type. 

\begin{lem}\label{lem:ess_fintype_criterion}
Let $f\colon X \to Y$ be a quasi compact morphism of algebraic spaces locally of finite type over a {quasi-}excellent base $S$. 
	
	Then $f$ satisfies the valuative criterion for properness for DVRs if and only if it satisfies the lifting criterion for DVRs essentially of finite type over $S$.
\end{lem}
\begin{proof}
Let $R$ be a DVR with fraction field $K$ and $$\xymatrix{ \Spec K \ar[d] \ar[r] & X \ar[d] \\ \Spec R \ar[r]\ar@{-->}[ur] & Y}$$ a lifting problem. As in \cite[TAG03K9]{stacks-project} the claim for algebraic spaces follows from the one for schemes by passing to a cover. We will thus assume that $X,Y$ are schemes.
	
	Denote by $\eta=f(\Spec K) \in X$ the image of $\Spec(K)$. The discrete valuation of $K$ then defines a valuation on $k(\eta)$ and the lifting problem has a solution if and only if it has a solution for the corresponding valuation ring. We can thus assume that $K$ is the residue field of some point of $X$, in particular this field is essentially of finite type over $S$.
	
	The lifting problem has a solution, if and only if there exists an affine open $\Spec B \to X$ containing $\eta$ such that for all $b\in B$ the valuation $b$ in  $K$ is non-negative.
	
	As the closure of $\eta$ is quasi compact we can take finitely many affine open subsets $\Spec B_i \to X$ that cover this closure.
	
	To prove the lemma it thus suffices to show that for finitely many elements $b_i \in k(\eta)$ with $\nu(b_i)=-n_i <0$ there exists a DVR $R^\prime$ essentially of finite type with fraction field $K$ that lies over the image of $\Spec R$ in $Y$ and does not contain the elements $b_i$. 
	
	To construct $R^\prime$ we choose a finitely generated algebra $A$ with fraction field $K$ such that\begin{enumerate}
		\item the morphism $\Spec R \to Y$ factors through $\Spec A$
		\item $A$ contains an element $\pi$ of valuation $1$
		\item $A$ contains the elements $a_i:=\frac{1}{\pi^{n_i-1}b_i}$ which have valuation $1$ and thus lie in the maximal ideal of $R$.
	\end{enumerate} 
	This is possible, by first choosing an affine neighborhood $\Spec B$ of the image of $\Spec R$ in $Y$, and then choosing a subalgebra of $R$ generated by the images of generators of $B$ and then adjoining the finitely many elements of $R$ listed in the second and third point. 
	
	Let $\cp \in \Spec A$ be the image of the closed point of $\Spec R$ so that we get a diagram:
	$$\xymatrix{ \Spec K \ar[d] \ar[rr] & &  X \ar[d] \\ \Spec R \ar[r] &  \Spec A_{\cp} \ar[r] & Y.}$$ 
	
	In the affine blow up algebra $A[\frac{\cp}{\pi}]$ the ideal $A[\frac{\cp}{\pi}] \cp = (\pi)$ becomes principal and contains the elements $a_i$.
	
	{Because A is quasi-excellent, it is Nagata \cite[Tag 07QV]{stacks-project}, and therefore }the normalization of $A[\frac{\cp}{\pi}]_{(\pi)}$ is a DVR, essentially of finite type in which the elements $b_i^{-1}=a_i\pi^{n_i-1}$ have positive valuation.
	
	In particular the lifting problem for this DVR, which is essentially of finite type over $S$ is obstructed.
	
	Similarly, for uniqueness it suffices to use DVRs essentially of finite type, because given a valuation on $K$ and finitely many elements of positive valuation the above construction shows that we can find a DVR essentially of finite in which these elements have positive valuation. This implies that two different lifts to affine neighborhoods can also be realized by a DVR essentially of finite type.
\end{proof}

\bibliographystyle{bibstyle}
\bibliography{refs-moduli}

\end{document}